\documentclass[amscd]{amsart}

\usepackage[latin1]{inputenc}
\usepackage{amsmath}
\usepackage{amssymb}
\usepackage[all]{xy}
\usepackage{mathrsfs}
\usepackage{amsthm}

\usepackage{amscd}

\newtheorem{Thm}{Theorem}
\newtheorem{Cor}[Thm]{Corollary}
\newtheorem{Lem}[Thm]{Lemma}
\newtheorem{Prop}[Thm]{Proposition}

\newtheorem{Claim}[Thm]{Claim} 
\newtheorem{Conj}[Thm]{Conjecture}

\theoremstyle{remark}
\newtheorem{Rem}[Thm]{Remark}

\newtheorem{Ex}[Thm]{Example}

\newcommand{\supp}{\mathop{\mathrm{supp}}}
\newcommand{\rank}{\mathop{\mathrm{rank}}}

\newcommand{\uHom}{{\underline{Hom} }}

\newcommand{\unc}{{\underline{c} }}

\newcommand{\bete}{{\mathrm{b\hat{e}te} }}

\newcommand{\BK}{{\Bbb K}}

\newcommand{\MH}{{\mathcal{MH}}}

\newcommand{\Ql}{{\overline {{\Bbb Q}_l} }}

\newcommand{\cal}{\mathcal}
\newcommand{\Gr}{{{\cal G}{\frak r} }}
\newcommand{\Fl}{{{\cal F}\ell}}

\newcommand{\bu}{\bullet}

\newcommand{\To}{\longrightarrow}
\newcommand{\iso}{{\widetilde \longrightarrow}}

\newcommand{\imbed}{\hookrightarrow}

\newcommand{\ot}{\otimes}

\newcommand{\epf}{\qed} 
\newcommand{\unV}{{\underline V}}

\renewcommand{\P}{{\cal P}}
\newcommand{\PP}{{\cal P}}

\newcommand{\fP}{{^f{\cal P}}}

\newcommand{\DIW}{D_{ IW}^I}
\newcommand{\DIWmon}{D_{IW}^{I^0}}

\newcommand{\Xib}{{\hat{\Xi}}}
\newcommand{\Xis}{\Xi}

\newcommand{\Zb}{{\hat{Z}}}

\newcommand{\Zs}{Z}

\newcommand{\ZZb}{{\hat{\ZZ}}}

\newcommand{\EE}{{\cal E}}
\newcommand{\TT}{{\cal T}}
\newcommand{\TThat}{{\hat{\TT}}}

\renewcommand{\H}{{\mathbb H}}
\newcommand{\Hb}{\overline{{\mathbb H}}}

\newcommand{\Phip}{\Phi_{perf}}
\newcommand{\Phid}{\Phi_{diag}}

\newcommand{\Phib}{\hatt{\Phi}}

\newcommand{\Phipb}{\hatt{\Phi}_{perf}}
\newcommand{\Phips}{\Phi_{perf}}

\newcommand{\Psis}{\Psi'}

\newcommand{\Psib}{{\widehat{\Psi}}}

\newcommand{\fPhi}{{^f{\Phi}}}

\newcommand{\Psit}{\tilde{\Psi}}

\newcommand{\N}{{\cal N}}
\newcommand{\Ntil}{{\tilde{\cal N}}}

\newcommand{\Nt}{{\tilde{\cal N}}}

\renewcommand{\j}{{\frak j}}
\newcommand{\St}{{St}}
\newcommand{\Stt}{{\bf St}}
\newcommand{\suml}{\sum\limits}
\newcommand{\oplusl}{\bigoplus\limits}
\newcommand{\cupl}{\bigcup\limits}

\newcommand{\ho}{{{\frak h}^0}}
\newcommand{\h}{{\frak h}}
\renewcommand{\t}{{\frak t}}
\newcommand{\ft}{{\frak t}}
\newcommand{\fa}{{\frak a}}
\newcommand{\fm}{{\frak m}}
\newcommand{\fz}{{\frak z}}
\newcommand{\fq}{{\frak q}}
\newcommand{\fp}{{\frak p}}

\renewcommand{\b}{{\frak b}}

\newcommand{\bq}{{\bf q}}

\newcommand{\Phicheck}{{\Phi\check{\ }}}
\newcommand{\gal}{\check{\ }}
\newcommand{\M}{{\cal M}}

\newcommand{\That}{{\hat T}}

\newcommand{\g}{{\frak g}}
\newcommand{\Lg}{{ \frak g\check{\ }}}

\newcommand{\Lb}{{ \frak b\check{\ }}}

\newcommand{\LG}{{ G\check{\ }}}
\newcommand{\LB}{{B\check{\ }}}
\newcommand{\LN}{{N\check{\ }}}
\newcommand{\LP}{{P\check{\ }}}
\newcommand{\LQ}{{Q\check{\ }}}
\newcommand{\LU}{{U\check{\ }}}
\newcommand{\LT}{{T\check{\ }}}

\newcommand{\Lt}{{{\frak t}\check{\ }}}

\renewcommand{\L}{{\cal L}}

\renewcommand{\O}{{\cal O}}
\newcommand{\OO}{{\cal O}}
\newcommand{\ZZ}{{\cal Z}}
\newcommand{\Ohat}{{\hat\O}}

\newcommand{\F}{{\cal F}}
\newcommand{\G}{{\cal G}}

\newcommand{\CK}{{\cal K}}

\newcommand{\A}{{\cal A}}
\newcommand{\B}{{\cal B}}
\newcommand{\C}{{\cal C}} \newcommand{\CC}{{\cal C}}
\newcommand{\E}{{\cal E}}
\newcommand{\LL}{{\cal L}}

\newcommand{\cons}{{\overline{\underline{{\Bbb Q}_l}}}}
\newcommand{\J}{{\cal J}}
\newcommand{\JJ}{{\cal J}}
\newcommand{\Jl}{{\cal J}^l}
\newcommand{\Jr}{{\cal J}^r}

\newcommand{\BB}{{\cal B}}
\newcommand{\Z}{{\cal Z}}

\newcommand{\Zet}{{\Bbb Z}}
\newcommand{\Ce}{{\Bbb C}}

\newcommand{\Lotimes}{\overset{\rm L}{\otimes}}
\newcommand{\unO}{{\underline{\cal O}}}
\newcommand{\bG}{{\bf G}}

\newcommand{\bI}{{\bf I}}

\newcommand{\GO}{{\bf G_O}}
\newcommand{\K}{\GO}
\newcommand{\GK}{{\bf G_F}}
\newcommand{\GF}{{\bf G_F}}

\newcommand{\Db}{D^b}

\newcommand{\bbA}{{\mathbb A}}

\newcommand{\Pn}{{\mathbb P}^n}

\newcommand{\Gm}{{\mathbb G}_m}

\newcommand{\Fq}{{\mathbb F}_q}

\newcommand{\Ga}{\Gamma}

\newcommand{\Flt}{\tii \Fl}

\newcommand{\Fltil}{\tii \Fl}

\newcommand{\Phat}{\hat\P}
\newcommand{\Dhat}{\hat D}
\newcommand{\Th}{\hat\TT}

\newcommand{\bs}{\backslash}

\newcommand{\la}{\lambda}

\newcommand{\del}{{\bf{d}}}
\newcommand{\don}{{\mathfrak{d}}}

\newcommand{\hw}{\varpi}

\newcommand{\La}{\Lambda}

\newcommand{\Del}{\Delta}
\newcommand{\nab}{\nabla}

\newcommand{\Ltimes}{\overset{\rm L}{\times}}

\newcommand{\gt}{{\tilde{\frak g}\check{\ }}}

\newcommand{\tii}{\widetilde}

\newcommand{\hatt}{\widehat}

\newcommand{\Cb}{\overline{C}}

\newcommand{\gth}{\widehat{\gt}}
\newcommand{\Lth}{{\widehat{\Lt}}}

\newcommand{\Sth}{\widehat{St}}

\title[Two geometric realizations of affine Hecke algebra]
{On two geometric realizations of an affine Hecke algebra}

\author{Roman Bezrukavnikov}

\address{\small
Department of Mathematics, Massachusetts Institute of Technology, 77 Massachusetts ave.,
Cambridge, MA 02139, USA
}
\address{\small
National Research University Higher School of Economics,
International Laboratory of Representation
Theory and Mathematical Physics,
20 Myasnitskaya st., Moscow 101000, Russia}
\email{
bezrukav@math.mit.edu
}
\begin{document}

\maketitle

\centerline{\em to A. S.-K.}

\bigskip

\begin{abstract}{The article is a contribution to the local theory of geometric Langlands
duality. The main result is a categorification of the isomorphism
between the (extended) affine Hecke algebra associated to 
a reductive 
group $G$ and Grothendieck group of equivariant coherent
sheaves on Steinberg variety of Langlands dual group $\LG$; this isomorphism due to Kazhdan--Lusztig and Ginzburg is a key step in the proof of tamely ramified local Langlands conjectures. 

The paper is a continuation of \cite{AB}, \cite{B}, it relies on technical material developed
in \cite{BY}.

{\em This version contains some post-publication corrections,
namely, typos in and around footnote 5 and mistakes in the proof of \eqref{II} at the end of \S 9.3 and in the statement and the proof of Theorem \ref{suppf}.
The author is deeply grateful to Kostiantyn Tolmachov, Harrison Chen and Ivan Loseu for pointing them out.}}
 \end{abstract}
 
\tableofcontents

\section{Introduction and statement of the result} \subsection{Affine
Hecke algebra and its two categorifications}

 Let $k$ be a field, and let $F=k((t))\supset
O=k[[t]]$ be the field of functions on the punctured formal disc over $k$ and its ring of integers.
Let  $G$
 be a split  reductive 
 linear algebraic group over $k$; let
 $B\subset G$ be a Borel subgroup, and $I\subset G(F)$ be the
 corresponding Iwahori subgroup
 (thus $I$ is the preimage of $B$ under the evaluation
map $G(O)\to G$). 

If $k$ is finite then the group $G(F)$ is a locally compact
topological group, $I$ is its open compact subgroup, and the space
$\Hb$ of $\Ce$-valued finitely supported functions on the two-sided quotient $I\backslash G(F)/I$
carries an algebra structure under convolution; this is the
Iwahori-Matsumoto Hecke algebra. Also $\Hb = \H\otimes_{\Zet[q^{\pm
1}]} \Ce$ where $\H$ is the (extended) affine Hecke algebra and the homomorphism
$\Zet[q^{\pm 1}]\to \Ce$ sends $q$ to $|k|$.

Based on Grothendieck "sheaf-function" correspondence principle, one
can consider  the category of $l$-adic complexes (or perverse sheaves)
on an $\Fq$-scheme (or on its base change to an algebraically closed
field) as the categorical  
counterpart, or categorification, of
 the space of functions on the set of $\Fq$-points of the
scheme; in particular, the space of functions is a quotient of the Grothendieck group of
the category.  
This approach yields a certain derived category of
\'  etale sheaves (respectively, constructible sheaves or $D$-modules) which should be viewed as a categorification of the affine
Hecke algebra $\H$. 

On the other hand, as was discovered by
Kazhdan and Lusztig (and independently by Ginzburg), 
the affine Hecke
algebra can be realized as the Grothendieck group of equivariant 
coherent sheaves
on the Steinberg variety of the Langlands dual group, 
thus the corresponding derived
category of coherent sheaves provides another categorification of $\H$.

The goal of the present paper is to construct an equivalence between the two
triangulated categories which categorify the affine Hecke algebra.
  A step in this direction has been made in the previous works \cite{AB},
\cite{B}, where a geometric theory of the anti-spherical
(Whittaker) module over $\H$ was developed;\footnote{In {\em loc. cit.} 
the group $G$ is assumed to be simple. However, its arguments apply also to the
case of a general reductive group $G$.}
 in the present paper we
extend this analysis to the affine Hecke algebra itself.

The possibility to realize the affine Hecke algebra $\H$ and the
"anti-spherical" module over it as Grothendieck groups of
(equivariant) coherent sheaves on varieties appearing in the Springer theory for $\LG$ plays
a key role in Kazhdan--Lusztig's proof of classification of irreducible
representations of $\Hb$, which constitutes a particular case of local Langlands
conjecture, see \cite{KL} and exposition in \cite{CG}.\footnote{In fact, some of the key ideas of this theory already appeared in an earlier work of 
Lusztig   \cite{L_pre1}, \cite{L_pre2} where certain modules over the affine Hecke algebra
were realized via $K$-groups of Springer fibers; also the relation between  the $q$-deformation
of the $K$-group and dilation equivariance was described in {\em loc. cit.}}
 Thus one
may hope that the categorification of these realizations
proposed here can contribute to the geometric Langlands program.
In fact, since the result of the paper was announced, it has been 
applied and generalized by several authors working in that area, see \cite{FG}, \cite{BZN}, \cite{BLin}. 
Let us point out that existence of (some variant of) such a
categorification was proposed as a conjecture by V.~Ginzburg, see
Introduction to \cite{CG}.

\subsection{Statement of the result}
Let us now describe our result in more detail.

\subsubsection{Categories of $l$-adic sheaves (the "Galois side")} Recall the
well known group schemes $\K\supset \bI$ over $k$ (of infinite type)
such that  $\K(k)=G(O)$, $\bI(k)=I$, and a group ind-scheme $\GF$
with $\GF(k)=G(F)$. We let  $\bI^0$ be  the pro-unipotent
radical of $\bI$; if $k={\mathbb F}_{p^a}$ then $I^0=\bI^0(k)$ is the pro-$p$ radical of $I$. 
We also have the quotient ind-varieties:
the affine Grassmanian $\Gr$,   the
affine flag variety $\Fl=\GF/\bI$ and the extended affine flag variety
$\tii \Fl=\GF/\bI^0$, see
e.g. \cite{KGB}, Appendix, \S A.5.
 Thus $\Gr$, $\Fl$, $\tii \Fl$ are direct limits of
finite dimensional varieties with transition maps being closed embeddings,
in the case of $\Gr$ and $\Fl$ all the finite dimensional varieties in the direct system are projective.
We have
 $\Gr(k)=G(F)/G(O)$, $\Fl(k)=G(F)/I$, $\tii \Fl(k)=G(F)/I^0$.


{\em From now on we assume that the base field $k$ is algebraically closed.}

  Let  $D(\tii \Fl)$, $D(\Fl)$, $D(\Gr)$ be the constructible
 derived categories of $l$-adic sheaves ($l\ne
char(k)$; see \cite[\S 1.1.2]{Weil2}, 
\cite[\S 2.2.14--2.2.18]{BBD}; and  \cite[\S A.2]{KGB},  for a
(straightforward) generalization of the definition of an $l$-adic complex to a certain class of
ind-schemes)
 on the respective spaces.\footnote{If $char(k)=0$ we can also consider the corresponding derived
 categories of $D$-modules and if $k=\Ce$ we can work with the derived category of constructible
 sheaves in the classical topology. 
  All our results, excepts for some statements in \S \ref{sec_fur} which explicitly involve a Frobenius action, hold in these settings, the proofs are identical.}

 The protagonists of this paper are as follows.
  Let
$D_{II}=D_I(\Fl)$  be the $\bI$-equivariant derived category of
$l$-adic sheaves on $\Fl$; $D_{I^0 I}=D_{I^0}(\Fl)$ be the $\bI^0$-equivariant
derived category of $l$-adic sheaves on $\Fl$, and let $D_{I^0I^0}$ be the full subcategory
in the $\bI^0$-equivariant derived category of $\tii \Fl$ consisting of complexes whose cohomology
is monodromic (or weakly equivariant, see \cite{Ver})
 with respect to the right $T=\bI/\bI^0$ action with unipotent monodromy.

The categories $D_{II}$ and $D_{I^0I^0}$ are equipped with an associative product operation
provided by convolution; $D_{II}$ is unital while $D_{I^0I^0}$ lacks the unit object.\footnote{Notice
 \label{footconv}
 that convolution with an object of $D_{I^0I^0}$ involves direct image under a non-proper
morphism, thus convolution could be defined in two different ways, using either direct image or direct image with compact support; we use the version with the ordinary direct image.
However, the convolution diagram involved in the definition of convolution in $D_{I^0I^0}$
is a composition of a $T$ bundle and a proper morphism, while the sheaves to which we apply
the direct image are $T$-monodromic. The direct image under projection to the base
of a $T$-monodromic sheaf on a principal $T$-bundle
can be expressed as derived invariants of monodromy, cf. Lemma \ref{toruseqDG}, 
while direct image with proper support
is a Verdier dual operation. Since 
derived invariants with respect to a symmetric algebra action is a self-dual operation, up to homological shift, we see that the two definitions produce equivalent
monoidal categories.}
We have commuting actions of $D_{I^0 I^0}$ and $D_{II}$ on $D_{I^0I}$ by left and right
convolution respectively. The convolution operation will be denoted by $*$.

Let $\P_{II}\subset D_{II}$, $\P_{I^0I}\subset D_{I^0I}$, $\P_{I^0I^0}\subset D_{I^0I^0}$
be the subcategories of perverse sheaves.
A standard argument (see e.g. \cite[Proposition 1.5]{BBM} for the first
equivalence, the second one  follows by a similar argument using e.g.
\cite[Corollary A.4.7]{BY})  
shows that
\begin{equation}\label{DisDP}
\begin{array}{ll}
D^b(\P_{I^0I})\cong D_{I^0I},\\
D^b(\P_{I^0I^0})\cong D_{I^0I^0},
\end{array}
\end{equation}
while the natural functor $D^b(\P_{II}) \to D_{II}$ is not an equivalence.

\subsubsection{The dual side} \label{dual_side} Let $\LG$ be  the
Langlands dual group over the field $\Ql$.
The goal of the paper is to provide a description for the above categories in terms of
$\LG$. To formulate the answer we need to recall the following construction.

Let $X\to Y$, $X'\to Y$ be  morphisms of algebraic varieties. 
We will
assume that $X,X', Y$ are varieties over a field $k$, $Y$ is smooth and morphisms $X\to Y$,
$X'\to Y$ are proper.

 One can consider the
{\em derived fiber product} $X\Ltimes _Y X'$ which is a {\em differential graded scheme}
(DG-scheme for short), and the triangulated category $DGCoh(X\Ltimes_Y X')$. 

If $Tor_i^{\O_Y}(\O_X, \O_{X'})=0$ for $i>0$ then the derived fiber product
reduces to the ordinary fiber product and $DGCoh(X\Ltimes _Y X')=D^b(Coh(X\times _Y X'))$.

The triangulated category $DGCoh(X\Ltimes _Y X)$ has a natural  monoidal structure
provided by convolution. The category $D^b(Coh(X))$ is naturally a module
category for the monoidal category $DGCoh(X\Ltimes _Y X)$. [For example, when $X$ is a finite set and $Y$ is
a point the induced structures on the Grothendieck group amount to matrix
multiplication and the action of $n\times n$ matrices on $n$-vectors respectively].
More generally, the category $DGCoh(X\Ltimes _Y X')$ has two commuting actions:
 the action of $DGCoh(X\Ltimes _Y X)$ on the left and an action of $DGCoh(X'\Ltimes _Y X')$
on the right. 

Given an action of an affine algebraic group $H$ on $X,\, X',\, Y$ compatible with the maps,
one gets equivariant versions of the above statements.

We will apply this in the following situation. We let $Y=\Lg$ be the Lie algebra of $\LG$,
$X=\gt=\{ (\b, x)\ |\ \b\in \BB, \, x\in \b\}$, $X'= \Nt=\{ (\b, x)\ |\ \b\in \BB, \, x\in rad(\b)\}$,
where $\BB$ is the flag variety for $\LG$ parametrizing Borel subalgebras in $\Lg$.

A standard complete intersection argument shows that 
$Tor_{>0}^{\O_Y}(\O_X, \O_{X'})=0$ for  $X=\gt$, $Y=\Lg$ and $X'=\gt$ or $X'=\Nt$,
thus the corresponding derived fiber products coincide with the usual fiber product of schemes. 
However, it fails for $X=X'=\Nt$, $Y=\Lg$, so the derived fiber product $\Nt\Ltimes _\Lg \Nt$
is essentially different from $\Nt\times _\Lg \Nt$.

We set $St=\gt\times _\Lg \gt$, $St'=\gt\times _\Lg \Nt$.

\subsubsection{Statement of the result}
We now formulate the main result of the paper.

For an algebraic variety  $X$ and a closed subset $Z\subset X$ we will let $Coh_Z(X)$
denote the full subcategory in $Coh(X)$ consisting of sheaves set-theoretically supported on $Z$.
For a map $f:X\to Y$ and a closed subset $Z\subset Y$ we will abbreviate  $Coh_{f^{-1}(Z)}(X)$
to $Coh_Z(X)$.

\begin{Thm}\label{mainthm}
There exist natural equivalences of categories:
\begin{equation}\label{I0I0}
\Phi_{I^0I^0}:D_{I^0I^0}\cong D^b(Coh^\LG_\N(St)),  
\end{equation}
\begin{equation}\label{I0I}
\Phi_{I^0I}:D_{I^0I}\cong D^b(Coh^\LG (St')), 
\end{equation}
\begin{equation}\label{II}
\Phi_{II}:D_{II}\cong DGCoh^\LG (\Nt\Ltimes_\Lg \Nt).
\end{equation}
Equivalences \eqref{I0I0} and \eqref{II} are compatible with the convolution product,
while \eqref{I0I} is compatible with the action of the categories from \eqref{I0I0} and \eqref{II}.
\end{Thm}

\subsection{The action on the Iwahori-Whittaker category}

It was pointed out above that the monoidal category of $DG$ coherent (equivariant)
sheaves on a fiber product $X\Ltimes _Y X$ admits a natural action on the derived category
of (equivariant) coherent sheaves on $X$. In particular, monoidal category
$DGCoh^\LG(\Nt\Ltimes _{\Lg} \Nt)$ acts on $D^b(Coh(\Nt))$, while $D^b(Coh^\LG(St))$
acts on $D^b(Coh^\LG(\gt))$.  

To describe the corresponding structures on the loop group side, recall the category
of {\em Iwahori-Whittaker} sheaves.  The quotient of $\bI^0$ by its commutant
is the sum of copies of the additive group indexed by vertices of the affine Dynkin graph.
Fix an additive character $\psi$ of  $\bI^0$ which is trivial on the summand of $\bI^0/(\bI^0)'$
corresponding to the affine root(s) and is non zero on the other summands.
We denote by $\DIW$ the $\bI$-equivariant derived category of $l$-adic sheaves
on the principal homogeneous space $\GF/(\bI^0)'$ which satisfies the $\psi$-equivariance
condition with respect to the right action of $\bI^0/(\bI^0)'$, see \cite{AB}.
(Conventions here differ from those of \cite{AB} by switching the roles
of left and right multiplication.) We let $\DIWmon$ denote the category
of $\bI$ monodromic sheaves with unipotent monodromy on $\GF/(\bI^0)'$ which are
 $\psi$-equivariant with respect to the right action of $\bI^0/(\bI^0)'$,  this is a particular
 case of the category considered in \cite{BY} (again, one needs to switch left with right to get from the present setting to that of \cite{BY}).

The categories $D_{I^0I^0}$, $D_{II}$ act on $D_{IW}^{I^0}$, $D_{IW}^I$ respectively
by convolution.


\begin{Thm}\label{mainthm2}
There exist equivalences of categories 
\begin{equation}\label{PhiIW}
\Phi_{IW}^I:D^b(Coh^\LG(\Nt))\iso
\DIW,
\end{equation}
\begin{equation}\label{PhiIWmon}
 \Phi_{IW}^{I^0}:D^b(Coh^\LG_\Nt(\gt))\iso \DIWmon,
 \end{equation} satisfying the following compatibilities:
the equivalence $\Phi_{IW}^{I^0}$ is compatible with the action
of $D^bCoh^\LG_{\N}(St)$ coming from 
the
action of $D_{I^0I^0}$ on $\DIWmon$ and equivalence \eqref{I0I0}. 

The equivalence $\Phi_{IW}^{I}$ is compatible with the action
of $DGCoh^\LG(\Nt\Ltimes _\Lg \Nt)$ coming from 
the
action of $D_{II}$ on $\DIW$ and equivalence \eqref{II}. 
\end{Thm}

Another useful compatibility between the equivalences in Theorems \ref{mainthm}
and \ref{mainthm2} is stated at the end of section \ref{monstr}.

A variant of equivalence \eqref{PhiIW} has been established in \cite{AB}, and 
\eqref{PhiIWmon} can obtained by a similar argument, see below. More precisely,
in \cite{AB} a functor $F:D^b(Coh^\LG(\Nt))\to D_{II}$ is constructed, below we construct
its "monodromic" counterpart $\Phid:D^b(Coh^\LG(\gth))\to \Dhat$,
where $\Dhat$ is a certain "pro"-completion of $D_{I^0I^0}$ and $\gth$ (a version of) the formal
neighborhood of $\Nt$ in $\gt$ (see \S \ref{furnot} for a precise definition). One can consider
the composition of $F$ with either left  or right Whittaker averaging, both compositions
turn out to be equivalences, the proofs of these two facts are parallel. In \cite{AB} we worked
with left Whittaker averaging, while here we work with the right one (this allows us to work
with modules over the monoidal category $D_{II}$ rather than modules over its opposite).






\subsection{Description of the strategy: the  Hecke algebra perspective}\label{Heckper}
Some of the constructions exploited here are sheaf-theoretic analogs of known results in the
theory of affine Hecke algebras.

 Recall that $\H$ has a standard basis $t_w$ indexed by elements $w$
 in the extended affine Weyl group $W$. 
 
Let $\Lambda$ be the coweight lattice of $G$ 
and $\Lambda^+\subset \Lambda$
be the set of dominant weights. There exists a unique collection of elements $\theta_\lambda
\in \H$, $\lambda \in \Lambda$, such that $\theta_\la\theta_\mu=\theta_{\la+\mu}$
for all $\la,\, \mu\in \La$ and $\theta_\la=T_\la$ for $\la\in \Lambda^+$. 
The categorification of the elements $\theta_\la$ are the so-called {\em Wakimoto
sheaves}, see \cite{AB} and section \ref{Wakimo} below. 

The elements   $\theta_\la$ span
  a commutative
subalgebra $A \subset \H$ which contains the center $Z(\H)$
of the affine Hecke algebra. Categorification of the center is provided
by the work of Gaitsgory \cite{KGB}. Categorification of the formula
expressing central elements as  linear combinations of $\theta_\la$
is the fact that central sheaves of \cite{KGB} admit a filtration whose
associated graded is a sum of Wakimoto sheaves, see \cite{AB} and section
\ref{moncen} below.  This filtration plays a key role in our construction, see \cite{AB} and
section \ref{monfundiag},
yielding a categorification of:

\begin{equation}\label{Ato}
A\cong K^0(Coh^\LG(\gt))\overset{\delta_*}{\To} K^0(Coh^\LG(St)),
\end{equation}
where $\delta:\gt\to St$ is the diagonal embedding. 

Another ingredient important to us is
the {\em $q$-analog of the Schur anti-symmetrizer, or anti-spherical projector} $\xi=\suml_{w\in
W_f} (-1)^{\ell(w)} t_w$. Its relevance to representation theory of $p$-adic
groups comes from the fact that the left ideal $\Hb \xi$ is canonically
isomorphic to  $I$ invariants in the space of {\em Whittaker functions} on $G(F)$, while
its relation to canonical basis in the affine Hecke algebra, thus to perverse sheaves on $\Fl$
goes back to \cite{LAst}.

The categorical counterpart of $\xi$ is the
maximal projective object in the
category of perverse sheaves on $G/B\cong \GO/\bI$ equivariant with respect
to $\bI^0$, it is discussed in section \ref{sectionXi}.
Under the equivalence with the coherent sheaves category that object
corresponds to the structure sheaf of Steinberg variety.

Let $\H_{perf}\subset \H$ be the 
two-sided ideal generated by $\xi$.
 The full subcategory
$D_{perf} ^\LG(St)\subset D^b(Coh^\LG(St))$ of perfect complexes can be considered as
a categorification of $\H_{perf}$. 
Furthermore, it is easy to see that 
$\H_{perf}$  is freely generated by $\xi$ as a module over $A\otimes _Z A$. This allows one to deduce
an equivalence between the two categorifications of $\H_{perf}$ from the categorification
of \eqref{Ato}. 
The  subcategory $D_{perf}^\LG(St)$ is dense in 
$D^b(Coh^\LG(St))$ in an appropriate sense, which allows to extend the equivalence
from the subcategory to the whole category.

\subsection{Acknowledgements} 
The initial  ideas of this paper were conceived during the Princeton IAS special year 1998/99 led by G.~Lusztig,  the first stages were carried out as a joint project with S.~Arkhipov \cite{AB}. 
Since then the material was discussed with many people; 
the outcome was particularly influenced by the input from A.~Beilinson and V.~Drinfeld,
many conversations with M.~Finkelberg, D.~Gaitsgory, D.~Kazhdan and I.~Mirkovi\' c and others were important
for keeping the project alive.
More recently I have benefited from the expert advice of D.~Arinkin,
L.~Positelskii and Z.~Yun. 
I would like to express my gratitude to these mathematicians.
I am also much indebted to the referee for a careful reading resulting in many corrections  
 and to G.~Lusztig 
for helpful comments on the text.

The author was partially supported by an NSF grant and a Simons Foundation Fellowship.

\medskip

{\small
In this text we follow the original plan conceived more than a decade ago and treat
the issues of homological algebra by ad hoc methods, using explicit DG models for triangulated
categories of constructible sheaves based on generalized tilting sheaves. 
While the properties of tilting sheaves established in the course of the argument
are (in the author's opinion) of an independent interest, it is likely that recent advances
in homotopy algebra can be used to develop an alternative approach.}

\section{Outline of the argument}

\subsection{Further notations and conventions}  \label{furnot}
We let 
 $B\supset N$ be a Borel subgroup in $G$ and its unipotent radical, and $\LN\subset \LB$ be similar
 subgroups in $\LG$; we assume that $B$ is the image of $\bI$ under the evaluation map
 $\GO\to G$.
 
Let $\Lambda$ be the coweight
lattice of $G$, i.e.  the coweight lattice of the 
abstract Cartan group of $G$;
thus $\Lambda$ is
 identified with the weight lattice of (the abstract Cartan group of) $\LG$.  
We let $W_f$ denote the Weyl group and $W=W_f\ltimes \Lambda$ the 
extended affine Weyl group;
 $\ell:W\to \Zet_{\geq 0}$ is the length function, and $\Lambda^+\subset \Lambda$
 is the set of dominant coweights, $w_0\in W_f$ is the longest element.
 
 We let $W^f\subset W$ be the subset of minimal length representatives
 of right cosets $W/W_f$. Notice that  
 $\La^+\subset W^f$.
 
 We let $\B=\LG/\LB$ be the flag variety. 
The set of isomorphism classes of $\LG$-equivariant line bundles on $\B$ is 
identified with $\La$; for $\la\in \La$
 we let  $\OO_\B(\la)$ be the corresponding line bundle. Recall that $\O(\la)$ is semi-ample
 iff $\la\in \La^+$. For $\la\in \La^+$ we let $V_\la$ denote the corresponding irreducible $G$-module,
thus  $V_\la=\Gamma(\B, \O(\la))$.

 The ind-schemes $\Fl=\GK/\bI$, $\Fltil=\GK/\bI^0$   and the categories $D_{II}\supset \P_{II}$, $D_{I^0I}\supset \P_{I^0I^0}$, $D_{I^0I^0}\supset \P_{I^0I}$
 were introduced above. We abbreviate $\P=\P_{I^0I}$ and let $\Phat$, $\Dhat$
 be the pro-completions of $D_{I^0I^0}$ and $\P_{I^0I^0}$ respectively, see section
 \ref{monpro}. 

Let $\pi:\Fltil \to \Fl$ be the projection.

  The $\bI$ orbits 
 on $\Fl$ are indexed by $W$, for $w\in W$ we let 
 $j_w:\Fl_w\to \Fl$ be the embedding of the corresponding orbit.
 We have $\dim(\Fl_w)=\ell(w)$.

  We have standard objects $j_{w!}:=j_{w!}(\cons[\ell(w)])$  and costandard object
  $j_{w*}=j_{w*}(\cons[\ell(w)] ) $ in $\P$. Their counterparts in $\Dhat$ are the free monodromic
  (co)standard objects $\nab_w$, $\Del_w$, see sections \ref{genonmo}, \ref{moreonmo}.
 
 We also consider the Iwahori-Whittaker categories $D_{IW}^I\supset \P_{IW}^I$, $D_{IW}^{I^0}
 \supset \P_{IW}^{I^0}$, the pro-completions $\Dhat_{IW}$, $\Phat_{IW}$ of, respectively,
 $D_{IW}^{I^0}$, $\P_{IW}^{I^0}$,
  (co)standard objects  
    $j_{w!}^{IW}, j_{w*}^{IW}\in \P_{IW}^I$
    and free monodromic (co)standard object $\Del_{w}^{IW}$, $\nab_{w}^{IW}\in \Phat_{IW}$,
    $w\in W^f$ (see \S \ref{monpro} for further details).  
 
 Recall that $St=\gt \times_\Lg \gt$, let $p_{Spr,1}:St\to \gt$, $p_{Spr,2}:St\to \gt$ be
 the two projections. Also $St'=\gt \times_\Lg\Nt$ with two projections $p_{Spr,1}':St'\to
 \gt$, $p_{Spr,2}':St'\to \Nt$.
 Let  $\gth=\gt \times _{\Lg} \hatt{\Lg}$, 
 $\Sth=St\times _{\Lg}  \hatt{\Lg}$, where $\hatt{\Lg}$ is the 
 spectrum\footnote{Alternatively \label{footform}
 we could work with completion defined as a formal scheme,
 the resulting category of coherent sheaves would be equivalent. In more detail,
 by \cite[Th\' eoreme 5.4.1]{EGA}
  the scheme $\gth$, $\Sth$ is the inductive limit in the category
of schemes over $\Lg$ of nilpotent thickenings of $\Nt$ in $\gt$ (respectively, $St'$ in $St$). 
By \cite[Th\' eoreme 5.1.4(1)]{EGA} the category of coherent sheaves on 
 $\gth$  is
equivalent to the category of coherent sheaves on the formal scheme completion of $\Nt$ in 
  $\gt$ and similarly for $\Sth$, this readily extends to the category of $\LG$-equivariant sheaves; cf. the discussion
  at the end of Introduction to \cite{BM}.}
  of the completion of the ring of functions $\O(\Lg)$ at the ideal of the subscheme $\N$.


For an algebraic group $H$ acting on an (ind)-scheme $X$ we let $D_H(X)$
denote the equivariant derived category of $H$-equivariant constructible sheaves
on $X$ and let
$Perv_H(X)\subset D_H(X)$ be the subcategory of perverse sheaves.
Given a subgroup $K\subset H$ we have the functor of restricting the equivariance
$Res^H_K:D_H(X)\to D_K(H)$ and  the left adjoint functor
 $Av_K^H:D_K(X)\to D_H(X)$ (the latter can be thought of as the $!$-direct image  for
 the morphism of stacks $X/K\to X/H$). 
  In particular, we have a functor $Av_{\bI^0}^\bI:D_{I^0I}\to D_{II}$ (to unburden
  typography we will write $Av_{I^0}^I$).
  
  In order to introduce a similar functor involving  Iwahori-Whittaker sheaves we fix an 
  Iwahori subgroup $\bI_-\subset \GO$ which is opposite to (in general position with)
 the subgroup $\bI$.  We also fix a nondegenerate additive character $\psi_-$ of $\bI_-^0$.
  The pair $(\bI_-^0,\psi_-)$ is conjugate to $(\bI^0,\psi)$ by an element in $G(F)$
  which is unique up to right multiplication by an element in $I^0$. 
  Thus the categories $D^I_{IW}$, 
  $D^{I^0}_{IW}$ are canonically equivalent to the derived categories 
  $D^I_{IW_-}$, $D^{I^0}_{IW_-}$ of right $(\bI_-,\psi_-)$-equivariant sheaves.  
  We define the functors $Av^{I^0}_{IW}: D^{I^0}_{IW_-}\cong D^{I^0}_{IW}\to
  D_{I^0I^0} $,
  $Av^{IW}: D_{I^0I^0}\to D_{IW}^{I^0}$
  by setting $Av^{I^0}_{IW}= Av_{\bI^0\cap\bI^0_-}^{\bI^0}$,
  $Av^{IW}=Av_{\bI^0\cap\bI^0_-}^{\bI^0_-,\psi_-}$.
  Here we used that $\psi |_{\bI^0\cap \bI^0_-}$ is trivial; the restriction of equivariance
  functor is omitted from notation, and $Av_{\bI^0\cap\bI^0_-}^{\bI^0_-,\psi_-}$
  is the left adjoint to the restriction of equivariance functor from $D^{I^0}_{IW_-}$
  to the corresponding $\bI^0\cap \bI^0_-$-equivariant category. The result
  of \cite{BeBrM} implies that we get the same functor $Av^{IW}$ if we replace
  $Av_{\bI^0\cap\bI^0_-}^{\bI^0_-,\psi_-}$ by the corresponding right adjoint
  to the restriction of equivariance functor: the Whittaker averaging functor is clean.
  We also have a similarly defined functor on $I$-equivariant categories:
  $^I Av_{IW}: D_{II^0}\to D^I_{IW}$.
      
  Notice that the definition of $D^{I^0}_{IW}$
  involves the left action of $I^0$, while $Av^{I^0}_{IW}$, $Av^{IW}$
have to do with the right action; when the action used may not be clear from the context
we use notation $Av_{I^0}^{left}$, $Av_{I^0}^{right}$ to distinguish between the two.


\subsection{Idea of the argument: structural aspects}\label{strasp}
The functor from the coherent category to the constructible one stems from certain natural structures on the constructible category. To describe the mechanism of obtaining such a functor from the additional structures on the target category it  is convenient
to use the concept of a triangulated category $\CC$ {\em over} a stack $X$.

\subsubsection{Linear structure over a stack}\label{linover}
We refer to  \cite{Gaover} for the notion of an abelian category over an algebraic stack, and to \cite{Gaoverder} 
 for a generalization to triangulated (or rather homotopy theoretic) context.
For our present purposes it suffices to use the following simplified version of this concept.
Let $S$ be an algebraic stack and $\CC$ a triangulated category (in all our example $S=X/G$
where $X$ is a quasi-projective algebraic variety and $G$ is a reductive algebraic group).
The subcategory of perfect complexes $D_{perf}(S)\subset D^b(Coh(S))$ is a triangulated
tensor category under the usual tensor product
of coherent sheaves. By an {\em $S$-linear structure} on $\CC$ we will mean an
action of the tensor category $D_{perf}(S)$ on $\CC$ compatible with the triangulated
structure. 

We now list basic classes of examples of such a structure to be used below.

\begin{enumerate}
\item If $S=Spec(R)$ is an affine scheme, then for an $R$-linear abelian category
$\A$ the triangulated category $D^b(\A)$ acquires a natural $S$-linear structure.

\item Let $S=pt/H$ where $H$ is a linear algebraic group. If  an abelian category
$\A$ is a module category for the tensor category $Rep(H)$ of algebraic (finite dimensional)
representations acting by exact functors, then
$D^b(\A)$ is an $S$-linear triangulated category.

\item  Combining the first two examples, 
assume now that $S=Spec(R)/H$ is a quotient of an affine scheme by a linear algebraic
group action. Let $\A$ be an abelian category which is a module category for $Rep(H)$ acting on $\A$ by exact functors. Then we can define a new (in general not abelian)
 {\em "deequivariantized"} category $\A_{deeq}$ by setting 
$Ob(\A_{deeq})=Ob(\A)$, $Hom_{deeq}(X,Y)=
Hom_{Ind(\A)}(X,\unO_H(Y))$ where $Ind(\A)$ stands for the category of Ind-objects in $\A$ and 
 $\unO_H\in Ind(Rep(H))$ denotes the space of regular functions on $H$ with $H$ acting by left translations,
 see  section \ref{deeqsubs}  for further details.
 
 Then $\A_{deeq}$ is a category enriched over the category of algebraic (not necessarily finite dimensional) representations of $H$.  
 Then  an $R$-linear structure on 
 $\A_{deeq}$ which is compatible with the $H$-action induces an $S$-linear structure
 on $D^b(\A)$. To see 
  this observe that $Hom_{Coh^H(Spec(R))}(V\otimes \O,V'\otimes \O)
 = (V'\otimes V^*\otimes R)^H$, thus an equivariant $R$-linear structure on 
$\A_{deeq}$ induces an action of the tensor category $Coh_{fr}(S)$ on $\A$ by exact functors;
 here $Coh_{fr}(S)
\subset Coh(S)=Coh^H(Spec(R))$ is the full subcategory consisting of objects $V\otimes
\O_{Spec(R)}$, $V\in Rep(H)$. 
Since $D_{perf}(S)$ is the Karoubi (idempotent) completion
 of the homotopy category of finite complexes  $Ho(Coh_{fr}(S))$, the action
 of $Coh_{fr}(S)$ on $\A$ induces an $S$-linear structure on $D^b(\A)$ (notice that
 $D^b(\A)$ is necessarily Karoubian).
 
\item Suppose we are given an open embedding of algebraic stacks $S\imbed S'$ and a category
$\CC$ with an $S'$-linear structure; assume for simplicity that $S'$ is a quotient
of a quasi-projective variety over a field of characteristic zero by an action of a reductive group
and $S$ comes from an invariant open subvariety therein. By results of \cite{TT} 
(see also  \cite[\S 2.1.4, esp. proof of Lemma 2.6]{Nee}\footnote{The result is only claimed in {\em loc. cit.} 
for a subscheme in a scheme
but the case of stacks of the described type follows by the same argument.}) we have
\begin{equation}\label{Idem}
D_{perf}(S)\cong Idem (D_{perf}(S')/D_{perf}(S')_{\partial S'})
\end{equation}
where $Idem$ denotes the Karoubi (idempotent)
completion and 
$D_{perf}(S')_{\partial S'}$ is the full subcategory of perfect complexes on $S'$ whose restriction
to $S$ vanishes. Thus if $\CC$ is a Karoubian (idempotent complete) category, then an $S'$-linear
structure on $\CC$ such that $D_{perf}(S')_{\partial S'}$ acts by zero induces an $S$-linear
structure on $\CC$.

\item One can use a variant of Serre's description of the category of coherent
sheaves on a projective variety as a quotient of the category of graded modules over
the homogeneous coordinate ring to devise a procedure for constructing an $S$-linear
structure for more general stacks $S$.

Suppose that $S=X/H$ where $X$ is a quasi-projective variety with an action 
of an affine algebraic group $H$. Assume that a linearization of the action, i.e.
a linear action of $H$ on the linear space $ {\mathbb A}^{n+1}$
together with an equivariant locally closed
embedding $X\to \Pn$ is fixed. Let $C\subset {\mathbb A}^{n+1}$ be the cone over the closure $\overline{X}$ of $X$
in $\Pn$. Then $C$ is an affine variety acted upon by $H\times \Gm$ and we have
an open embedding  $S\to S'=C/(H\times \Gm)$. 
Using \eqref{Idem} we see that an $S$-linear structure on $\CC=D^b(\A)$
can be constructed by providing $\A$ with a $Rep(H\times \Gm)$ action by exact functors,
 introducing an $R$-linear structure on $\A_{deeq}$ where $R=\O(C)$ is the homogeneous coordinate
 ring of the projective variety $\overline{X}$, and verifying that the resulting
 $S'$-linear structure sends $D_{perf}(S')_{\partial S'}$ to zero.
  \end{enumerate}

\begin{Rem}
 Most of the statements in the main Theorem of the paper assert an equivalence between (a subcategory of) $D^b(Coh(S))$ for an algebraic stack $S$ and
$D^b(\A)$ for an abelian category $\A$ (with the exception of \eqref{II} which involves
coherent sheaves on a DG-stack and an equivariant derived category of constructible sheaves).

We first construct the $S$-linear structure on $\CC=D^b(\A)$ and then consider the action
on a particular object of $\CC$ to get an equivalence. The construction of the action 
almost follows the pattern of example (5). The difference is as follows. We have $S=X/\LG$
where $X$ admits an affine equivariant map to $\BB^2$. Though $\BB^2$ is a projective variety
there is no preferred choice of an equivariant projective embedding, so to keep things
more canonical we work with the "multi-homogeneous" coordinate ring and consider open 
embeddings of our stacks into $Y/(\LG \times \LT^2)$ for an appropriate affine variety
$Y$. A more essential difference is that while $Rep(\LG)$ acts by exact functors on our abelian
category $\A$, the action of $Rep(\LT^2)$ is only defined on the triangulated category
$\CC$, it is not compatible with the natural $t$-structure on $\CC=D^b(\A)$. 

An additional argument based on properties of tilting modules is needed to deal with this issue (see subsection \ref{extecompl}).
\end{Rem}

\subsubsection{The list of structures}
We concentrate on the equivalence \eqref{I0I0}, the equivalence \eqref{I0I} is similar,
and \eqref{II}  will be deduced formally from \eqref{I0I0}.

Consider the following sequence of maps
$$St/\LG \rightrightarrows \gt/\LG \to \Lg/\LG \to pt/\LG.$$

Moving from right to left in this sequence, we successively equip $\Dhat$ with the linear structure 
 for the corresponding stack.

The $pt/\LG$-linear structure comes from an action of the tensor category $Rep(\LG)$
on the abelian category $\P_{I^0I}$. Such an action was defined in \cite{KGB} where the {\em central sheaves}
categorifying the canonical basis in
 the center of the affine Hecke algebra were constructed; an extension 
of the action to $\P_{I^0I^0}$ is sketched in section \ref{moncen} below.

By a version of the Tannakian formalism, 
lifting an action of the tensor category $Rep(\LG)$ to a $\Lg/\LG$-linear structure amounts to equipping the $Rep(\LG)$ action with a tensor endomorphism. Such an endomorphism comes
from the {\em logarithm of monodromy} acting on central sheaves: recall that the central sheaves are constructed by nearby cycles which carry a monodromy automorphism.

We now discuss the two structures of a stack over $\gt/\LG$. 
The starting point here is the familiar observation that for a representation
$V$ of $\LG$ the trivial vector bundle $V\otimes \O_\BB$ 
with fiber $V$ on the flag variety $\BB=\LG/\LB$ carries a canonical filtration
whose associated graded is a sum of line bundles. This filtration can be lifted
 to a similar filtration for $V\otimes \O_{\gt}$. 
Under our equivalences this filtration corresponds to a filtration on (monodromic) 
central sheaves by (monodromic) {\em Wakimoto sheaves} (the non-monodromic version
was presented in \cite{AB}, and the monodromic generalization is presented below
in \S \ref{Wakimo}). 
It turns out that the filtration defines a monoidal functor $D^b(Coh^\LG(\gt))\to \Dhat$.
 We then get two commuting actions of $D^b(Coh^\LG(\gt))$ on $D_{I^0I^0}$ from the left
 and the right 
  action of the monoidal category $\Dhat$ on itself; combining the two actions we see that
$\Dhat$ is naturally a category over $\gt^2/\LG$.
  Since  $Rep(\LG)$ acts by {\em central} functors and the tensor endomorphism is compatible
with the central structure, we conclude
that the $\gt^2/\LG$ linear structure factors through the one for the fiber square
 $(\gt\times_\Lg \gt)/\LG=St/\LG$.

More precisely, we get the monoidal functor $D^b(Coh^\LG(\gt))\to \Dhat$
from the filtration following a strategy similar to the one in Example (5) above.
The first term of the filtration (the "lowest weight arrow") determines a 
functor from $D_{perf}^{\LG\times \LT}(C)$ where $C$ is a certain
affine scheme with an action of $\LG\times \LT$ with an open 
$\LG\times \LT$-equivariant embedding $\LG/\LU\to C$. The fact that the lowest weight arrow extends
to a filtration satisfying certain properties implies that complexes supported on $\partial ( \LG/\LU)=
C\setminus (\LG/\LU)$ act by zero. These ideas have already been used in \cite{AB}.

The fact that the action of the log monodromy endomorphisms on the category $D_{I^0I^0}$ of monodromic sheaves is  {\em nilpotent}, allows us to show that the $St/\LG$-linear structure on $D_{I^0I^0}$, $\Dhat$ factors through a canonical $\Sth/\LG$-linear
structure, where $\Sth$ is formal completion of $St$ at the preimage of $\N$. 

Once the $\Sth/\LG$-linear structure on $\Dhat$ is constructed,  any object $M\in \Dhat$
defines a functor
$D_{perf}^\LG(\Sth)\to \Dhat$, $\F \mapsto \F(M)$. We use the functor (denoted by $\Phi_{perf}$) corresponding to
the choice $M=\Xib$ where $\Xib$ is a certain tilting pro-object discussed in section \ref{sectionXi}.
This choice 
 can be motivated by the requirement of compatibility with
the equivalence $\Phi_{IW}^{I^0}$: the object $\Xib$ is obtained from the unit
object in $\Dhat$ by projection to $\Dhat_{IW}$ composed with its adjoint,
on the dual side this corresponds to the sheaf $pr_{Spr,2}^*pr_{Spr,2*}(\delta_*(\O_{\gth}))
\cong \O_{\Sth}$, where $\delta:\gt\to St$ is the diagonal embedding.
Thus the compatibility implies that $\Phi_{perf}(\OO)\cong \Xib$. The object $\Xib$
can also be thought of as a categorification of the element $\xi$ in the affine Hecke algebra,
thus it is closely related to Whittaker model, see \S \ref{Heckper}. 
 
 The fact that $\Phi_{perf}$ constructed this way is compatible with projection to $\DIW$ 
 follows from the properties of $\Xib$. 
 
 We then establish the equivalence $\Phi_{IW}^{I^0}$ as in \cite{AB}. Together with compatibilities
 between $\Phi_{perf}$ and $\Phi_{IW}^{I^0}$ 
 this implies that $\Phi_{perf}$ is a full embedding. 
  
Once $\Phi_{perf}$ is constructed and shown to be full, we get functors in the opposite direction
 $\Psib:\Dhat \to D^b(Coh^\LG(\Sth)) $, $\Psi:D_{I^0I^0} \to D^b(Coh^\LG(\St)) $.
 (The logic here is reminiscent of arguments in functional analysis where a map between
 spaces of (smooth rapidly decreasing) test functions induces a map between spaces of generalized functions going in the opposite direction).
 We show existence of $\Psi$, $\Psib$ and check that they are equivalences based on
  a general result relating the categories $D^b(Coh(X))$ and $D_{perf}(X)$ for an algebraic
   stack $X$.
We show that $D^b(Coh(X))$ embeds into the category of functors $D_{perf}(X)\to Vect$
and characterize the image of this embedding. The characterization makes use
of the standard $t$-structure on the derived category of coherent sheaves. In order
to apply the general criterion in our situation we show that, although the functor
$\Phi_{perf}$ is not $t$-exact with respect to 
 the natural $t$-structures on the two triangulated categories, it satisfies
 a weaker compatibility (see section \ref{with_t}).

At this point the equivalence \eqref{I0I0} is constructed, it remains to check its compatibility with
the convolution monoidal structure. 
 We use presentation of $\Dhat$ as homotopy category of complexes
of free-monodromic tilting (pro)sheaves introduced in \cite{BY} and recalled below. Using the observation that convolution of two free monodromic tilting sheaves is also a free monodromic tilting sheaf
we get an explicit monoidal structure on the category of tilting complexes, which is identified
with 
the monoidal structure on $\Dhat$. It turns out that $\Psib$ sends a free monodromic tilting
sheaf to a coherent sheaf (rather than a complex).
Thus the monoidal structure on the equivalence $\Phi_{I^0I^0}$ follows
 from
compatibility with the action on $\Dhat_{IW}$, since 
a sheaf in $D^b(Coh^\LG(\Sth))$
 can be uniquely reconstructed from the endo-functor of $D^b(Coh^\LG_\Nt(\gt))$ 
given by convolution with $\F$.

\subsection{Description of the content}
Sections \ref{monpro} and \ref{sectionXi} mostly recall the results of \cite{BY} while
section \ref{Funcon} recalls the material of \cite{AB} and extends it to the present slightly more general setting. 


As was indicated above, it is technically convenient to enlarge both categories in \eqref{I0I0} and construct
the equivalence
\begin{equation}\label{hatI0I0}
\Dhat \cong D^b(Coh^\LG(\Sth )).
\end{equation}

In section \ref{monpro} we recall the definition of $\Dhat$ and 
an extension of the formalism of tilting sheaves to this setting. We also present
a "monodromic" generalization of central sheaves \cite{KGB}.

Section \ref{Funcon} provides a generalization of the main result of \cite{AB} to the monodromic
setting. Namely, it establishes a monoidal functor $\Phid$ from the derived category of
equivariant coherent sheaves on the formal completion $\hatt{\gt}$ of $\gt$ at $\Nt$ to $\Dhat$.
(The composition of this functor with the equivalence \eqref{hatI0I0}
which will be established later is the direct image
under the diagonal embedding $\gt\to \St$, see Lemma \ref{commPsi}(b).) 
A variation of the argument allows us to define the action of
the tensor category $(D_{perf}^\LG(\St),\otimes_\O)$ on $D_{I^0I^0}$ and
$\Dhat$.

We also consider the projection of $\Dhat$ to the  Iwahori-Whittaker
category $\Dhat _{IW}$ and show that the composition of $\Phid$ with this projection induces an equivalence $\Phi^{I^0}_{IW}:D^b(Coh^\LG (\hatt{\gt}))\iso \Dhat_{IW}$.

Section \ref{sectionXi} is devoted to a particular object  $\Xib\in\hatt \P$
which will correspond to the structure sheaf of $\Sth$ under the equivalence.

In section \ref{Phip} we define a functor $\Phip$ from the subcategory
of perfect complexes $D_{perf}^\LG(\Sth)\subset
D^b(Coh^\LG(\Sth))$ to $\Dhat$ by sending an object
in the tensor category $(D_{perf}^\LG(\Sth),\otimes_\O)$ to the result of 
its action on $\Xib$. 

We then make a step towards establishing monoidal structure on our functors:
the functor $\Phip$ allows to define an action of $D_{perf}^\LG(\Sth)$ 
 on $D^{I^0}_{IW}$, while the category $D^b(Coh^\LG_\Nt(\gt))$ also carries such an action;
 we use properties of $\Xib$ to show that $\Phi^{I^0}_{IW} $ is compatible with these module structures. 
 Here (in contrast with the previous paragraph) $D_{perf}^\LG(\Sth)$ is equipped
 with the convolution product (notice that the subcategory
 $D_{perf}^\LG(\Sth)\subset D^b(Coh^\LG(\Sth))$ is easily seen to be closed
 under convolution).

This compatibility allows us to deduce that $\Phip$ is a full embedding and endow it
with the structure of a monoidal functor.

In  section \ref{with_t} we check a  property of $\Phip$ with respect to the natural
$t$-structures on the two categories.
In section \ref{gencrit} we give a general criterion allowing to extend
an equivalence from the category of perfect complexes to the bounded derived category of coherent sheaves.

In section \ref{theeq} we show that the criterion of section \ref{gencrit} applies,
by virtue of properties established in section \ref{with_t}, to the present
situation yielding \eqref{hatI0I0}. We then deduce
 \eqref{I0I} and \eqref{II} by means of a general lemma describing the equivariant
 constructible category via the monodromic one.
Section \ref{monstr} deals with technicalities on $DG$-models for convolution monoidal 
categories of sheaves needed to equip our equivalences with a monoidal structure.
The final section \ref{sec_fur} describes additional properties of our functors in relation to the
Frobenius automorphism (where $k=\bar{\Fq}$) and $t$-structures, as well as conjectural
generalizations and relation to Hodge $D$-modules.

\section{Monodromic sheaves and pro-object} \label{monpro}
\subsection{Generalities on monodromic sheaves}\label{genonmo}
Objects of $\P_{I^0I^0}$ are by definition perverse sheaves monodromic
with respect to both the left and the right action of $T$ on $\tii \Fl$. Thus we get
two actions of the group
$\Lambda\times \Lambda$  by automorphisms of the identity functor of $\P_{I^0I^0}$
coming respectively from the left and the right monodromy. 
Both actions  on each object are unipotent. 

Let $\Phat$ be the category of pro-objects $M$ in $\P_{I^0I^0}$ such that the coinvariants of
the left (equivalently, right) monodromy action belongs to $\P$. 
It is easy to see from the definitions that $\Phat$ is identified with the heart of
the natural $t$-structure on the pro-completion of the derived category $D_{I^0I^0}$
 introduced in \cite[Appendix A]{BY}.
Furthermore,  \cite[Corollary A.4.7]{BY} implies that the category
 $\Dhat=D^b(\Phat)$ is identified with that completion.\footnote{This way to define
 $\Dhat$ relies on the formalism of triangulated subcategories in the category of pro-objects 
in the derived category of constructible sheaves developed in  \cite[Appendix A]{BY} by Z.~Yun.
Alternatively  one could first define the category $\Phat$ as a subcategory in the category
of pro-objects in $\P$ and use free monodromic tilting objects to equip $\Dhat:=D^b(\Phat)$
with a monoidal structure.} 
  Thus  
 the construction of {\em loc. cit.} shows that 
$\Dhat$ is monoidal and contains $D_{I^0I^0}$ as a full  subcategory closed under the convolution product.
 An object $\F\in \Dhat$ belongs to $D_{I^0I^0}$ iff the monodromy automorphisms
 of $\F$ are unipotent. The formalism of {\em loc. cit.} applies also to Iwahori-Whittaker sheaves
 yielding the definition of an abelian category $\Phat_{IW}$ and triangulated category $\Dhat_{IW}\cong D^b(\Phat_{IW})$,
 so that $\Phat_{IW}$ is a full subcategory in the category of pro-objects in $\P_{IW}^{I^0}$ consisting
 of pro-objects with finite coinvariants of monodromy automorphisms, while $D_{IW}^{I^0}$ is a full
 subcategory on $\Dhat_{IW}$ consisting of objects where monodromy automorphisms are unipotent.
 Convolution action of $D_{I^0I^0}$ on $D^{I^0}_{IW}$ extends to an action of $\Dhat$ on $\Dhat_{IW}$.

 Let $\EE$ be the
{\em free prounipotent} rank one local system on $T$ (see \cite{BY}),  
thus $\EE=
\varprojlim \EE_n$ where $\EE_n$ is the local system whose fiber
at the unit element $1_T\in T$ is identified with the quotient of the group
algebra of tame fundamental group $\pi_1^{tame}(T)$ by the $n$-th power 
of augmentation ideal, where the action of monodromy coincides with the 
natural structure of $\pi_1^{tame}(T)$ module. Let $\Flt_w$ denote the preimage of $\Fl_w$
in $\Flt$.
The quotient 
$ \bI^0\bs \Flt_w$ is a torsor over $T$, choosing an arbitrary trivialization of the torsor
we get a projection $\tii \Fl_w\to T$ which we denote $pr_w$.
 Set
$\Delta_w=j_{w!}pr_w^*(\EE)[\dim \Flt_w]$, $\nabla_w=j_{w*}pr_w^*(\EE)[\dim \Flt_w]$.
The objects $\Delta_w$, $\nabla_w$ are defined uniquely up to a non-unique isomorphism, 
we call them
 a free-monodromic standard and costandard object respectively. One similarly defines
 $\Del_w^{IW}$, $\nab_w^{IW}\in \Phat_{IW}$.

\subsection{More on monodromic (co)standard pro-sheaves}\label{moreonmo}
A free prounipotent local system on $\tii \Fl_w$ is defined uniquely up to a non-unique isomorphism,
thus so are the (co)standard sheaves $\Delta_w$, $\nabla_w$. We now present 
geometric data allowing to fix these objects up to a canonical isomorphism.


Fix a maximal torus $T\subset B$ (recall that $\bI$ maps to $B$ under the evaluation
map $\GO\to G$); we get a canonical identification of 
$T$ with the abstract Cartan group of $G$, thus the group of coweights of $T$ is identified with 
$\Lambda$.
Thus for $w=\la\in \La\subset W$
the choice of a uniformizer $t\in F$
 defines an element $t_\la=\la(t)\in T_F\subset G(F)$; its image in $\Flt=\bG/\bI^0$
lies in the orbit of $I$ corresponding to $\lambda$.
This yields the choice of a point   $\overline{\lambda(t)}\in \bI^0\bs \Flt_\la$ which
gives a trivialization of the $T$-torsor, and hence
the choice of objects $\Delta_\la$, $\nab_\la$ defined uniquely up to a unique isomorphism.
We use the same notation for those canonically defined objects and the objects
defined earlier uniquely up to a non-unique isomorphism.

\begin{Lem}\label{DelDel}
a) We have isomorphisms $\Del_{w_1}*\Del_{w_2}\cong \Del_{w_1w_2}$,
$\nab_{w_1}*\nab_{w_2}\cong \nab_{w_1w_2}$ when $\ell(w_1w_2)=\ell(w_1)+\ell(w_2)$.

b) Assume that $w_1=\la_1,\, w_2=\la_2\in \La^+$ and let $\Del_{\la_i}$,
$\nab_{\la_i}$, $(i=1,2)$ be the canonically defined objects as above.
We have canonical isomorphisms  $\Del_{\la_1}*\Del_{\la_2}\cong \Del_{\la_1+\la_2}$,
$\nab_{\la_1}*\nab_{\la_2}\cong \nab_{\la_1+\la_2}$, which satisfy
the associativity identity for a triple $\la_1$, $\la_2$, $\la_3$.

c) $\Del_0=\nab_0$ is the unit object in $\Dhat$; and we have
a canonical isomorphism $\nab_w*\Del_{w^{-1}}\cong \nab_0$. 

d) We have $\Del_{w_1}*\nab_{w_2}\in \Phat$,
$\nab_{w_1}*\Del_{w_2}\in \Phat$ for all $w_1,w_2\in W$.

e) $\pi_*(\nab_w)\cong j_{w*}$, $\pi_*(\Del_w)\cong j_{w!}$ canonically.
\end{Lem}

\proof (a) and the first claim in (c) follow from \cite[Lemma 4.3.3]{BY},
\cite[Corollary 4.2.2 ]{BY}. 
A noncanonical isomorphism in the second statement in (c) follows from
the similar non-mondromic statement $j_{w*}*j_{w^{-1}!}\cong j_{e!}\cong j_{e*}$
by using \eqref{picomp1}, \eqref{picomp2} below and the observation that
any object $X$ in $\Dhat$ with $\pi_*(X)\cong j_{e*}$ is isomorphic to $\Delta_0$.
The non-monodromic statement is standard, in the case of a finite dimensional flag variety it
goes back to \cite{BeBe}; to check it directly one can reduce to the case when  $w$  is a simple reflection, then it amounts to an easy calculation based on the fact that 
$H_c^*({\mathbb A}^1\setminus
 \{0\} )=0$. 

Now given $\la_1$, $\la_2\in \La^+$ consider the locally closed subvariety in the convolution
diagram: $\Fltil_{\la_1}\boxtimes ^{\bI^0} \Fltil_{\la_2}\to \Fltil_{\la_1+\la_2}$.
Using the above trivializations of the 
 $T$ torsors
$\bI^0\bs \Flt_{\la_1}$, $\bI^0\bs \Flt_{\la_2}$, $\bI^0 \bs \Flt_{\la_1+\la_2}$
we can identify the quotient of  $\Fltil_{\la_1}\boxtimes ^{\bI^0} \Fltil_{\la_2}$
by $\bI^0$ with $T\times T$ and the quotient of $\Fltil_{\la_1+\la_2}$ by $\bI^0$ with $T$;
the quotient of the convolution map is readily seen to be the multiplication map
$T\times T\to T$. Since the convolution $\EE *_T \EE$ is canonically isomorphic
to $\EE[-\dim T]$ (here $*_T$ denotes convolution of sheaves on the group $T$) we get
the desired canonical isomorphism. Verification of the associativity identity is straightforward. 

Part (d) follows once we know that the functors $\M\mapsto \M*\nab_w$,
$\M\mapsto \nab_w*\M$ are right exact, while  $\M\mapsto \M*\Del_w$,
$\M\mapsto \Del_w*\M$ are left exact. These follow from their nonmonodromic analogues
by virtue of \eqref{picomp1}, \eqref{picomp2}, which are standard, see e.g. \cite[\S 5.1]{AB}
(in {\em loc. cit.} only $w_1$, $w_2$ of a special form are considered, but the argument applies generally). Alternatively, the statement follows from (4.4), (4.5) in the proof of 
\cite[Proposition 4.3.4]{BY}.

Finally, part (e) easily follows from the fact that cohomology of the free prounipotent
local system on $T$ is zero in degrees other than $r=\dim(T)$ and $r$-th cohomology 
is one dimensional. 
\epf

\subsection{Wakimoto pro-sheaves}\label{Wakimo}
Recall {\em Wakimoto sheaves} $J_\la\in \P_{II}$ characterized by:
$J_\la*J_\mu\cong J_{\la+\mu}$ for $\la,\mu\in \Lambda$ and $J_\la=j_{\la*}$
for $\la\in \La_+$, see \cite[3.2]{AB}. The following 
monodromic version follows directly from   Lemma \ref{DelDel}(b,c). 

\begin{Cor}\label{monfunTheta}
 There exists a monoidal functor $\Theta:Rep(\LT)\to \Dhat$ sending a dominant character $\la$
to $\nab_\la$ and an anti-dominant character $\mu$ to $\Delta_\mu$. Such a functor
is defined uniquely up to a unique isomorphism.  \epf

\end{Cor}

The image of a character $\la$ of $\LT$ under this functor will be called a 
free monodromic Wakimoto
sheaf and will be denoted by $\J_\la$.

Some of the basic properties of Wakimoto sheaves are as follows.

\begin{Lem}\label{proJ} We have:

a)   $\J_\la\in \Phat\subset \Dhat$.

b) $Hom^\bu(\J_\la,\J_\mu)=0$ for $\mu \not \preceq \la$ where $\preceq$ is 
the standard partial order on (co)weights.

c) $\pi_*(\J_\la)\cong J_\la$ canonically.
\end{Lem}

\proof a)  follows from Lemma \ref{DelDel}(d). 
b) is clear since 
$$Hom^\bu(\J_\la,\J_\mu)\cong Hom^\bu(\J_{\la+\eta},\J_{\mu+\eta})=
Hom^\bu (\nab_{\la+\eta}, \nab_{\mu+\eta}),$$
where $\eta\in \La$ is chosen so that $\la+\eta,\, \mu+\eta\in \La^+$.
The latter $Hom$ space vanishes when $\mu\not\preceq \la$ because
in this case $\Fltil_{\mu+\eta}$ is not contained in the closure of $\Fltil_{\la+\eta}$.  
The special case of part (c) when $\pm \la \in \La^+$
is contained in  Lemma \ref{DelDel}(e). To deduce the general case
we use isomorphisms:
\begin{equation} \label{picomp1}
\F* \pi_*(\G)\cong \pi_*(\F*\G)\in D_{I^0I},\ \ \ \ \ \F,\, \G\in D_{I^0I^0},
\end{equation}
\begin{equation}\label{picomp2}
\F* Res^I_{I^0}(\G)\cong \pi_*(\F)*\G\in D_{I^0I},\ \ \ \ \  \ \F\in D_{I^0I^0}, \, \G\in D_{II},
\end{equation}
which are easily checked using base change and transitivity of direct image.
Given $\la\in \La$ we write it as $\la=\la_+-\la_-$, $\la_\pm\in \La^+$ and apply
\eqref{picomp1} to $\F=\J_{\la_+}$, $\G=\J_{-\la_-}$. Applying then \eqref{picomp2} to 
$\F=\J_{\la_+}$, $\G=J_{-\la_-}$ we get statement (c).
\epf

\subsection{Generalized tilting pro-objects}
Recall that an object of $ \P$ is called {\em tilting}
if it  carries a standard and also a costandard filtration; here
a filtration is called (co)standard if its associated graded is a sum of (co)standard objects,
see e.g. \cite{BBM}.

An object  of $\Phat$ is called {\em free-monodromic tilting} if it carries a free-monodromic
standard and also a free-monodromic costandard filtration; here
a filtration is called (co)standard if its associated graded is a sum of
free-monodromic (co)standard objects, see \cite{BY}.

 Let $\TT\subset \P$ be the full subcategory of  tilting objects and
$\Th\subset \Phat$ denote the full subcategory of 
 free-monodromic tilting objects \cite[Definition A.7.1(1)]{BY}.

Let $Ho(\TT)$, $Ho(\TThat)$ denote the homotopy category of bounded complexes
of objects in $\TT$, $\TThat$ respectively.

The next Proposition summarizes the properties of tilting objects that will be used
in the argument.

\begin{Prop}\label{tiltprop}
a) The natural functors $Ho(\TT)\to D^b(\P)=D$,  $Ho(\TThat)\to D^b(\Phat)=\Dhat$
are equivalences. 

b) The convolution of two object in $\TThat$ lies in $\TThat$, thus $Ho(\TThat)$
has a natural monoidal structure. The natural functor $Ho(\TThat)\to \Dhat$
is a monoidal equivalence.

c) More generally, assume that $\F,\, \G\in \Dhat$ are represented by bounded complexes
 $\F^\bu$, $\G^\bu$ of objects in $\Phat$,
such that $\F^i*\G^j\in \Phat$. Then $\F*\G$ is represented by the total complex of the bicomplex
$\F^i*\G^j$. The same statement holds for $\F\in \Dhat$ and $\G\in \Dhat_{IW}$
or $\G\in D_{I^0I}$ 
represented by $\F^\bu$, $\G^\bu$ with $\F^i*\G^j\in \Phat_{IW}$ (respectively,
$\F^i*\G^j\in \P_{I^0I}$). 
 Given three  complexes $\F^\bu_1$, $\F^\bu_2$, $\G^\bu$ the two isomorphisms
 between $\F_1*\F_2*\G$ and the object represented by the complex
 $\CC^d=\oplusl_{i+j+l=d}\F^i_1*\F^j_2*\G^l$ coincide.
\end{Prop}

\proof The first statement in (a) appears in \cite[Proposition 1.5]{BBM}, the second
one (whose proof is similar) is a particular case of \cite[Proposition B.1.7]{BY}. 

The first statement in  (b) follows from \cite[Proposition 4.3.4]{BY}, while the second statement in  (b) 
and (c)
follow from \cite[Proposition B.3.1]{BY} applied to the functor of push-forward under the
convolution (or triple convolution) map  and the twisted product of the corresponding class of sheaves, cf. also 
\cite[Remark  B.3.2]{BY}. \epf

\begin{Rem}\label{extvan}
Implicit in Proposition \ref{tiltprop}(a) is  Ext vanishing:
 $$Ext^{>0}(\That_1,\That_2)=0=Ext^{>0}(T_1,T_2)$$
for $T_1,T_2\in \TT$, $\That_1,\That_2\in \Th$. 
 A stronger statement will be used later:
 $$Ext^{>0}(\hat M_1,\hat M_2)=0=Ext^{>0}(M_1,M_2)$$
where $M_1,\,M_2\in \P_{I^0I}$, $M_1$ admits a standard filtration
while $M_2$ admits a costandard filtration, $\hat M_1,\, \hat M_2\in \Phat$,
$\hat M_1$ admits a free-monodromic standard filtration, while $\hat M_2$ admits
a free-monodromic costandard filtration. The proof is immediate from
$Ext^{>0}(\Del_{w_1},\nab_{w_2})=0=Ext^{>0}(j_{w_1!}, j_{w_2*})$.
\end{Rem}

\begin{Prop}\label{dirim} 
 An object $M\in \Dhat$ admits a free-monodromic
(co)standard filtration iff $\pi_*(M)\in D_{I^0I}$ lies in $\P$ and admits
a (co)standard filtration. 
\end{Prop}

\proof The "only if" direction follows from Lemma \ref{DelDel}(e), while the "if"
direction is checked in \cite[Lemma A.7.2]{BY}. \epf

\begin{Cor}
 An object  $M\in \Dhat$ lies in $\Th$ iff $\pi_*(M)\in \TT$.
 
\end{Cor}

\begin{Prop}\label{TTT}  
a) For $T\in \Th$ the functors $\F\mapsto T*\F$ and $\F\mapsto \F*T$ are $t$-exact
(i.e. send $\P_{I^0I^0}$ to $\P_{I^0I^0}$ and $\Phat$ to $\Phat$).

b) For any $w\in W$ there exists a unique (up to an isomorphism) indecomposable
object  $T_w\in \TT$ whose support is the closure of $\Fl_w$.
There also exists a unique indecomposable object $\That_w\in \Th$
whose support is the closure of $\tii \Fl_w$.
We have $\pi_*(\That_w)\cong T_w$.

c)  
For $T\in \Th$ and $w\in W$ the objects $\Delta_w*T$, $T*\Delta_w\in \Phat$
 have a free-monodromic standard filtration, 
while the objects $\nabla_w*T$, $T*\nabla_w\in \Phat$ have a
free-monodromic costandard filtration.
\epf
\end{Prop}

\proof 
Parts (a,c) follows from the proof of \cite[Proposition 4.3.4]{BY}.
The first statement in (b) is standard, see e.g. \cite[Proposition 1.4]{BBM}.
The second one then follows from  \cite[\S A.7]{BY} which shows that
the functor $\M \mapsto \pi_*(\M)$ induces a bijection between isomorphism
classes of indecomposable objects in $\Th$ and $\TT$: \cite[Lemma A.7.2]{BY}
shows that $\pi_*:\Th\to \TT$, by \cite[Lemma A.7.3]{BY} it induces a surjective
map on isomorphism classes of objects, and \cite[Lemma A.7.4]{BY} implies that
this map is also injective, as it shows that for $\That_1,\, \That_2\in \Th$
an isomorphism $\pi_*(\That_1)\cong \pi_*(\That_2)$ can be lifted to an isomorphism
$\That_1\cong \That_2$. Alternatively, the second statement in (b) follows from 
\cite[Corollary 5.2.2]{BY}. \epf

\begin{Cor}\label{prescost}
 Convolution with a free-monodromic tilting object 
preserves the categories of objects
admitting a free-monodromic (co)standard filtration.
\end{Cor}

\subsection{Monodromic central sheaves}\label{moncen}
We need to extend the central functors of \cite{KGB} to the monodromic
setting. 

\subsubsection{A brief summary of \cite{KGB}}\label{KGBbrief}
Recall first the main result of \cite{KGB}. In our present notation it reads as follows.

For $V\in Rep(\LG)$ one defines an exact functor $\ZZ_V:\P_{II}\to \P_{II}$. 
One then constructs canonical isomorphisms
\begin{equation}\label{ZGa}
Z_V*\F\cong \ZZ_V(\F)\cong \F*Z_V,\ \ \ \ \F\in \P_{II};
\end{equation}
 \begin{equation}\label{ZGamon}
\ZZ_{V\otimes W}\cong \ZZ_V\circ \ZZ_W,
\end{equation}
where $Z_V=\ZZ_V(\delta_e)$, where $\delta_e=j_{e!}=j_{e*}$ is the skyscraper at 
the point $\Fl_e$. 

The two isomorphisms satisfy natural compatibilities (some are demonstrated in 
\cite{BOap}) which amount to saying that $V\mapsto Z_V$
is a tensor functor from $Rep(\LG)$ to Drinfeld center of $D_{II}$.

The goal of this subsection is to extend these results to the monodromic
setting.

Construction of the functor $Z_V$ is based on existence of a certain
deformation of the affine flag variety $\Fl$ and the convolution diagrams.

Let $C$ be a smooth algebraic curve over $k$ and fix a point $x_0\in C(k)$
and set $C^0=C\setminus \{x_0\}$.
The ind-schemes $\Fl_C$, 
$\Fl_C^{(2)}$, $Conv_C$, $Conv_C'$ were constructed in \cite{KGB}.
They are defined as moduli spaces parametrizing the following collections of data:

$\Fl_C=\{ (x, \, \E, \, \phi, \, \beta)\}$, where $x\in C$, $\E$ is a $G$-bundle on $C$, $\phi$
 is a trivialization
of $\E$ on $C\setminus \{x\}$ and $\beta\in (G/B)^\E_{x_0}$ is a point in the fiber of the associated fibration  
with fiber $G/B$ at $x_0$.

$\Fl_C^{(2)}=\{ (x,\, \E,\, \phi', \, \beta )\}$,
where $x$, $\E$, $\beta$ are as above and $\phi'$
 is a trivialization
of $\E$ on $C\setminus \{x,x_0\}$.

$Conv_C = \{ (x,\, \E,\,  \E',\, \phi, \, \psi, \beta,\beta' ) \}$ where $x$, $\E$,  $\beta$, $\phi$
are as above, $\E'$ is another $G$ bundle on $C$,  $\psi$ is an isomorphism $\E|_{C \setminus\{x\}} \cong \E'|_{C \setminus\{x\}}$,
while $\beta'\in (G/B)^{\E'}_{x_0}$.

 $Conv_C' = \{ (x,\, \E,\,  \E',\, \phi, \, \psi', \beta,\beta' ) \}$ where $x$, $\E$, $\phi$, $\E'$, $\beta$,
 $\beta'$ are as above, while $\psi'$ is an isomorphism $\E|_{C \setminus\{x,x_0\}} \cong \E'|_{C \setminus\{x,x_0\}}$.

These
 ind-schemes 
 come with a map to $ C$ satisfying the following properties.

The preimage of $x_0$ in $\Fl_C$
is identified with $\Fl$, while the preimage of $C\setminus \{x_0\}$ is identified
with $G/B \times \Gr_{C^0}$, where $\Gr_{C^0}$ is the {\em Beilinson-Drinfeld 
global Grassmannian}; thus the fiber of $\Fl_C$ over $y\in C^0(k)$ is (noncanonically)
isomorphic to $G/B\times \Gr$.

The preimage of $x_0$ in $\Fl_C^{(2)}$
is  identified with $\Fl$, while the preimage of $C^0$ is identified
with $\Fl \times \Gr_{C^0}$; thus the fiber of $\Fl_C^{(2)}$ over $y\in C^0(k)$ is (noncanonically)
isomorphic to $\Fl\times \Gr$.

To spell out the properties of $Conv_C$, $Conv_C'$ recall the convolution
space $\Fl \times^{\bI} \Fl $, which is the fibration over $\Fl$
with fiber $\Fl$ associated with the natural principal $\bI$ bundle over $\Fl$
using the action of $\bI$ on $\Fl$. We have the projection map $pr_1: \Fl\times ^\bI \Fl \to \Fl$
and the convolution map $conv:\Fl\times^\bI\Fl\to \Fl$ coming from multiplication
map of the group $\GF$. 

The fiber of both $Conv_C$ and $Conv_C'$ over $x_0$ is $\Fl\times ^\bI \Fl$;
the preimage of $C^0$ in $Conv_C$ is the product  $((G/B)\times ^\bI \Fl) \times
\Gr_{C^0}$, while the preimage of $C^0$ in $Conv_C'$ is identified with
$(\Fl\times^\bI (G/B) ) \times \Gr_{C^0}$.

One has canonical ind-proper morphisms $conv_C:Conv_C\to \Fl^{(2)}$,
$conv_C':Conv_C'\to \Fl^{(2)}$ whose fiber over $x_0$ is the convolution
map $conv$.

Starting from $V\in Rep(\LG)$ one can use the geometric Satake  isomorphism 
to produce a semi-simple perverse sheaf $S(V)$ on $\Gr_{C^0}$. 
For $\F\in Perv(\Fl)$ one gets a sheaf $\F\boxtimes S(V)$
on $\Fl\times \Gr_{C^0}\subset \Fl^{(2)}_C$. Taking nearby
cycles of that sheaf with respect to a local coordinate at $x_0$ 
one obtains a sheaf $\ZZ_V(\F)$ on $\Fl$.  

The spaces $Conv_C$, $Conv_C'$ and the maps $conv_C$, 
$conv_C'$  are used in \cite{KGB} to show that the functor $\ZZ_V |_{D_{II}}$
is isomorphic to both left and right convolution with a certain
object $Z_V\in \P_{II}$.

\subsubsection{The monodromic case}\label{KGBmon}
A straightforward modification of the definition from  \cite{KGB} 
yields spaces $\tii{\Fl}_C$, $\tii{\Fl}_C^{(2)}$, $\tii{Conv}_C$, $\tii{Conv}_C'$,
whose definition repeats the definition of ind-schemes in \S \ref{KGBbrief}
with the only difference that $\beta$, $\beta'$ are replaced by $\tilde \beta
\in (G/U)^\E_{x_0}$, $\tilde \beta'\in (G/U)^{\E'}_{x_0}$.
The following facts about these ind-schemes are proven by an argument similar
to that of  \cite{KGB} which deals with parallel statements about 
 ind-schemes from \S \ref{KGBbrief}.

The ind-schemes $\Fltil_C$, $\Fltil_C^{(2)}$,  $\tii{Conv}_C$, $\tii{Conv}_C'$ come with a map to $ C$
satisfying the following properties.

The preimage of $x_0$ in $\Fltil_C$
is identified with $\Fltil$, while the preimage of $C\setminus \{x_0\}$ is identified
with $G/U \times \Gr_{C^0}$; thus the fiber of $\Fltil_C$ over $y\in C^0(k)$ is (noncanonically)
isomorphic to $G/U\times \Gr$.

The preimage of $x_0$ in $\tii{\Fl}_C^{(2)}$
is  identified with $\tii{\Fl}$, while the preimage of $C\setminus \{x_0\}$ is identified
with $\tii{\Fl} \times \Gr_{C^0}$; thus the fiber of $\tii{\Fl}_C^{(2)}$ over $y\in C^0(k)$ is (noncanonically)
isomorphic to $\tii{\Fl}\times \Gr$. 

We will now use the convolution
space $\tii{\Fl} \times^{\bI^0} \tii{\Fl} $, which is a fibration over $\tii{\Fl}$
with fiber $\tii{\Fl}$ associated with the natural principal $\bI^0$ bundle over $\tii{\Fl}$
using the action of $\bI^0$ on $\tii{\Fl}$. We have the projection map $pr_1: \tii{\Fl}\times ^{\bI^0}
\tii{ \Fl} \to \tii{\Fl}$
and the convolution map $\tii{conv}:\tii{\Fl}\times^{\bI^0}\tii{\Fl}\to \tii{\Fl}$ coming from multiplication
map of the group $\GF$.

The fiber of both $\tii{Conv}_C$ and $\tii{Conv}_C'$ over $x_0$ is $\tii{\Fl}\times ^{\bI^0} \tii{\Fl}$;
the preimage of $C^0$ in $\tii{Conv}_C$ is the product  $((G/U)\times ^{\bI^0} \tii{\Fl}) \times
\Gr_{C^0}$, while the preimage of $C^0$ in $\tii{Conv}_C'$ is identified with
$(\tii{\Fl}\times^\bI (G/U) ) \times \Gr_{C^0}$.

One has canonical morphisms $\tii{conv}_C:\tii{Conv}_C\to \tii{\Fl}^{(2)}_C$,
$\tii{conv}_C':\tii{Conv}_C'\to \tii{\Fl}^{(2)}_C$ whose fiber over $x_0$ is the convolution
map $\tii{conv}$.

The main technical difference with the setting of \cite{KGB} recalled
in the previous subsection is that in contrast with the maps $conv_C$,
$conv_C'$ the maps $\tii{conv}_C$, $\tii{conv}'_C$ are {\em not ind-proper}.

For $V\in Rep(\LG)$ and $\F\in D_{I^0 I^0}$ we can form a complex
 $\F\boxtimes S(V)$ on $\tii{\Fl}\times \Gr_{C^0}\subset
 \tii{\Fl}^{(2)}_C$. Taking nearby cycles with respect to a local
 coordinate on $C$ near $x_0$ we get a complex which we denote
 $ \ZZb_V(\F)$.
 
The functor $ \ZZb_V$ obviously extends to the category $\Dhat$.
We set
 $\Zb_V=\ZZb_V(\Del_e)$.

\begin{Prop}\label{Zb}
a) Recall that $\pi: \tii \Fl\to \Fl$ is the projection. Then we have 
$\ZZb_V(\pi^*\F)\cong \pi^*(\ZZ_V(\F))$ canonically. 

b) $\ZZb_V$  is canonically isomorphic to the functors of both left and right convolution
with $\Zb_V$.

c) The map $V\mapsto \Zb_V$ extends to a central functor
$Rep(\LG)\to \Dhat$, i.e. to a tensor functor from $Rep(\LG)$ to the Drinfeld center
of $\Dhat$. 

d) We have a canonical isomorphism $\pi_*(\Zb_V)\cong Z_V$.
\end{Prop}

\proof a) follows from the fact that nearby cycles commute with pull-back
under a smooth morphism.

 The proof of (b,c) is parallel to the argument of \cite{KGB} and \cite{BOap}
respectively, with the following modification. The argument of {\em loc. cit.} uses
that the convolution maps and its global counterparts (denoted presently by
$conv_C$, $conv_C'$) are proper in order to apply the fact that nearby
cycles commute with direct image under a proper map.
The maps $\tii{conv}$, $\tii{conv}_C$, $\tii{conv}_C'$ are not proper,
thus we do not  a'priori have an isomorphism between the direct image under
$\tii{conv}_C$ or $\tii{conv}_C'$
of nearby cycles of a sheaf and nearby cycles of its direct image.
However, we do have a canonical map in one direction. 
If we start from a sheaf on $\tii \Fl$ which is the pull-back of a sheaf on $\Fl$,
then the map is an isomorphism because the sheaves in question 
are pull-backs under a smooth map of ones considered in \cite{KGB}.
Since all objects of $D_{I^0I^0}$   can be obtained from objects in the image of the pull-back
functor $D_{II}\to D_{I^0I^0}$ by successive extensions, the map in question is
an isomorphism for any $\F\in D_{I^0I^0}$, and claims (b,c) follows.

Using (b) we see that $\ZZb_V(\pi^*\F)\cong \pi^*(\F*\pi_*(\Zb(V)))$;  
thus (d) follows from (a). \epf

\subsubsection{Monodromy endomorphisms}\label{monodr}
Being defined as (the inverse limit of) nearby cycles sheaves, the objects $\Zb_V$, $V\in Rep(\LG)$
carry a canonical {\em monodromy automorphism}. It is known that the monodromy
automorphism acting on the sheaf $\Zs_V$ is unipotent, it follows that the one
acting on $\Zb_V$ is pro-unipotent. We let $m_V:\Zb_V\to \Zb_V$ denote the logarithm
of monodromy. 

It will be useful to have an alternative description of this endomorphism.
Consider the action of $\Gm$ on $\Fltil$ by {\em loop rotation}.
Since each $\bI\times \bI$ orbit on $\Fltil$ is invariant under this action,
every object of $\P_{I^0I^0}$ is $\Gm$ monodromic 
with unipotent monodromy. Thus every $\F\in \P_{I^0I^0}$ acquires a canonical
logarithm of monodromy endomorphism which we denote by $\mu_\F$. By passing
to the limit we also get a definition of $\mu_\F$ for $\F\in \Phat$.

\begin{Prop}\label{mmu}
a) We have $m_V=-\mu_{\Zb_V}$.

b) The logarithm of monodromy  defines  a tensor endomorphism of the functor 
$\Zb$, i.e. we have $m_{V\otimes W}=m_V* Id_{\Zb_W}+Id_{\Zb_V}*m_W$. 

\end{Prop} 
 
\proof  a) follows by the argument of \cite[5.2]{AB}, while (b) is parallel 
to \cite[Theorem 2]{KGB}.
 \epf

\subsubsection{Filtration of central sheaves by Wakimoto sheaves}
 
 It will be convenient to fix a total ordering on $\La$ compatible with addition
 and the standard partial order. This allows to make sense of an object in an abelian
 category with a filtration indexed by $\La$ and of its associated graded.
   
Recall that  the object $\J_\la$ was defined canonically up to a unique isomorphism starting from
 a fixed uniformizer $t$ of the local field $F$, while the central functor $\ZZ_V$ was defined
 using an algebraic curve $C$ with a point $x_0$ together with a fixed
 isomorphism between $F$ and the field of functions on the punctured formal
 neighborhood of $x_0$ in $C$.
 In the next Proposition we assume that $t$ is given by a local \' etale coordinate.
   We abbreviate $Z_{V_\la}$, $\Zb_{V_\la}$ to $Z_\la$, $\Zb_\la$ respectively. 
   
 \begin{Prop}\label{hwla}
a) For any $\lambda$ there exists a canonical surjective
morphism $\varpi_\lambda:\Zb_\la\to \J_\la$.
It is compatible with convolution in the following way: 
the composition of $\varpi_{\la+\mu}$ with the canonical 
map $\Zb_\la*\Zb_\mu\to \Zb_{\la+\mu}$ coming from the canonical
map $V_\la\otimes V_\mu=\Gamma(\O_\BB(\la))\otimes \Gamma(\O_\BB(\mu))\to V_{\la+\mu}
=\Gamma(\O_\BB(\la+\mu))$ equals $\varpi_\la*\varpi_\mu$.

b) The surjection $\varpi_\la$ extends to a unique filtration on $\Zb_\la$
indexed by $\La$ with associated graded isomorphic to a sum of Wakimoto sheaves
$\J_\mu$.

c) The filtration on $\Zb_\la$ is compatible with the monoidal structure on 
the functor $V\mapsto \Zb_V$, making $V \mapsto gr(\Zb_V)$ a monoidal functor.

 \end{Prop}

\proof
a) follows from the following standard geometric facts. 
Let $(\Fl_C^{(2)})_\la$ be the closure of $\Fl_e \times (\Gr_{C^0})_\lambda 
\subset \Fl^{(2)}_C$,
where $e\in W$ is the unit element and
 $(\Gr_{C^0})_\lambda$ is the locally closed subscheme in the Beilinson-Drinfeld
 global Grassmannian $\Gr_{C^0}$ whose intersection with a fiber of the projection
 to $C^0$ is the $\GO$ orbit $\Gr_\la$ (recall that such a fiber
 is identified with $\Gr$).
 Then $\Fl_\lambda\subset \Fl$ (where $\Fl$ is identified with the fiber of 
 $\Fl_C^{(2)}$ over $x_0$)
  is contained in the smooth locus of $(\Fl_C^{(2)})_\la$, it is open
in $(\Fl_C^{(2)})_\la\times_C \{x_0\}$. It follows that
$Z_\la$ which is by definition  the nearby cycles of $\delta_e\boxtimes
IC_\la$ (where $\delta_e$ is the skyscraper at $\Fl_e$) is constant on $\Fl_\la$ which is open in its support; see \cite[3.3.1, Lemma 9]{AB}.
Likewise, considering the preimage $(\Fltil_C^{(2)})_\la$ of $(\Fl_C^{(2)})_\la$
in $\Fltil _C^{(2)}$ we see that $\Fltil_\la$ is open in the support of 
$\Zb_\la$ and the restriction of $\Zb_\la$ to $\Fltil_\la$  is a free pro-unipotent
local system (shifted by $\dim(\Fltil_\la)$). This yields a surjection as in (a).
To see  existence of a canonical choice of the surjection it suffices
to see that the stalk of $\Zb_\la$ over the point $\overline{\lambda(t)}$ has
a canonical generator as a topological $\pi_1(T)$ module.
This follows from the fact that 
 the section $(1_{\Fltil},\la_\Gr):C^0\to \tii \Fl_C$ extends to $C$ and its value
at $x_0$ is $\lambda(t)_{\tii \Fl}$. 

Uniqueness of the filtration in (b) follows from the fact that $Hom^\bu(\J_\la,\J_\mu)=0$
for $\mu\not \preceq \la$ (Lemma \ref{proJ}(b)). Together with the isomorphism $\J_\la*\J_\mu\cong \J_{\la+\mu}$ 
this also implies compatibility with convolution and the monoidal property. 
Existence of the filtration is equivalent to the fact that $\J_\mu *\Zb_\la$ admits a free-monodromic
costandard filtration when $\mu$ is deep in the dominant chamber (more precisely,
when $\mu +\nu$ is dominant for any weight $\nu$ of $V_\lambda$).
This follows from Proposition \ref{dirim} and the corresponding fact about the sheaves $Z_\la$ established in \cite[\S 3.6]{AB}. \epf 

{\bf Remark.}
It is shown in \cite{AB} that the multiplicity of $J_\mu$ as a subquotient of $Z_V$ equals
the multiplicity of the weight $\mu$ in representation $V$. It is clear that the same multiplicity
also equals the multiplicity of $\J_\mu$ as a subquotient of $\Zb_V$. This is also a consequence
of \eqref{hatI0I0}, since that equivalence sends $\Zb_V$ to $V\otimes \O_{\gth}$ which admits
a filtration whose associated graded is a direct sum of line bundles on $\gth$ with the above multiplicities. 

The objects $Z_\la$, $\Zb_\la$ can be thought of as a categorification of the central elements $S_\la$ in the affine Hecke algebra introduced by Lusztig in \cite{LAst}; the filtration by Wakimoto sheaves
with the above multiplicities  categorifies formula (8.2) of
{\em loc. cit.}
\subsubsection{Torus monodromy}\label{tormonsec}

Every sheaf in $\P_{I^0I^0}$ is monodromic with respect to $T \times T$ with unipotent
monodromy, since every irreducible object in $\P_{I^0I^0}$ is equivariant.
Thus taking logarithm of monodromy we get an action of $Sym(\t \oplus \t)$ on
$\P_{I^0I^0}$ by endomorphisms of the identity functor.

\begin{Lem}\label{torusmon}
a) 
The action of the two copies of $\t$ on $\Del_w$, $\nab_w$
differ by twist with the element $\bar{w}\in W_f$, where we use the notation $w\mapsto \bar{w}$
for the projection $W\to W_f$. In particular, the left action of $\t$ on the objects $\Del_\la$, $\nab_\la$,
$\la \in \La$,
coincides with the right one. 

b) The left action of $\t$ on the objects
 $\J_\la$, $\la\in \La$, $\Zb_\mu$, $\mu\in \La^+$ coincides with the right one.
 
 c) The action of loop rotation monodromy on $\Del_\la$, $\nab_\la$, $\J_\la$
coincides with the image of coweight $d\la\in \t$ under the above
action of $\t$.
\end{Lem}

\proof 
Let $x_w$ be a point in $\Fltil_w$ such that the orbit of $x_w$ with respect to the left
and the right action of $T$ coincide. Then restriction from $\Fltil_w$ to $T(x_w)$ is an equivalence 
between $\bI\times \bI$ unipotently monodromic sheaves on $\Fltil_w$ and unipotent
local systems on $T(x_w)$. Also for $t\in T$
 we have $t (x)=x (\bar{w}(t))$
in the self-explanatory notation. This implies (a).

The statement about $\J_\la$   in (b) for $\pm \la\in \La^+$ follows from (a), this yields
the general case because of compatibility of torus monodromy with convolution.

The statement about $\Zb_\la$ in (b) follows from the construction with 
nearby cycles, since the action of $T^2$ on $\Fltil\times \Gr_{C^0}$ (where
$T^2$ acts trivially on the second factor) extends to an action 
on $\Fltil_C$.

Finally, part (c) is a consequence of the following observation.
Let $R$ denote the loop rotation action of $\Gm$ on $\Fltil$.
  Then for  $\la\in \Lambda$
  let $h_\la:\Gm \to T$ be the corresponding homomorphism (see \S \ref{moreonmo}).
 Then  we have $R(s)(x_\lambda)=h_\la(s)(x_\lambda)$. \epf

\section{Construction of functors} \label{Funcon}

\subsection{A functor from $D^b(Coh^\LG(\gt))$}
Recall that $\gth$ denotes the formal completion of $\gt$ at $\Nt$.

In this subsection we construct a monoidal functor $\Phid:D^b(Coh^\LG(\gth))\to \Dhat$. 
The functor we presently construct is compatible with the equivalence
$\Phi: D^b(Coh^\LG(\hatt{St}))\cong \Dhat$ that will be established in section \ref{theeq}
as follows: $\Phid\cong \Phi\circ \delta_*$,
where $\delta: \gth\to \hatt{St}$ is the diagonal embedding, see Lemma \ref{commPsi}(b). 

The construction is parallel to that of \cite[\S 3]{AB}, so we only recall the main ingredients of the construction referring the reader to \cite{AB} for details. 

Following
the strategy outlined in section \ref{strasp},
we first list compatibilities satisfied by the functor $\Phid$ which characterize it uniquely.

\subsubsection{Line bundles and Wakimoto sheaves} 
Recall that for $\la\in \La$ the corresponding line bundle on $\BB$ is denoted
by $\O_\BB(\la)$, while $\OO_{\gth}(\la)$ is its pull-back to $\gth$.
The functor $\Phid$ satisfies:
 $$\Phid(\OO_{\gth}(\la))\cong \J_\la.$$
This isomorphism is compatible with the monoidal structure on the two categories,
i.e. it provides a tensor isomorphism between the functor $\Theta$ 
(see Corollary \ref{monfunTheta}) and
 the composition of $\Phid$ with the tensor
functor $\la\mapsto \O_{\gth}(\la)$.

\subsubsection{Twists by representations and central functors}\label{cenfun}

We have a tensor functor $Rep(\LG)\to Coh^\LG(\gth)$ sending a representation 
$V$ to $V\otimes \O$. 
Composition of $\Phid$ with this functor is isomorphic to the tensor functor $V\mapsto \Zb_V$
(see section \ref{moncen}).

\subsubsection{The lowest weight arrow}\label{highwei}
We have a familiar morphism
  of $\LG$-equivariant vector bundles on $\BB$:
$\O_\BB\otimes V_\la\to \O_\BB(\la)$.
We can pull it back to
$\gth$ to get a morphism in $Coh^\LG(\gth)$. The functor $\Phid$ sends this 
arrow to the map $\varpi_\la$ (notations of Proposition \ref{hwla}). 

\subsubsection{Log monodromy endomorphism}\label{logmon}
 Notice that for $x\in \Lg$, $\F\in Coh^\LG(\Lg)$ the centralizer 
of $x$ in $\LG$ acts on the fiber $\F_x$ of $\F$ at $x$. Differentiating this action one
gets the action of the Lie algebra of the centralizer $\fz(x)$. In particular, $x\in \fz(x)$
produces a canonical endomorphism of $\F_x$, it is easy to see
that it comes from a uniquely defined endomorphism of $\F$, which we denote
by $\fm_\F$ (in \cite{AB} we used notation $N_\F^{taut}$). 
It is clear that restricting $\fm$ to sheaves of the form $\F=V\otimes \O_{\gth}$
one gets a tensor endomorphism of the tensor functor $V\mapsto V\otimes \OO_{\gth}$.

We require that $\Phid$ sends $\fm_{V\otimes \OO}$ to the monodromy endomorphism 
$m_V$.  

\subsubsection{Projection to $\Lt^2$ and torus monodromy}\label{t2mon}
 We have a canonical
map $\gt\to \Lt$, thus the category $D^b(Coh^\LG(\gt))$ is canonically an $\O(\Lt)$-linear category, i.e. $\Lt^*=\t$ acts on it by endomorphisms of the identity functor.
This induces a pro-nilpotent action of $\t$ on $D^b(Coh^\LG(\gth))$. 

According to section \ref{tormonsec}, we have two commuting pronilpotent $\t$ actions on 
$\Phat$ and hence on $\Dhat$. The functor $\Phid$ intertwines the action of $\t$
described in the previous paragraph with either of the two monodromy actions.

\subsection{Monoidal functor from sheaves on the diagonal}\label{monfundiag}

We use a version of homogeneous coordinate ring construction and Serre description
of the category of coherent sheaves on a projective variety.

Let $C_{\gt}$ be the preimage of $\gt\subset \Lg \times \BB$ under the morphism
$\Lg\times \LG/\LU\to \Lg \times \BB$. Let $\overline{\LG/\LU}$ denote the affine closure of $\LG/\LU$. Notice that $\LG/\LU$ can be realized as a locally closed subscheme, namely as the orbit of a highest weight vector in the space $V$ of a representation of $G$.
Moreover, if the representation $V$ is chosen appropriately, the closure of $\LG/\LU$ in $V$ is isomorphic to $\overline{\LG/\LU}$. Define the action of the abstract Cartan $\Lt$ on $V$ such that $t\in \Lt$ acts on an irreducible summand with highest weight $\lambda$ by  the scalar $\langle \la, t\rangle$. Then define a closed
subscheme\footnote{Here notations diverge from that of \cite{AB}, there "hat" was used to denote the affine cone, while in the present paper it is used to denote completions.}
$\Cb_{\gt}\subset \Lg \times \Lt \times \overline{\LG/\LU}$ by the equation $x(v)=t(v)$,
$x\in \Lg$, $t\in \Lt$, $v\in \overline{\LG/\LU}\subset V$. It is easy to see that
$C_{\gt}$ is an open subscheme in $\Cb_{\gt}$. More precisely, without loss of generality we can assume
that representation $V$  is multiplicity free, i.e. it is a sum of pairwise non-isomorphic irreducible
representations; then $C_{\gt}$ is identified with the intersection of $\Cb_{\gt}$ with the open set of vectors which
have a nonzero projection to each irreducible factor. 

We leave the proof of the following statement to the reader.

\begin{Prop}\label{HomCb}
A) The scheme $\Cb_{\gt}$ does not depend on the choice of $V$ subject to the above conditions.

B) Consider the category of commutative rings over $\O(\Lt)$ equipped with a $\LG$ action which fixes the image of $\O(\Lt)$.

The following two functors on that category are canonically isomorphic:

\begin{enumerate}
\item
$R\mapsto Hom (Spec(R), \Cb_{\gt})$ where $Hom$ stands for maps compatible with the $\LG$ action and the map  to $\Lt$.

\item $R \mapsto \{(E_V, \iota_V)\ |\ V\in Rep(\LG)\}$. Here for $V\in Rep(\LG)$,
$E_V\in End_R(V\otimes R)$ and $\iota_V$ is a map of $R$-modules $R\to V\otimes R$.
This data is subject to the requirements:

\begin{enumerate}
 \item functoriality in $V$;

\item  $E_{V\otimes W}=E_V\otimes Id_W + Id_V\otimes E_W$;

\item  $\iota_{V\otimes W}=\iota_V\otimes \iota_W$;

\item     The action of $E_{V_\la}$ on the image of $ \iota_{V_\la}$ 
coincides with the action of the element in $R$ which is the image of 
$\la\in \Lt^*$ under the map $\Lt^*\to R$.

\end{enumerate}
\end{enumerate}

\end{Prop}

\subsubsection{Deequivariantization}\label{deeqsubs} (cf. \S \ref{linover}(3)) 
We will make use of the following construction. Let
$\C$ be an additive category linear over the field $k$, with an action of the tensor
category $Rep(H)$ of (finite dimensional algebraic) representation of $H$, where $H$ is a reductive algebraic group over $k$.
(Recall that $k$ is algebraically closed of characteristic zero; the definition is
applicable under less restrictive assumptions). 

We can then define a new category $\C_{deeq}$ by setting $Ob(\C_{deeq})=
Ob(\C)$, $Hom_{\C_{deeq}}(A,B)=Hom_{Ind(\C)}(A,\unO(H)(B))$, where $Ind(\C)$
is the category of Ind-objects in $\C$ and $\unO(H)$ is the object of $Ind(Rep(H))$
coming from the module of regular functions on $H$ equipped with the action of $H$ by left
translations. 
Using that $H$ is reductive over an algebraically closed field of characteristic zero we can
write the Ind-object $\unO(H)$ as $\oplusl _{\unV\in IrrRep(H)} \unV\otimes V^*$,
where  $IrrRep(H)$ is a set of representatives for isomorphism classes of irreducible
$H$ modules and for a representation $\unV\in IrrRep(H)$ we let $V$ denote the underlying vector space. Thus we have
 $$Hom_{deeq}(X,Y)=\oplusl_{\unV\in IrrRep(H)}  Hom(X,\unV(Y))\otimes V^*.$$

For example, if $\C=D^b(Coh^H(X))$ where $X$ is a scheme equipped with an $H$ action
then for $\F,\G\in \C$ we have $Hom_{deeq}(\F,\G)=Hom_{D^b(Coh(X))}(\F,\G)$. 

When we need to make the group $H$ explicit in the above definition we write
$Hom_{deeq}^H$ instead of $Hom_{deeq}$. 

The category $\C_{deeq}$ is enriched over $H$-modules, i.e. every $Hom$ space
carries the structure of an $H$-module compatible with composition.
We refer the reader to \cite{AG} for further details and to \cite{Gaover} for a more general
construction (cf. also \cite{AB}, proof of Proposition 4).

This formalism comes in handy for deducing the following statement.

Let $Coh_{fr}^{\LG\times \LT}(\Cb_{\gt})$ be the full subcategory in $Coh^{\LG\times \LT}(\Cb_{\gt})$ consisting of objects of the form $V\otimes \O$, $V\in Rep(\LG\times \LT)$. In other words, objects of $Coh_{fr}^{\LG\times \LT}(\Cb_{\gt})$ are representations of $\LG\times \LT$ and morphisms are given by $Hom(V_1,V_2)=Hom_{Coh^{\LG\times \LT}(\Cb_{\gt})}(V_1\otimes \O, V_2\otimes \O)$. This is a tensor category under the usual tensor product of vector bundles.

\begin{Cor}\label{extefun}
Let $\CC$ be a $k$-linear additive monoidal category. Suppose we are given


1) A tensor functor $F:Rep(\LG\times \LT)\to \CC$.

2) A tensor endomorphism $E$ of $F|_{Rep(\LG)}$, $E_{V_1\otimes V_2}=E_{V_1}\ot
Id_{F(V_2)} + Id_{F(V_1)}\ot E_{V_2}$. 

3) An action of $\O(\Lt)$ on $F$ by endomorphisms, so that for $f\in \O(\Lt)$
we have $f_{V_1\ot V_2}=f_{V_1}\ot Id_{F(V_2)}=Id_{F(V_1)}\ot f_{V_2}$.

4) A "lowest weight arrow" $\hw_\la:F(V_\lambda) \to F(\lambda)$ making the following
diagrams commutative: 
$$\begin{CD}
F(V_\la\otimes V_\mu) @>>> F(V_{\la+\mu} )\\
@V{\hw_\la\ot_\CC \hw_\mu}VV @VV{\hw_{\la+\mu}}V \\
F({\la})\otimes_\CC F(\mu) @>{\sim}>> F(\la+\mu)
\end{CD}$$

$$\begin{CD}
F(V_\la)@>{\hw_\la}>> F(\la)\\
@V{E_{V_\la}}VV @VV{\la}V \\
F(V_\la)@>{\hw_\la}>> F(\la)
\end{CD} $$
where the right vertical map is the action of the element $\la\in \t\subset \O(\Lt)$ coming from (3).

Then the tensor functor $F$ extends uniquely to a tensor functor
$Coh^{\LG\times \LT}_{fr}(C_{\gt})\to \CC$, so that $E$ goes to the tautological endomorphism $\fm$
(see section \ref{logmon}), the action of $\ft$ comes from the projection $\Cb_{\gt}\to \Lt$ and 
the lowest weight arrow comes from the map described in section \ref{highwei}.
\end{Cor}

\proof Extending the functor $F$ to a functor $Coh^{\LG\times \LT}_{fr}(\Cb_{\gt})\to \CC$ is equivalent
to providing a $\LG\times \LT$-equivariant homomorphism $\O(\Cb_{\gt})\to Hom_{deeq}^{\LG\times \LT}
(1_\CC,1_\CC)$, where we used the action of $Rep(\LG\times \LT)$ on $\CC$ given by
$V:X\mapsto F(V)X$. We now apply  Proposition \ref{HomCb} to the ring 
$R:=Hom_{deeq}^{\LG\times \LT}
(1_\CC,1_\CC)$. The action described in (3) provides it with a structure of a ring over $\O(\Lt)$;
the tensor endomorphism $E$ from (2) yields the collection $(E_V)$  
and the arrows $\varpi_\lambda$ induce the maps $\iota_V$ as in  Proposition \ref{HomCb}.
The commutative diagrams in part (4) of the Corollary imply identities (c,d) in Proposition \ref{HomCb}. Thus existence of a unique functor $F$ with above properties follows
from Proposition \ref{HomCb}. \epf

\subsubsection{The functor $\Phi_{diag}$} 
We now construct a monoidal functor $ \Phi_{diag}^{fr}:Coh_{fr}^{\LG\times \LT}(\Cb_{\gt})\to \Phat$
(more precisely, a monoidal functor to $\Dhat$ taking values in $\Phat$).

The functor is provided by Corollary \ref{extefun}:
we have a tensor functor from $Rep(\LG)$ to $\Dhat$ coming from the central functors (subsection \ref{cenfun}), and another commuting one from $Rep(\LT)$ to $\Dhat$ coming from Wakimoto sheaves (\ref{Wakimo});
the logarithm of monodromy endomorphisms (subsection \ref{logmon}) provide endomorphism 
$E$ while the torus monodromy (section \ref{tormonsec}) gives an action of 
$\t =\Lt^*$ (notice that due to Lemma \ref{torusmon}(b) we get the same action by using
either left or right torus action).
The morphisms described in subsection \ref{highwei} yield arrows $\varpi_\la$.
The conditions of Corollary \ref{extefun} are checked as follows.
Condition (2) follows from Proposition \ref{mmu}(b). 
Condition (3) is clear from
 compatibility of the convolution map with the torus action. 
 The first commutative diagram in condition (4) follows from Proposition \ref{hwla}(a), while
the second one is obtained   by comparing
Proposition \ref{mmu}(a) with Lemma \ref{torusmon}(c).  \epf

\subsection{"Coherent" description of  the anti-spherical (generalized Whittaker) category}
Consider the composition $Ho(Coh_{fr}^{\LG\times \LT}(\Cb_{\gt}))\to Ho(\Phat)\to \Dhat$
where $Ho$ denotes the homotopy category of complexes of objects in the given
additive category and the first arrow is induced by $\Phi_{diag}^{fr}$; this composition will 
be denoted by $\Phi_{diag}^{Ho}$.

Let $Acycl\subset  Ho (Coh_{fr}^{\LG\times \LT}(\Cb_{\gt}))$ be the subcategory
of complexes whose restriction to the open subscheme $C_{\gt}$ is
acyclic.

\begin{Prop}\label{Acyclto0}
 The functor $\Phi_{diag}^{Ho}$ sends the subcategory $Acycl$ of acyclic complexes to  zero.
\end{Prop}

\proof Proposition follows from existence of a filtration on $\Zb_\la$ with associated graded
being the sum of Wakimoto sheaves (Proposition \ref{hwla}(b)) by an argument parallel
to \cite[3.7]{AB}. \epf

The perfect  derived category of modules over a positively graded algebra  (by which we mean
a $\Zet$-graded algebra with vanishing negative components and component of degree zero
generated by the unit element) over a field of characteristic zero is well known
to be equivalent to the homotopy category of free graded modules, the same applies to 
equivariant modules, where an algebra is assumed to be equipped with an action of a reductive group. Applying this to $\O(\Cb_{\gt})$ we see that $D_{perf}^{\LG\times \LT}(\Cb_{\gt})\cong
Ho (Coh_{fr}^{\LG\times \LT}(\Cb_{\gt}))$.  Since $\gt$ is smooth,
$D^b(Coh^\LG(\gt))=D_{perf}^\LG(\gt)$,
thus \eqref{Idem} shows that\footnote{In fact, $K^0(Coh^\LG(\gt))=K^0(Coh^\LB(\Lb))=K^0(Coh^\LT(pt))$ is generated
by the classes of equivariant line bundles, thus the functor 
$Ho (Coh_{fr}^{\LG\times \LT}(\Cb_{\gt}))\to D^b(Coh^\LG(\gt))$ induces a surjection on $K^0$.
 By a standard argument (see e.g. \cite[Corollary 0.10]{Nee1}) this  implies:
 $D^b(Coh^\LG(\gt))\cong Ho (Coh_{fr}^{\LG\times \LT}(\Cb_{\gt}))/Acycl.$} 
$$D^b(Coh^\LG(\gt))\cong Idem(Ho (Coh_{fr}^{\LG\times \LT}(\Cb_{\gt}))/Acycl) .$$
 Thus
 Proposition \ref{Acyclto0} yields a functor $D^b(Coh^\LG(\gt))\to \Dhat$. The log monodromy
 action of $\t$ on the identity functor of $\Dhat$ is pro-unipotent, thus it extends canonically
 to an action of the completion of $\O(\Lt)$ at the maximal ideal of $0$.
  It is easy to deduce
 that the functor factors canonically through a functor $D^b(Coh^\LG(\gth))\to \Dhat$,
 we denote the latter functor by $\Phi_{diag}$.  

A closely related functor $F:D^b(Coh^\LG(\Nt))\to D_{II}$ was constructed
in \cite[\S 3]{AB}.

\begin{Lem}\label{sogl} Let $i:\Nt\to \gth$ be the embedding.
The following diagrams commute up to a natural isomorphism:
$$ \begin{CD}
D^b(Coh^\LG(\gth)) @>{ \Phi_{diag}}>> \Dhat \\
@V{i^*}VV @VV{ \pi_*}V \\
D^b(Coh^\LG(\Nt)) @>{Res^I_{I^0} \circ F}>> D_{I^0I}
\end{CD}
$$
$$ \begin{CD}
D^b(Coh^\LG(\Nt)) @>F>> D_{II}\\
@V{i_*}VV @VV{Res^I_{I^0}\circ \pi^*[r]}V \\
D^b(Coh^\LG(\gth)) @>{\Phi_{diag}}>> \Dhat
\end{CD}
$$
where
 $Res$ stands for restriction of equivariance,  
 and $r={\mathrm{rank}} (\g)$.
\end{Lem}

\proof  To check commutativity of the first diagram it suffices to prove the similar
commutativity for functors on the categories of finite complexes in 
$Coh_{fr}^{\LG\times \LT}(\Cb_{\gt})$.
This follows 
from the isomorphisms $\pi_*(\Zb_V)\cong Z_V$ (Proposition \ref{Zb}(d)),
$\pi_*(\J_\la)\cong J_\la$ (Lemma \ref{proJ}(c)) which are easily seen to be compatible with monodromy endomorphism and lowest weight arrows.

Now using commutativity of the first diagram we get a natural transformation between the two compositions in the second diagram due to the isomorphisms:
\begin{multline*}
Hom(Res^I_{I^0}\pi^*(F(\F))[r], \Phi_{diag}(i_*(\F)))\cong Hom(Res^I_{I^0} F(\F)[r],
\pi_*\Phi_{diag}(i_*(\F)))\cong \\ Hom(Res^I_{I^0} F(\F)[r], Res^I_{I^0}F(i^*i_*(\F))),
\end{multline*}
which yield the desired arrow since $Id[r]$ is a canonical direct summand in 
the functor $i^*i_*$. The constructed arrow is nonzero, hence it is an isomorphism for 
$\F=\O_{\Nt}$, as both compositions are then isomorphic to the skyscraper 
$\delta_e=j_{e!}=j_{e*}$ and $Hom(\delta_e,\delta_e)$ is
one dimensional. Also it is easy to see that the arrow is compatible with action of 
the tensor category $Coh_{fr}^{\LG\times \LT}(\Cb_{\gt})$, thus it is an isomorphism for $\F$ in a generating
set of $D^b(Coh^\LG(\Nt))$, hence it is an isomorphism for all $\F$.
\epf

\subsubsection{Equivalences $\Phib_{IW}$, $\Phi_{IW}^{I^0}$} We are now ready to establish \eqref{PhiIWmon}.

The functor $Av^{IW}:D_{I^0I^0}\to D_{IW}^{I^0}$ introduced at the end of \S \ref{furnot} extends to a functor between the completed
categories $\Dhat$, $\Dhat_{IW}$ introduced in \S \ref{genonmo}, we will use the same
notation for this extension.

\begin{Prop}
a) The functor $\Phib_{IW}:=Av^{IW}\circ \Phi_{diag}:D^b(Coh^{\LG}(\gth))\to
\Dhat_{IW}$ is an equivalence.

b) The functor $\Phib$ restricts  to an equivalence $\Phi_{IW}^{I^0}:D^b(Coh_\Nt(\gt))
\to D_{IW}^{I^0}$.
\end{Prop}

\proof 
We first show that $Av^{IW}\circ \Phid$ is fully faithful. It suffices to show that
$$Hom(\F,\G)\iso Hom (Av^{IW} \Phid(\F),Av^{IW} \Phid(\G))$$
when $\F=i_*(\F')$, $\F'\in D^b(Coh^\LG(\Nt))$, then the statement follows
since the image of $i_*$ generates $D^b(Coh^\LG_\Nt(\gt))$ under extensions,
so we get the isomorphism for $\F,\,\G\in D^b(Coh_\Nt(\gt))$. Passing to the limit
we then get the isomorphism for all $\F,\, \G\in D^b(Coh^\LG(\gth))$.


Using
 the parallel statement in the non-monodromic setting proved in \cite[\S 4]{AB}
 and the first
 commutative diagram in Lemma \ref{sogl} (or rather the statement obtained from it
 by left--right swap)
we get: 
\begin{multline*}
Hom (i_*(\F'),\G)
\cong 
Hom(\F',i^*\G[-r])
\cong 
Hom_{D_{IW}^{I}}\left( ^I Av^{IW} F(\F')), ^I Av^{IW} F(i^*\G)[-r]\right)
 \cong \\
Hom _{D_{IW}^{I}}\left( ^I Av^{IW}F(\F'), (Av_{I^0}^I)^{left}_* Av^{IW}(\Phid(\G)) [-r]\right)
\cong \\
Hom_{D_{IW}^{I^0}}\left( (Res^I_{I^0})^{left}(^I Av^{IW}F (\F')) ,  Av^{IW}
(\Phid(\G))[-r] \right) \cong\\
Hom \left( Av^{IW}\Phid(\F), Av^{IW}\Phid(\G') \right),
\end{multline*}
where  we used that $i^*[-r]$ is {\em right} adjoint to $i_*$. Here  $(Av_{I^0}^I)^{left}_*$
is the {\em right} adjoint to the restriction of equivariance functor which can be thought 
of as a direct image under the morphism of stacks $\bI^0 \bs \tii \Fl \to \bI \bs \tii \Fl$
(recall that $Av$ stands for  the $!$-direct image under that morphism).

This shows that the functor is fully faithful. Again using
 the parallel statement in the non-monodromic setting 
 and Lemma \ref{sogl} we see that the essential image of $\Phid$ contains
 the image of the functor of restricting the equivariance $D_{IW}^I\to D^{I^0}_{IW}$,
 since $D^{I^0}_{IW}$ is generated by irreducible perverse sheaves which are $I$-equivariant,
 the essential image contains $D^{I^0}_{IW}$, this proves part (b).
 Any object in $\Dhat_{IW}$ is an inverse  limit of objects in 
 $D_{IW}^{I^0}$, moreover, its image under the functor between the categories
 of pro-objects
   $Pro(D_{IW}^{I^0})
 \to Pro(D_{IW}^I)$
 induced by the averaging functor
 lies in $D_{IW}^I\subset Pro(D_{IW}^I)$. This shows that such an object is isomorphic
 to the image of a pro-object in $D^b(Coh^\LG_\Nt(\gt))$
 whose image under $i^*:Pro(D^b(Coh^\LG_\Nt(\gt)))\to Pro(D^b(Coh^\LG(\Nt)))$ lies
 in $D^b(Coh^\LG(\Nt))$. An object in $ Pro(D^b(Coh^\LG_\Nt(\gt)))$ satisfying the latter
 property 
 is easily seen to lie in $D^b(Coh^\LG(\gth))
 \subset Pro(D^b(Coh^\LG_\Nt(\gt)))$, this
  implies essential surjectivity in part (a). \epf

\subsection{$\Dhat$ is a category over $\Sth/\LG$, $D_{I^0I}$ is a category
over $St'/\LG$}

The goal of this section is to construct 
an action of the tensor category $D_{perf}^{\LG}(\Sth)$ on $\Dhat$
and of $D_{perf}^\LG(\St')$ on $D_{I^0I}$, both categories are equipped
with the tensor structure coming from tensor product of perfect complexes.

\subsubsection{The action of the tensor categories $Coh_{fr}^{\LG\times \LT}(\Cb_{St})$,
$Coh_{fr}^{\LG\times \LT}(\Cb_{St'})$}

\ 

We let $\Cb_{St}$ be the preimage of diagonal under the map $\Cb_{\gt}\times \Cb_{\gt}\to
\Lg \times \Lg$, and let $\Cb_{St'}$ be the preimage of $0$ under the second projection to 
$\Lt$.  We have open subsets $C_{St}\subset \Cb_{St}$ and $C_{St'}\subset \Cb_{St'}$
where the action of $\LT\times \LT $ is free and $St=C_{St}/\LT^2$, $St'=C_{St'}/\LT^2$.

Notation $Coh_{fr}^{\LG\times \LT}(\Cb_{\gt})$,
 was introduced in \S \ref{deeqsubs}, tensor categories $Coh_{fr}^{\LG\times \LT^2}(\Cb_{St})$,
$Coh_{fr}^{\LG\times \LT^2}(\Cb_{St'})$ etc. are similarly defined as full subcategories
in the categories of equivariant coherent sheaves whose objects are obtained from the structure sheaf by tensoring with a representation.

We apply Corollary \ref{extefun} in the following setting: the group $\LG$ is replaced
by $\LG^2$ and $\C$ is the category of functors $\Dhat\to \Dhat$
(respectively $D_{I^0I}\to D_{I^0I}$). 

We have two actions of $Rep(\LT)$ coming from, respectively, left and right convolution
with Wakimoto sheaves. We consider the action of $Rep(\LG^2)$ obtained as composition
of restriction to the diagonal copy of $\LG$ and the action by central functors. 
The nearby cycles monodromy acting on the cental functor defines a tensor endomorphism 
$E$ of the $\LG^2$ action, while the torus monodromy defines an action of $\Lt^2$.
It is not hard to see that conditions of Corollary \ref{extefun}
are satisfied. Thus we get an action of $Coh^{\LG^2\times \LT^2}_{fr}(\Cb_{\gt}^2)$ on $\Dhat$, $D_{I^0I}$.

The fact that the action of $\LG^2$ factors through restriction to diagonal is 
easily seen to imply that the action factors canonically through a uniquely defined
action of $Coh^{\LG\times \LT^2}(\Cb_{\gt}^2)$. Furthermore, since the isomorphism between
the two actions of $\LG$ is compatible with the tensor endomorphism $E$, both actions
factor through a uniquely defined action of $Coh_{fr}^{\LG\times \LT^2}(\Cb_{St})$.
Finally, since the second (right monodromy) action of $\Lt$ on $D_{I^0I}$ vanishes,
the action of  $Coh_{fr}^{\LG\times \LT^2}(\Cb_{St})$ factors through $Coh_{fr}^{\LG\times \LT^2}(\Cb_{St'})$. We denote the two actions by $\Phi_{fr}$, $\Phi_{fr}'$ respectively.

\subsubsection{Extending the actions to the perfect derived categories}
\label{extecompl}
Our next goal is to extend the action described in the previous subsection
to complexes. We encounter the standard non-functoriality of cone issue, which we
circumvent in the following way. 

We use the equivalences $Ho(\Th)\iso \Dhat$, $Ho(\TT)\iso D_{I^0I}$.

Assume given a finite complex $\F^\bu$ of objects in $Coh^{\LG\times \LT^2}(\Cb_{St})$,
where 
each term $\F^i$ is a trivial vector bundle twisted by a representation $U^i$ of $\LG\times \LT^2$.
Pick $\la, \,\mu\in \La$ so that for each character $(\la_j,\mu_j)$ of $\LT^2$ appearing in one of the representations $U^i$ we have $\la+\la_j\in (-\La^+)$, $\mu+\mu_j\in \La^+$. 

In view of Corollary \ref{prescost} and Lemma \ref{DelDel}(d) the functor $\Phi_{fr}(\F^i)
\circ \J_\la^l\circ \J_\mu^r$ sends $\Th$ to $\Phat$,
where $\J_\la^l: X\mapsto \J_\la*X$, $\J_\mu^r:X\mapsto X*\J_\mu$; thus
one gets a functor $Ho(\Th)\to Ho(\Phat)$ sending a complex $T^\bu$
to the total complex of the bicomplex $\Phi_{fr}(\F^\bu)\circ
\J_\la^l\circ \J_\mu^r(T^\bu)$.

We now define a functor $\Dhat \to \Dhat$ as the composition:
$$\begin{CD}
\Dhat @>{\J_{-\la}^l\circ \J_{-\mu}^r}>> \Dhat @<{\sim}<< Ho(\Th)@>{\Phi_{fr}(\F^\bu)\circ
\J_\la^l\circ \J_\mu^r}>> Ho(\Phat)\to \Dhat.
\end{CD}$$

We claim that different choices of $\la,\,\mu$ produce canonically isomorphic
functors. This follows from existence of a canonical up to homotopy quasi-isomorphism 
$\J_{-\la}*T\to T'$, $T'\to T*\J_\mu$, where $T$, $T'$ are finite complexes of objects in 
$\Th$ and $T'$ representing the object in the derived category corresponding to $\J_{-\la}*T$
(respectively, $T*\J_\mu$), $\la,\mu \in \La^+$.

Thus
we get a well defined functor $Ho(Coh_{fr}^{\LG\times \LT^2}(\Cb_{St}))\to
End(\Dhat)$. 
 It is not hard to see from the definition that the last arrow carries
a natural monoidal structure. 

Let $Acycl_{St}\subset Ho(Coh_{fr}^{\LG\times \LT^2}(\Cb_{St}))$
be the subcategory of complexes whose restriction to $C_{St}$
is acyclic.  
As in Proposition \ref{Acyclto0}, the fact that the lowest weight arrow
$\hw_\la$ extends to a filtration by Wakimoto sheaves compatible with convolution 
implies that $Acycl_{St}$ acts on $\Dhat$ by zero. 
In view of \eqref{Idem} we have  
$$Idem(Ho(Coh_{fr}^{\LG\times \LT^2}(\Cb_{St}))/Acycl)\cong D_{perf}^\LG(St).$$
Thus we obtain an action 
of $ D_{perf}^\LG(St)$ on $\Dhat$.
Finally, since the action of the log monodromy endomorphism is pro-nilpotent,
we conclude that the action factors through $D_{perf}^\LG(\Sth)$.

A parallel argument (with the last sentence omitted) endows $D_{I^0I}$
with an action of $D_{perf}^\LG(St')$.


\subsubsection{Compatibility between the two actions}
For future reference we record a compatibility between the two actions.

\begin{Lem}\label{compaact}
For $\F\in D_{perf}^\LG(\Sth)$, $X\in \Dhat$ and $Y\in D_{I^0I}$ we have
 canonical isomorphisms
$$\pi_*(\F(X))\cong i_{St}^*(\F)(\pi_*(X)),$$
$$\pi^*(i_{St}^*(\F)(Y))\cong \F(\pi^*(Y)),$$
where $i_{St}$ denotes the closed embedding $St'\to \Sth$. The isomorphism
is functorial in $\F$, $X$, $Y$ it is also compatible with the monoidal
structure of the action functor.
\end{Lem}

\proof Comparing the procedures of extending the action to the category
of complexes for $\Dhat$ and $D_{I^0I}$ and using that
 $\pi_*$ sends $\Th$ into $\TT$ we see that to get the first isomorphism
 it suffices to 
construct a functorial isomorphism for $\F\in Coh_{fr}^\LG(\Sth)$.
This follows from $\pi_*(\Zb_V)=Z_V$, $\pi_*(\J_\la)=J_\la$, where the second
isomorphism is compatible with the log monodromy endomorphism and the last two isomorphisms
are compatible with the lowest weight arrows. 
The second isomorphism can be deduced using the adjunction
$$Hom(\F(X),X')\cong Hom(X,\F^*(X'))$$
which holds for both actions; here $\F^*=R\uHom(\F,\O)$. 
In view of the isomorphism $i^*(\F^*)\cong (i^*(\F))^*$
 the second isomorphism follows from the first one.
 \epf

\section{The anti-spherical projector}\label{sectionXi}

Set $\Xib=\That_{w_0}$, $\Xis=T_{w_0}$.

Recall that $D_{IW}^{I^0}$ is the derived category of Iwahori-Whittaker
sheaves on $\Fltil$. We have {\em averaging functors}  
$Av^{IW}: D_{I^0I^0}\to D_{IW}^{I^0}$ and 
$Av^{I^0}_{IW}:D_{IW}^{I^0}\to D_{I^0I^0}$.

\subsection{$\Xib$ and Whittaker averaging}
\begin{Prop}\label{onXi} 
a) Right convolution with $\Xib$ is isomorphic to $Av^{I^0}_{IW}\circ Av^{IW}$.

b)  Convolution with $\Xib$ is isomorphic to its left and right adjoint.

c) The full subcategory in $\Th$ consisting of direct sums of copies of $\Xib$
is a  subcategory closed under the convolution product. It is tensor
 equivalent to the full subcategory in
$Coh(\widehat{\t^*\times_{\t^*/W_f} \t^*})$ whose objects are sheaves isomorphic
to $\O^{\oplus N}$ for some $N$; 
here "hat" stands for completion at zero.

d) Consider the full subcategory $\P_{I^0I^0}^{fin}\subset \P_{I^0I^0}$ of sheaves
supported on $G/U\subset \Fltil$, and let $\overline{\P_{I^0I^0}^{fin}}$ be its Serre quotient
by the Serre subcategory generated by all irreducible objects except for $L_e$, the irreducible
object supported on the closed cell $\Flt_e$. 
Then the endofunctor of $\overline{\P_{I^0I^0}^{fin}}$  induced by the functor 
$\F\mapsto \F*\Xib$ is isomorphic to the functor $\O(\Lt)\otimes _{\O(\Lt)^{W_f}}\F$, where 
$\O(\Lt)^{W_f}\subset \O(\Lt)$ acts on $\overline{\P_{I^0I^0}^{fin}}$ by log monodromy with respect
to the right $T$ action. 
\end{Prop}

\proof Part a) follows from \cite[Lemma 4.4.11(3)]{BY}. To check
part (b)  we use part (a) and adjunctions
in \cite[Lemma 4.4.5]{BY}, which show that the right adjoint to the functor
$Av^{I^0}_{IW}\circ Av^{IW}$ is isomorphic to $(Av^{I^0}_{IW})_*\circ Av^{IW}$;
here $(Av^{I^0}_{IW})_*$ is the right adjoint to the pull-back functor 
defined using the $*$ direct image. The isomorphism 
$Av^{I^0}_{IW}\circ Av^{IW}\cong (Av^{I^0}_{IW})_*\circ Av^{IW}$
follows from the relation between the two convolutions on $D_{I^0I^0}$
explained in footnote \ref{footconv} and \cite[Corollary 5.4.3]{BY}. 

Part (c) is a consequence of \cite[Proposition 4.7.3]{BY}.

Part (d) follows from parts (b,c) since $\Xib$ is the projective cover of $L_e$ in
$\P_{I^0I^0}^{fin}$.

{\bf Remark.} Beilinson-Bernstein Localization Theorem identifies $\P_{I^0I^0}^{fin}$
with category $O$ of modules over the Lie algebra $\g$ with a regular
integral generalized central character $\la$. It is easy to see that convolution 
with $\Xib$ is a projective functor isomorphic to the composition of two translation 
functors: translation  from $\la$ to the singular central character $-\rho$
followed by translation from $-\rho$ to $\la$. Properties (b,d) of this functor are well known
and play a central role in Soergel bimodules method. 


\subsection{Tilting property of $\Xib*Z_\la$}

\begin{Prop}
\label{XiZT}  For $V\in Rep(\LG)$ we have

a) $\Xis*Z_V\in \TT$.

b) $\Xib*\Zb_V\in \Th$.
\end{Prop}
\proof Recall that $\P_{II}\subset D_{II}$ is the category of perverse sheaves, let
$^f\P_{II}$ be the Serre quotient of $\P_{II}$ by the Serre subcategory generated
by  irreducible objects with support $\overline{\Fl_w}$, where $w$ is not the  minimal
length element in its coset $W_fw$. Let $\bar{Z_V}$ be the image of $Z_V$ in $^f\P_{II}$.
It follows from \cite[Theorem 7]{AB} together with \cite[Theorem 2]{AB} that $\bar{Z_V}$ admits 
a standard and a costandard filtration. Here a filtration is 
called (co)standard if its associated graded is a sum of
$\overline{j_{w!}}$, (respectively, $\overline{j_{w*}}$), where $\overline{j_{w!}}$,
$\overline{j_{w*}}$ is the image
of $j_{w!}$, respectively $j_{w*}$, under the projection $\P_{II}\to ^f\P_{II}$. 
It is easy to see (either by combining \cite[Lemma 4(a)]{AB} with a 
"left--right swap" of Proposition \ref{onXi}(a), or directly) that the functor
$\P_{II}\to \P_{I^0I}$, $\F\mapsto  \Xi* \F$ factors through $^f\P_{II}$.
It follows that  $\Xis*Z_V$ admits
a filtration whose subquotients are of the form $\Xis*j_{w!}$ and another
one with subquotients  of the form $\Xis*j_{w*}$. 
It is also easy to see that  $\Xis*j_{w!}$ carries a filtration such that
$gr(\Xis*j_{w!})\cong \oplusl_{v\in W_f} j_{vw!}$ and similarly for
 $\Xis*j_{w*}$. This proves part (a).

Part (b)
follows from Proposition \ref{dirim}, compatibility of central functors with direct image 
and 
part (a) of this Proposition. \epf

{\bf Remark.} The proof of Proposition \ref{XiZT} is the only place in this article where we use
the results of \cite{AB} directly, without applying the "left--right swap". 

\begin{Cor}\label{noExtZT}  For $T\in \TT$, $\That\in \Th$ we have 
$$Ext^{\ne 0}(\J_\la*\Xi *J_\mu*Z_\nu, T) =0,$$
$$Ext^{\ne 0}(\J_\la*\Xib *\J_\mu*\Zb_\nu, \That) =0$$ provided $(-\la),(-\mu) \in \Lambda^+$,
i.e.  $\la$, $\mu$ are anti-dominant.
\end{Cor}

\proof We have $$Ext^{\bu}(\J_\la*\Xib *\J_\mu*\Zb_\nu, \That) \cong Ext^{\bu}(\J_\la*\Xib *\Zb_\nu, \That*\J_{-\mu}).$$
Comparing Proposition \ref{XiZT} with Corollary \ref{prescost} we see
that  $\J_\la*\Xib *\Zb_\nu$ admits a free monodromic standard filtration,
while $\That*\J_{-\mu}$ admits a free-monodromic costandard filtration, which implies
the second vanishing. The first one is similar. \epf

\subsection{Convolution with $\Xib$ and the Springer map}

We let  $p_{Spr}$ denotes the projection $\gt\to\Lg$. 

\begin{Prop}\label{Xiaction} The equivalence  \eqref{PhiIWmon} 
intertwines the endo-functor $\F\mapsto \Xib*\F$ with the endo-functor $p_{Spr}^*p_{Spr*}$.
\end{Prop}


We start with 

\begin{Lem}\label{PhiXi}
 Recall that $\Phib_{IW}$ denotes  the equivalence $D^b(Coh^{\LG}(\gth))\cong \Dhat_{IW} $.

a) The object $(\Phib_{IW}^{-1}(Av^{IW}(\Xib))$ is canonically isomorphic to $\O(\Lt)\otimes _{\O(\Lt)^{W_f}}\O$. 

b) The composed  functor $\F\mapsto \Xib*(\Phib_{IW} \circ p_{Spr}^*(\F))$ is isomorphic to the functor
$\F\mapsto \Phib_{IW}(\O(\Lt)\otimes _{\O(\Lt/W_f)}p_{Spr}^*)$.

c) For $\F\in \Dhat_{IW}$ we have $\Xib*\F=0$ iff $(\Phib_{IW})^{-1}(\F)\in Ker(p_{Spr*})$.
\end{Lem}

{\em Proof} of the Lemma. a) 
The restriction of the functor $Av^{IW}$ to the category $\P_{I^0I^0}^{fin}$  factors through the category  $\overline{\P_{I^0I^0}^{fin}}$ (notations of Proposition \ref{onXi}).
Thus Proposition \ref{onXi}(d) shows that $Av^{IW}(\Xib)\cong \O(\Lt)\otimes _{\O(\Lt)^W}
Av^{IW}(\Delta_e)$. Since $ Av^{IW}(\Delta_e)\cong \Phib_{IW}(\O_{\gt})$, the claim follows.


b) The functor $\Phib_{IW}\circ p_{Spr}^*:D^b(Coh^\LG(\hat{\Lg}))\to \Dhat_{IW}$ comes from the {\em central} action of $D^b(Coh^\LG(\hat{\Lg}))$ on $\Dhat$. Since this action commutes with the functor of convolution with $\Xib$, (b) follows from (a).

c) The kernel of $p_{Spr*}$ is the (right) orthogonal to the objects
 $\O\otimes V$,
$V\in Rep(\LG)$. So we need to show that  $\Xib*\F=0 \iff 
Hom_{\Dhat_{IW}}(Av^{IW}(\Zb_\la), \F)=0$ for all $\la\in \La^+$.
First, if $\Xib*\F=0$ then by self-adjointness of convolution with 
$\Xib$,  $Hom(\Xib*Av^{IW}(\Zb_\la), \F)=0$. We have $\Xib*\Zb_\la\cong
 \Zb_\la*\Xib$ and $Av^{IW}( \Zb_\la*\Xib)$
admits a filtration where each subquotient is isomorphic to 
$Av^{IW}(\Zb_\lambda)$. By a standard argument (see e.g. \cite[Lemma 5]{BNcone}) it follows 
that $Hom_{\Dhat_{IW}}(Av^{IW}(\Zb_\la),\F)=0$. 
Conversely, suppose that $\Xib*\F\ne 0$. We need
to show that $Hom _{\Dhat_{IW}}(Av^{IW}(\Zb_\la),\F)\ne 0$ for some $\la$.
 Without loss of generality we can assume that $\F\in \Phat_{IW}$
 (recall that convolution with $\Xib$ is exact).
 Then, since
 $Hom_{\Dhat_{IW}}(\nabla_w^{IW},\Xib*\F)$ 
depends only on the 2-sided coset $W_f wW_f$, we see that
 $Hom_{\Dhat_{IW}}(\nabla_w^{IW},\Xib*\F)\ne 0$  for some $w$ 
which is maximal in its 2-sided $W_f$-coset.
 Using the tilting property of $\Xib*\Zb_\la$ one sees that for such $w$ the 
object
$\nabla_w^{IW}$ is a quotient of $\Xib*Av^{IW}(\Zb_\la)$ if $\la \in W_f w W_f$.
 Thus $Hom (\Xib*Av^{IW}(\Zb_\la),\Xib*\F)\ne 0$, hence 
$Hom(Av^{IW}(\Zb_\la*\Xib),\F)\ne 0$ and $Hom(Av^{IW}(\Zb_\la), \F)\ne 0$.  \epf

{\em Proof } of Proposition \ref{Xiaction}.
Set $F_\Xib:\F\to \Phib_{IW}^{-1}(\Xib*\Phib_{IW}(\F))$. Our goal is to show
that $F_\Xib\cong F_{Spr}$, where $F_{Spr}:= p_{Spr}^*p_{Spr*}$. Notice that both functors
are self-adjoint: for $F_{\Xib}$ this is Proposition \ref{onXi}(b) and for $F_{Spr}$ this follows from 
the fact that both $\gt$ and $\Lg$ have trivial canonical bundle and their dimensions coincide,
which yields an isomorphism $p_{Spr}^*\cong p_{Spr}^!$.

Lemma \ref{PhiXi}(b) shows that
$F_\Xib\circ p_{Spr}^*\cong F_{Spr}\circ p_{Spr}^*$ (notice that $p_{Spr*} p_{Spr}^*(\F)
\cong \O(\Lt)\otimes _{\O(\Lt)^{W_f}} \F$ canonically),
which implies $F_\Xib\circ F_{Spr}\cong F_{Spr}\circ F_{Spr}$. Self-adjointness of 
$F_{Spr}$ yields the adjunction arrow $F_{Spr}\circ F_{Spr}\to Id$, 
thus we get an arrow  $F_\Xib\circ F_{Spr}
\to Id$. Applying self-adjointness of $F_\Xib$ we get an arrow $c:F_{Spr}\to F_\Xib$. 

Lemma \ref{PhiXi}(b) 
provides an isomorphism $F_\Xib \circ p_{Spr}^*\cong F_{Spr}\circ p_{Spr}^*$,
a diagram chase  shows that this isomorphism coincides with the arrow induced by $c$.
Also, Lemma \ref{PhiXi}(c) shows that $F_\Xib|_{Ker(F_{Spr})}=0$. Thus $c_\F:F_{Spr}(\F)\to 
F_\Xib(\F)$ is an isomorphism when $\F\in Im(F_{Spr})$ or $\F\in Ker(F_{Spr})$. Again using 
self-adjointness of $F_{Spr}$, $F_\Xib$ we see that for any $\F$ the object
$Cone(F_{Spr}(\F)\overset{c}{\To} F_\Xib(\F))$ lies in the left  orthogonal to both $Ker(F_{Spr})$ and
$Im(F_{Spr})$. However, $^\perp Im(F_{Spr}) = Ker(F_{Spr})$ due to self-adjointness
of $F_{Spr}$;  thus 
 $^\perp Ker(F_{Spr}) \cap ^\perp Im(F_{Spr}) = 0$,
 which shows that $F_{Spr}(\F)\, \underset{c}{\iso} \, F_\Xib(\F)$
   for all $\F$. \epf

\section{Properties of $\Phip$}\label{Phip}

 Recall the actions defined in subsection \ref{extecompl}
  and objects
  $\Xib= \hat{T}_{w_0}$, $\Xis=T_{w_0}$. We define $\Phipb: D_{perf}^{\LG}(\Sth)\to \Dhat$, $\Phipb(\F)=\F(\Xib)$ and \\
  $\Phips:  D_{perf}^{\LG}(St')\to D_{I^0I}$, 
$\Phips(\F)=\F(\Xis)$.


\subsection{Compatibility of $\Phip$ with projection $St\to \gt$}
We start by recording some of the compatibilities following directly from the definitions.

\begin{Lem}\label{YYY} The following diagrams commute up to a natural 
isomorphism:
  $$\begin{CD} 
D_{perf}^\LG(\Sth) @>{i^*}>> D_{perf}^\LG(St')\\
@V{\Phipb}VV @VV{\Phip}V \\
\Dhat @>{\pi_*}>> D
\end{CD}$$

 $$\begin{CD} 
D^b(Coh^\LG(\gth)) @>{pr_{Spr,1}^*}>> D_{perf}^\LG(\Sth)\\
@V{\Phib_{IW}}VV @VV{\Phipb}V \\
\Dhat_{IW} @>{Av_{I^0}^{right}}>> \Dhat
\end{CD}$$

 $$\begin{CD} 
D^b(Coh^\LG(\gth)) @>{(pr_{Spr,1}')^*}>> D_{perf}^\LG(St')\\
@V{\Phib_{IW}}VV @VV{\Phip}V \\
\Dhat_{IW} @>{\pi_*\circ Av_{I^0}^{right}}>> D
\end{CD}$$

\end{Lem}

\proof Commutativity of the first diagram follows from
the corresponding compatibility for action  (Lemma \ref{compaact})
and the isomorphism $\pi_*(\Xib)\cong \Xi$ (Proposition \ref{TTT}(b)).
To see commutativity of the second one observe that
 the functor $Av_{I^0}^{right}$ of averaging with respect to the right
action of $I^0$ commutes with convolution on the left. 
For $\F\in
D^b(Coh^\LG(\gth))$ the object 
 $pr_{Spr,1}^*(\F)\in D_{perf}^\LG(\Sth)$ acts on $\Dhat$ by the left convolution
with $\Phi_{diag}(\F)$, thus
the required 
commutativity  follows the isomorphism $Av_{I^0}^{right}(\Delta_e^{IW})\cong \Xib$, which
is  a consequence of Proposition \ref{onXi}(a). The third diagram is obtained
by concatenation of the first two: the left (respectively, right) vertical arrow in the third diagram coincides with the left (right) arrow in the second (respectively, first) one, the horizontal arrows are compositions of the corresponding horizontal arrows in the first two. 
\epf

The goal of this subsection is the following

\begin{Prop}\label{compa_act}
The functor $\Phib_{perf}$ is compatible with the convolution action of $D^b(Coh^\LG(\Sth))$
on $D^b(Coh^\LG(\gth))$ and the action of $\Dhat$ on $\Dhat_{IW}$; i.e. for 
$\F\in D_{perf}^\LG(\Sth)$, $\G\in D^b(Coh^\LG(\gth))$ we have an isomorphism
\begin{equation}\label{PhiFPhiG}
\Phib_{IW}(\F*\G)\cong \Phipb(\F)*\Phib_{IW}(\G)
\end{equation}
functorial in $\F$, $\G$.
\end{Prop}
\proof
When $\F\cong \O$, so that $\Phib(\F)\cong \Xib$,
 the isomorphism (for any $\G$) is provided by Proposition \ref{Xiaction}. Since the functors commute with twist either by a line bundle or by a representation of $\LG$ we get an isomorphism for $\F$ of the form  $\O(\la,\mu)\otimes V$, $V\in Rep(\LG)$, this isomorphism is functorial in $\F$, $\G$.

By a standard argument\footnote{More generally, for a reductive group $H$ acting linearly on $\bbA^{N+1}$ 
and an $H$-invariant  locally closed subscheme
$X\subset {\mathbb P}^N$   every object $\F$ in the perfect equivariant
derived category $D_{perf}^H(X)$
 is a direct summand in an object represented by a finite complex of equivariant bundles of the form $\oplus V_i\otimes \O_X(n_i)$, $V\in Rep(H)$. To see this one constructs a bounded above complex $\F^\bu$ 
whose terms are finite sums of bundles
$V_i\otimes \O_X(n_i)$ representing $\F$, then denoting by $\F_{\geq -N}$
the "stupid" truncation of $\F^\bu$ we get for $N\gg 0$ a distinguished triangle
\eqref{truntrun}. For large $N$ we have $Ext^{N+1}(\F, \F_N)=0$, so the triangle splits.}
 (attributed, in particular, to Kontsevich, see also 
\cite{Nee}, Theorem 2.1 and Example 1.10)  
 any object in $D_{perf}^\LG(\Sth)$ is a direct summand 
in one represented by a finite complex of sheaves
of the form $\O(\la_i,\mu_i)\otimes V_i$, where $\la_i,\, \mu_i$ are antidominant,
thus we can assume without loss of generality
 that $\F$ is of this form.
 Pick $\nu\in \La^+$ such that $\mu_i+\nu \in \La^+$ for all $i$.  We can choose a finite
 complex of free-monodromic tilting objects in $\Phat_{IW}$ representing
 $\J_{-\nu}*\Phib_{IW} (\G)$, then 
  $\Phib(\G)$ is represented by a finite complex of objects $\J_\nu*\That_j$, 
  where $\That_j\in \Phat_{IW}$ is free-monodromic tilting. 

We claim that
 $$(\J_{\la_i}*\Xib*\J_{\mu_i})*(\J_{\nu}*\That_j)\in \Phat_{IW},$$
$$\Phib_{IW}^{-1}\left( (\J_{\la_i}*\Xib*\J_{\mu_i})*(\J_{\nu}*\That_j) \right) \in Coh^\LG(\gth) .$$

Here the first claim follows from Lemma \ref{DelDel}(d)
the second one follows from Lemma \ref{XiJT}(a) below.

Now \eqref{PhiFPhiG} follows by comparing Proposition \ref{tiltprop}(c) to Corollary 
\ref{convcompl}(c) below. 
 \epf

\begin{Lem}\label{XiJT}
a) $(\Phib_{IW})^{-1}(\Xib*\J_\mu*\That)\in Coh^\LG(\gth)$ for $\mu\in \La^+$ and
 $\That\in \Phat_{IW}$ a free-monodromic tilting object.

b) $(\Phi_{IW}^I)^{-1}(j_{w*}^{IW})\in Coh^\LG(\Nt)$ for any $w\in W/W_f$.

c)  $(\Phib_{IW})^{-1}(\nab_{w}^{IW})\in Coh^\LG(\gth)$ for any $w\in W/W_f$.

d) $\Phi_{IW}^{-1}: \P _{IW}\to D^{\geq 0}(Coh^\LG(\Nt))\cap D^{\leq \dim(\Nt)} (Coh^\LG(\Nt))$;

$\Phib_{IW}^{-1}: \Phat_{IW}\to D^{\geq 0}(Coh^\LG(\Nt))\cap D^{\leq \dim(\gt)} (Coh^\LG(\gth))$.

\end{Lem}

\proof a) An object  $\F\in D^b(Coh^\LG(\gth))$ lies in the abelian heart iff
for large $\la$ we have $R^i\Gamma(\F\otimes \O(\lambda))=0$ for $i\ne 0$.
Since $R^i\Gamma(\F)=Hom_{deeq}^\LG(\O,\F)$, 
 it suffices to show that $Hom_{\Phat_{IW}}^i(\J_{-\nu}*Av^{IW}(\Zb_\la), \Xib*\J_\mu*\That)=0$ for $i\ne 0$ and $\la,\mu,\nu\in \La^+$.
Using Proposition \ref{XiZT}(b) and  Proposition \ref{TTT}(c) we see that
$\J_{-\nu}*Av^{IW}(\Zb_\la)$ has a  free-monodromic standard filtration,
while $\J_\mu*\That$ and hence  $\Xib*\J_\mu*\That$
 has a free-monodromic costandard filtration, this implies the desired vanishing. 

Similarly, the first statement in (b) follows from $Ext^i_{D_{IW}^I}(^I Av^{IW}(J_{-\la}*Z_\mu), j_{w*}^{IW})=0$ for $i\ne 0$,
$\la\in \La^+$.
The latter $Ext$ vanishing is clear from the fact that $^I Av^{IW}(Z_\mu)$ is tilting in $\P_{IW}^I$
\cite[Theorem 7]{AB},
hence $^I Av^{IW}(J_{-\la}*Z_\mu)$ admits a costandard filtration.
The proof of (c) is parallel to that of (b), with (co)standard
replaced by  free monodromic (co)standard.
The inclusion $\Phi_{IW}^{-1}(\P _{IW})\subset D^{\geq 0}(Coh^\LG(\Nt))$ follows from part (b).
To check the other inclusion we use that $\Phi_{IW}(\O_\Nt(\la)\otimes V)\in \P_{IW}$ for all 
$\la\in \La$, $V\in Rep(\LG)$ thus
for $\F\in \P_{IW}$ we have $Ext^{<0}_{Coh(\Nt)} (\Phi_{IW}^{-1}(\F), \O_\Nt(\la))=0$. Applying  Grothendieck-Serre duality we conclude that $Ext^{<0}_{Coh(\Nt)}(\O(\la), S(\Phi_{IW}^{-1}(\F)))=0$,
where $S(\G)=R\underline{Hom}(\G,\O)$. Thus $S(\Phi_{IW}^{-1}(\F))\in D^{\geq 0}(Coh^\LG(\Nt))$,
so $\Phi_{IW}^{-1}(\F)=S(S(\Phi_{IW}^{-1}(\F)))\in D^{\leq \dim \Nt}$. This proves the first formula in 
(d), the second one is checked in a similar way.
\epf

\subsection{The functors $\Phip$, $\Phipb$ are fully faithful}
In this subsection we establish full faithfulness of $\Phip$, $\Phipb$. Since
$D^b(Coh_{\St'}^\LG(St)$ is a full subcategory in $D^b(Coh^\LG(\Sth))$, while
$D_{I^0I^0}$ is a full subcategory in $\Dhat$, it is enough to do so for $\Phipb$ only.

It suffices to show that the map
$$Hom^\bu (V\otimes \O_{\Sth}(\la,\mu), V'\otimes \O_{\Sth}(\la',\mu'))\to
Hom^\bu(\Zb_V*\JJ_\la*\Xib*\JJ_\mu, \Zb_{V'}*\JJ_{\la'}*\Xib*\JJ_{\mu'})$$
induced by $\Phipb$ is an isomorphism.

 The functor $\Phipb$ sends twisting by a line bundle to convolution by Wakimoto sheaves, and twisting by a representation of $\LG$ to the central functor.
  Since adjoint to such a twist is twist by the dual representation, and similar adjunction holds for the central functors and convolution by Wakimoto sheaves, we see that
it suffices to consider the case when $\la=0=\mu'$ and $V$ is trivial.

Then we have:
\begin{multline*}Hom_{\Dhat}(\Xib*\JJ_\mu, \JJ_{\la'}*\Zb_{V'}*\Xib)
\cong Hom_{\Dhat_{IW}}(Av^{IW}(\Xib * \JJ_\mu), Av^{IW}(\J_{\la'}*\Zb_{V'}))\cong \\
 Hom_{D^b(Coh^\LG(\gth))} (p_{Spr}^*p_{Spr*}(\O_{\gth}(\mu)), \O_{\gth}(\la')\otimes V')\cong Hom_{D^b(Coh^\LG(\gth))}(p_{Spr,2*}p_{Spr,1}^*(
 \O_{\gth}(\mu)), \\
 \O_{\gth}(\la')\otimes V').
\end{multline*}

Here the first isomorphism comes from the fact
 that right convolution with $\Xib$ is isomorphic to
$Av^{I^0}_{IW}\circ Av^{IW}$ (Proposition \ref{onXi}(a)). The second isomorphism uses the "coherent" description of the Iwahori-Whittaker category \eqref{PhiIWmon}
along with the fact that left convolution with $\Xib$
corresponds to $p_{Spr}^*p_{Spr*}$ on the coherent side (Proposition \ref{Xiaction}). 
Finally, the last isomorphism
comes from:   $p_{Spr}^*p_{Spr*}\cong p_{Spr,2*}p_{Spr,1}^*$, which
follows from base change for coherent sheaves and the fact that 
$Tor_{>0}^{\O(\Lg)}(\O_{\gt},\O_{\gt})=0$.

Using adjointness we get:
\begin{multline*}  Hom_{D^b(Coh^\LG(\gt))}(pr_{2*}pr_1^*(
 \O_{\gt}(\mu)),
 \O_{\gt}(\la')\otimes V')\cong Hom_{D^b(Coh(St))}(pr_1^*(
 \O_{\gt}(\mu)), \\ pr_2^*(\O_{\gt}(\la')\otimes V')),
 \end{multline*}
 where we used that $pr_2^*\cong pr_2^!$ since the target of $pr_2$ is smooth, while both its source and target have trivial dualizing complexes  
(more precisely, in both cases the dualizing complex is isomorphic to $\O[d]$, $d=\dim (\gt)=\dim(St)$).

Since $\Phipb: pr_1^*(
 \O_{\gt}(\mu)) \mapsto \J_\mu*\Xib$, \ $\Phipb: \O_{\gt}(\la')\otimes V'\mapsto \JJ_{\la'}
*\Zb_{V'}*\Xib$, 
we have constructed an isomorphism between the two $Hom$ spaces.
A routine diagram chase shows that this isomorphism coincides with the map induced by $\Phipb$.
\epf

\section{Extending an equivalence from the subcategory of perfect complexes} \label{gencrit}

\subsection{A criterion for representability}\label{subseccrit}
Let algebraic stack $X$ be given by $X=Z/H$ where $Z$ is a quasiprojective scheme over an algebraically closed field of characteristic
zero and $H$ is a reductive group. [The results of this section are likely valid in 
greater generality but we present the setting needed for our applications]. 
We fix a $H$-equivariant ample line bundle $L$ on $Z$, such a bundle
exists by Sumihiro embedding Theorem (though in examples considered in this paper
$Z$ comes equipped with a supply of such line bundles).

Set
$D=D^b(Coh(X))$ and let $D_{perf}(X)\subset D$ be the subcategory of perfect complexes.
Set  $D^{\leq n}_{perf}=D^{\leq n}(Coh(X))\cap D_{perf}(X)$, and let
$D^{\geq n}_{perf}\subset D_{perf}(X)$ be the full subcategory of objects represented
by complexes of locally free sheaves placed in degree $n$ and higher, and their direct summands.

\begin{Rem}\label{rem_finitist}
It is obvious that $D^{\geq n}_{perf}\subset D^{\geq n}(Coh(X))\cap D_{perf}(X)$. Using
\cite[Theorem 3.2.6]{RG} ("finiteness of finitistic dimension")
 one can also show that $D^{\geq n}_{perf}(X)
\supset D^{\geq n-\dim (Z)}\cap D_{perf}$. This implies that most of the statements
below hold with $D^{\geq n}_{perf}$ replaced by $D^{\geq n}(Coh(X))\cap D_{perf}(X)$. We neither prove nor use this point.
\end{Rem}

\begin{Prop}\label{fullfaithrep}
a) The natural functor from $D^b(Coh(X))$ to the category of contravariant
functors from
 $D_{perf}(X)^{op}$ to vector spaces is fully faithful.

\medskip

b) A cohomological functor $F$ from $D_{perf}(X)$ to vector spaces is represented by an object of $D^b(Coh(X))$ if and only
if the following conditions hold:

i) For any $n$ the functor $F|_{D^{\geq n}_{perf}}$ is represented by an object
of $D_{perf}(X)$ (not necessarily by an object of $D^{\geq n}_{perf}$).

ii) There exists $m$ such that $F|_{D^{\leq m}_{perf}} =0$.

\end{Prop}

\proof
 Fix $\F,\G\in D^b(Coh(X))$ and let $\phi_\F$, $\phi_\G$ be the corresponding functors
 on $D_{perf}(X)$.
  Fix a bounded above complex $\F^\bu$ of locally free sheaves representing $\F$.
  Let $\F_{\geq -n}=\tau^{\bete}_{\geq -n}(\F^\bu)$ denote  the stupid truncation.

Given a natural transformation $\phi_\F\to \phi_\G$ we get
morphisms $\F_{\geq -n}\to \G$, compatible with the arrows $\F_{\geq -n}\to \F_{\geq -(n+1)}$.
Choose $n$ such that $\F\in D^{>-n}(Coh(X))$. Then
for $N>n$ we have a canonical isomorphism $\F\cong \tau_{\geq -n}(\F_{\geq -N})$. Assuming
also that $\G\in D^{>-n}(Coh(X))$,
 we get an arrow $\F=\tau_{\geq -n}(\F_{\geq -N})\to \tau_{\geq -n}(\G)=
\G$. A standard argument shows that bounded above complexes representing a given
$\F\in D^b(Coh(X))$ form a filtered category (i.e. given two such complexes
$\F^\bu_1$, $\F^\bu_2$, there exists a complex $\F^\bu_0$ with maps of complexes
$\F^\bu_0\to \F^\bu_1$, $\F^\bu_0\to \F^\bu_2$ inducing identity maps in the derived category).
 This implies that the  arrow $\F\to \G$ does not depend on the choice of $\F^\bu$.

Thus we have constructed a map $Hom(\phi_\F,\phi_\G)\to Hom(\F,\G)$. It is clear from
the construction
that the composition $Hom(\F,\G)\to Hom(\phi_\F,\phi_\G)\to Hom(\F,\G)$ is the identity map.
It remains to see that the map  $Hom(\phi_\F,\phi_\G)\to Hom(\F,\G)$ is injective.

Let $h\in Hom(\phi_\F,\phi_\G)$ be a nonzero element.
Thus for some $\P\in D_{perf}(X)$ and $\varphi:\P\to \F$ we have $0\ne h(\varphi):\P\to
\G$. Fix again a complex $\F^\bu$, for $N$ as above we get a distinguished triangle
\begin{equation}\label{truntrun}
\F_N[N]\to \F_{\geq -N}\to \F\to \F_N[N+1]
\end{equation}
for some $\F_N\in Coh(X)$. For large $N$ we have $Hom(\P,\F_N[N+1])=0$, thus $\varphi$
factors through an arrow $\P\to \F_{\geq -N}$. It follows that for $N\gg 0$ applying $h$
to the tautological map $\F_{\geq -N}\to \F$ we get a nonzero arrow $\F_{\geq -N}\to \G$.
Since $Hom(\F_N[N],\G)=0=Hom(\F_N[N+1],\G)$ for large $N$, we see that the induced arrow
$\F\to \G$ is nonzero. This proves (a).

We now prove (b). We first check the "only if" direction.  Condition (ii) is clear, and to check
condition (i) let $\F$ be the representing object, and choose a bounded above complex $\F^\bu$ representing $\F$; we can and will choose $\F^\bu$ so that its terms are locally free sheaves. Setting again $\F_{\geq N}=\tau^{\bete}_{\geq N} (\F^\bu)\in D_{perf}(X)$, we claim
that $Hom (\G,\F)\iso Hom(\G,\F_{\geq N})$ when $\G\in D_{perf}^{\geq m}$, $N<m-d$, where $d=\dim (Z)$.
This follows from the fact that $Ext^i(\EE, {\mathcal K})=0$ for $i>d$,
 where $\EE,\, {\mathcal K}\in Coh(X)$ and $\EE$ is locally free.

To check the "if" direction,
given a functor $F$ satisfying the conditions
take $n$ in (i) satisfying $n< m-d$ where $m$ is as in (ii) and $d=\dim (Z)$.
Let $\F'\in D_{perf}(X)$ be a representing object for $F|_{D_{perf}^{\geq n}}$.
 We claim that $\F=\tau_{\geq n}(\F')$ represents $F$.

 First observe that
 \begin{equation}\label{Ftoright}
 \F\in D^{>m}(Coh(X)),
 \end{equation}
 to check this we need to see that $H^i(\F')=0$ for $i=n,\dots,m$. If $H^i(\F')\ne 0$ for such an $i$,
 we can find a locally free sheaf $\EE$ such that $Hom(\EE,H^i(\F'))\ne 0$ and $Ext^{>0}(\EE,H^j(\F'))=0$
 for all $j$ (in fact, we can take $\EE=L^{\otimes N}\otimes V$ where $L$ is an anti-ample
 $H$-equivariant line bundle on $Z$ and $V$ is a representation of $H$). Then we get
  $Hom(\EE[-i], \F')=F(\EE[-i])\ne 0$, which contradicts  (ii).

 We now construct a functorial isomorphism $F(\G)\cong Hom(\G,\F)$,
  $\G\in D_{perf}(X)$. Fix such $\G$, and fix  a finite complex $\G^\bu$
 of locally free sheaves representing $\G$.
 The desired isomorphism is obtained  as the following composition:
 $$Hom (\G,\F)\cong Hom(\tau^{\bete}_{\geq m}(\G^\bu), \F)\cong Hom(\tau^{\bete}_{\geq m}(\G^\bu), \F')
 \cong F(\tau^{\bete}_{\geq m}(\G^\bu))\cong F(\G).$$
 Here the first isomorphism follows from \eqref{Ftoright}, which implies that
 $Hom(\tau_{<m}^{\bete}(\G^\bu),\F)=0=Hom(\tau_{<m}^{\bete}(\G^\bu)[-1], \F)$.

 The second isomorphism follows
  from  the distinguished triangle
 $\tau_{<n}(\F')\to \F'\to \F \to \tau_{<n}(\F')[1]$ and the fact
 that
 $Hom(D_{perf}^{\geq m}(Coh(X)),
 D^{\leq n}(Coh(X)))=0$, since $m-n>d$ and $Ext^i(\EE, {\mathcal K})=0$ for $i>d$,
 where $\EE,\, {\mathcal K}\in Coh(X)$ and $\EE$ is locally free.

 The third isomorphism is the assumption on $\F'$, and the last isomorphism follows from
 (ii). It is easy to see that the constructed isomorphism is independent on the auxiliary
 choices and is functorial.
 \epf

\medskip

 Let $X=Z/H$ be as in the previous Proposition. We assume that $Z$ admits a
 projective $H$-equivariant morphism $Z\to Y$ where $Y$ is affine. Let $L$ be
 an $H$-equivariant ample line bundle on $Z$.  We have the homogeneous
 coordinate ring $\Ohat(Z)=\oplusl_{n\geq 0} \Gamma (L^{\otimes n})$. The assumptions
 on $Z$ imply that $\Ohat(Z)$ is Noetherian.


We now assume that $\CC$ is a triangulated category with a fixed full triangulated embedding $i:D_{perf}(X)\to \CC$.

For $M\in \CC$ we can form a module for the homogeneous coordinate ring
$$\Psit(M)=\oplusl_{n\geq 0,\la }Hom (i(L^{\otimes -n}\otimes \O(H)_\la), M),$$
where $\la$ runs over the set of dominant weights of $H$ and $\O(H)_\la$ denotes the corresponding
isotypic component of the translation  action of $H$ on $\O(H)$. 
A section of $L^{\otimes n}$ defines an element in $Hom_{Coh(X)}(L^{\otimes m}\otimes \O(H)_\la,
 \oplusl_\mu L^{\otimes m+n}\otimes \O(H)_\mu)$ for every $\la$ and $m$, thus
 $\Psit(M)$ does carry a natural action of the homogeneous coordinate ring.

Notice that
if $\CC$ is equipped with a $Rep(H)$ action making $i$ a functor of module categories
for $Rep(H)$ then we have:  $\Psit(M)=\oplusl_{n\geq 0 }Hom_{deeq}^H (i(L^{\otimes -n}), M).$

We also set: $\Psit_m(M)=\oplusl_{n\geq m,\la}Hom_{deeq}^H(i(L^{\otimes -n}\otimes \O(H)_\la), M)$.

\begin{Prop}\label{fingencrit} For $M\in \CC$ the following are equivalent.

a) For any $m$ the
  functor on $D_{perf}^{\geq m}(X)$,
 $\F\mapsto Hom(i(\F),M)$  is represented by an object of $D_{perf}(X)$.

b) The module $\Psit(M[n])$ is finitely generated for all $n$
and $\Psit(M[n])=0$ for $n\gg 0$.

c) We have $\Psit(M[n])=0$ for $n\gg 0$ and for any $n$ there exists $m$, such that  
$\Psit_m(M[n])$ is finitely generated.

\end{Prop}

The proof of the Proposition is based on the following
\begin{Lem}\label{PsiM}
If $\Psit(M[n])=0$ for $n\geq s$, then $Hom(i(\F),M)=0$ for $\F\in D_{perf}^{>s+d}$,
$d=\dim (Z)$.
\end{Lem}

\proof We claim that any object in $\F\in D_{perf}^{>s+d}$ is isomorphic to a direct summand
in an object represented by a complex placed in degree $s$ and higher, with each term isomorphic to
$L^{\otimes i}\otimes V$, $i\leq 0$, $V\in Rep(H)$. This clearly implies the Lemma.

It remains to check that claim. Let $\F\in  D_{perf}^{>s+d}$.
 By a standard argument there exists a bounded above complex $\F^\bu$
representing $\F$ whose terms are of the form $L^{\otimes n}\otimes V$, $n\leq 0$.
Then using the fact that $Ext^i$ from a locally free sheaf to any sheaf vanishes
for $i>d$, we conclude the argument by a standard trick:
consider the distinguished triangle
 $\F_s[s]\to \tau^{\bete}_{\geq -s}(\F^\bu)\to \F$ and use that $Hom(\F,\F_s[s+1])=0$. \epf

{\em Proof\ } of the Proposition. (a) clearly implies (b), while (b) implies (c).
 We proceed to prove that (c) implies (a).

Assume that (c) holds.
In view of the Lemma, it suffices to find for every $m$ an object $\F_{M,m}\in D_{perf}(X)$ and a morphism
$c_m:i(\F_{M,m})\to M$ so that $\Psit(Cone(c_m)[l])=0$ for $l\geq m$. Moreover, it suffices
to do so after possibly replacing the full embedding $i$ by the functor $i':\F\mapsto i\circ (\F\otimes L^{\otimes p})$ for some $p\in \Zet$ (notice that conclusion of Lemma \ref{PsiM} is not affected by such a
substitution).


Let $d_0$ be the largest integer such that $\Psit(M[d_0])\ne 0$.
We argue by descending induction in $d_0$.
Using the finite generation condition we find a locally free sheaf $\EE\in Coh(X)$
and a morphism $i(\EE)[-d_0]\to M$, such that the induced map $\Psit_m(i(\EE))\to
 \Psit_m(M[d_0])$ is
surjective for some $m\in \Zet$. Fix $m_0\geq 0$ such that $R^{>0}\Gamma (L^{\otimes i}\otimes \EE)=0$ for $i\geq m_0$.
We can assume without loss of generality that $m_0\geq m$.
 Then upon replacing the embedding $i$ by $i':\F\mapsto i(\F\otimes L^{\otimes -m_0})$ we get that 
  $M':=Cone\left(i( \EE)\to M \right)$ satisfies: $\Psit(M'[i])=0$
 for $i\geq d_0$. Also it is clear that the finite generation condition 
is satisfied for $M'$, $i'$. Thus we can assume that the statement is true for $M'$ by the induction assumption. Then
 the statement about $M$ follows from the octahedron axiom. 
 \epf

\subsection{A characterization of $D^b(Coh(X))$ as an ambient category of $D_{perf}(X)$}
We continue working under the assumption that
 $i $ 
is fully faithful.
Assume also that equivalent conditions of Proposition \ref{fingencrit} hold, thus
the condition of Proposition \ref{fullfaithrep}(b)(i) is satisfied.
Assume also that assumption (b,ii) holds.
In view of Proposition \ref{fullfaithrep} we get a functor $\Psi:\CC\to D^b(Coh(X))$ sending $M\in \CC$ to $\F\in D^b(Coh(X))$
representing the functor $\G\mapsto Hom(i(\G), M)$ on $D_{perf}(X)$.


It is not hard to see that $\Psi$ is a triangulated functor.

We now assume that $\CC$ is equipped with a bounded $t$-structure $\tau$. Consider the
 following properties of the functor $\Psi$ in relation to the $t$-structures.

A) The functor $\Psi$ is of bounded amplitude, i.e. there exists
$\del$ such that $\Psi:D^{\tau,\leq 0}\to D^{< \del}(Coh(X))$, $\Psi:D^{\tau,\geq 0}\to D^{> -\del}(Coh(X))$.

B) There exists $d\in \Zet$ such that for $\F\in \CC$ we have: 
$\Psi(\F)\in  D^{\leq 0}(Coh(X)) \Rightarrow \F\in \CC^{\tau,\leq d}$.


C)
There exists $\don>0$ such that 
for $\F\in \CC$ we have: 
$H^i(\Psi(\F))=0$ for $i\in [-\don,\don]\    \Rightarrow H^{\tau,0}(\F) =0$. Here $H^i(\Psi(\F))
\in Coh(X)$ is the cohomology with respect to the standard $t$-structure on $D^b(Coh(X))$.

\begin{Prop}\label{ABsidelinatrube}
 a) Property (B) implies that $\Psi$ is fully faithful.

b) Properties (A), (C) imply that $\Psi$ is an equivalence.

\end{Prop}

\proof a) To unburden notation we assume without loss of generality that $d=0$, this can be achieved
by replacing the  $t$-structure $\tau$ with its shift by $-d$.

Recall that $i$ is assumed to be a full embedding, which implies that
 $\Psi\circ i\cong Id_{D_{perf}(X)}$. It follows from the definition of $\Psi$ that
 for $\F\in D_{perf}(X)$ we have:
$$Hom(i(\F), M)\cong Hom (\F, \Psi(M))\cong Hom(\Psi i(\F),\Psi(M)).$$
Thus the map $Hom(M_1,M_2)\to Hom(\Psi(M_1), \Psi(M_2))$ is an isomorphism when $M_1\in Im(i)$.

Fix $M_1,M_2\in \CC$. Fix $n$ such that $M_2\in \CC^{\tau,>n}$ and $\Psi(M_2)\in D^{>n}(Coh(X))$.
Fix a bounded above complex $\F^\bu$ of locally free sheaves representing $\F=\Psi(M_1)$,
and let $\F_{\geq N}\in D_{perf}(X)$ be the naive truncation as above.
We have an exact triangle $\F_N[-N]\to \F_{\geq N}\to \F$ for some $\F_N\in Coh(X)$.

Assuming $N<n$, we get
$$Hom (\Psi(M_1),\Psi(M_2))\cong Hom(\F_{\geq N},\Psi(M_2))\cong Hom(i(\F_{\geq N}),M_2).$$
We have a morphism $i(\F_{\geq N})\to M_1$ whose cone lies in $D^{\tau, \leq N-1}$ in view of
condition (B). [Notice that $\Psi$ sends this cone to $\F_N[-N+1]$.]
Thus $Hom(M_1,M_2)\cong Hom (i(\F_{\geq N}),M_2)$, so composing the above
isomorphisms we get that $Hom(M_1,M_2)\cong Hom(\Psi(M_1),\Psi(M_2))$. It is easy to see that
this map coincides with the map induced by $\Psi$, so (a) is proved.

b) Property (C) implies (B), thus $\Psi$ is fully faithful by (a), 
 it remains to show that it is essentially surjective. Fix $\F\in D^b(Coh(X))$
and a bounded above complex of locally free sheaves $\F^\bu$ representing $\F$.
Let $n$ be such that $\F \in D^{\geq n}(Coh(X))$. Fix $N< n$ 
and let $A=i(\F_{\geq N})$. We assume as we may that $N<n-2\don$, then condition (C) implies
that  $H^{\tau,m}(A)=0$ for $m\in [N+\don+1,n-\don-1]$. Pick such an $m$ and set
$B=\tau_{\geq {m}}(A)$. We have an exact triangle $A\to B\to C$ where $A\in \CC^{\tau, \geq n-\don}$
and $C\in \CC^{\tau, < N+\don}$. Thus applying condition (A) we get
an exact triangle $\Psi(A)\to \Psi(B)\to \Psi(C)$, where $\Psi(B)\cong \F_{\geq N}$,
$\Psi(A)\in D^{\geq n-\don-\del}(Coh(X))$ and $\Psi(C)\in D^{<N+\don+\del}(Coh(X))$.
Assuming as we may that $N< n-2(\don+\del)$, we see that $\Psi(A)\cong \F$ which
proves that $\Psi$ is essentially surjective.
 \qed

\section{Compatibility between the $t$-structures and construction of
the functor from constructible to  coherent category}\label{with_t}

\subsection{Almost exactness of $\Phi_{perf}$}
\begin{Prop}\label{gdekogo} For some $d>0$ the following holds.
 If $\F\in D$ is such that $Hom^{i}( \Phip(\O(\la,\mu)\otimes V),\F)=0$ for all 
$\la,\mu\in \La$, $V\in Rep(\LG)$ and $i\in [-d,d]$, then $H^{p,0}(\F)=0\in \P$.


\end{Prop}


The proof of Proposition is preceded by some auxiliary results.

\begin{Lem}\label{obnu} For $X\in D_{I^0I}$
 there exists a finite subset $S\subset W$, such that
for  $\lambda\in \Lambda^+$  we have
\begin{equation}\label{obnueq}
j_w^!(X*J_{\lambda} )\ne 0\Rightarrow w\in S\cdot (\lambda)
\subset
 W.
\end{equation}
\end{Lem}

\proof
By the  $!$-support of an object $X\in D$ we mean the
set of points $i_x:\{x\}\imbed \Fl$ such that $i_x^!(X)\ne 0$.
Proper base change shows that the $!$-support of
$X*J_{\lambda} $  lies in the convolution of sets
$\supp(X)$ and $\Fl_{\lambda}$. This implies 
\eqref{obnueq}. \epf

\begin{Lem}\label{gdesemii}
Let $\F$ be as in Proposition \ref{gdekogo}.

 For large $\lambda$ and $n\in [-d+2\dim (\gt), d-2\dim (\gt)]$ we have
\begin{equation}\label{sloi}
 Ext^{n}(j_{w!}, \F*J_{\lambda})=0
\end{equation} for all $w$.


\end{Lem}

\proof  According to Lemma \ref{obnu} there exists a finite set $S\subset W$ such that
 for large  $\lambda$  the left
hand side of \eqref{sloi} vanishes for all $n$  unless $w\in
S\cdot  (\lambda)$. Also  for large $\lambda$ we have $S\cdot ( \lambda) 
\subset W_f\cdot (\Lambda^+)$ and each element in this set is the minimal length representative of its right $W_f$ coset. Hence for all $w\in W$
we have
\begin{equation}\label{sloi1}
Ext^{p}_D(j_{w!}, \F*J_{\la} )\cong
Ext^{p}_D(\Delta_w *  \Xi, \F* J_{\la}),
\end{equation}
$$ \mathrm{or}\ \ \ \ \ Ext^{p}_D(j_{w!}, \F*J_{\la} )=0,$$
which follows
from the fact that $\Delta_{w}*\Xi $ admits a filtration with associated graded
$\oplusl_{w_f\in W_f} j_{ w w_f !}$, and
 for $w_f\ne e$ we have $Ext^\bu(j_{w w_f !}, \F*J_{\lambda})=0$ provided
 that $Ext^\bu(j_{w !}, \F*J_{\lambda})\ne 0$.
 We can rewrite the right hand side of \eqref{sloi1} as
 $$ Ext^p_{\Dhat}(\Delta_w*\Xib,\pi^*(F*J_\la)[r])\cong Ext^{p}(Av_{IW}(\Delta_w),
 Av_{IW}(\pi^*(F*J_\la)[r]) ),$$
$r=\rank(G)$, where we used Proposition \ref{onXi}(a) and isomorphisms $\pi_*(\Xib)\cong
 \pi_!(\Xib)[r]\cong \Xi$, $\pi^!\cong \pi^*[2 r]$. 
 
 Now $\Phi_{IW}^{-1}(Av_{IW}(\Delta_w))\in D^{\geq 0}(Coh^\LG(\gth))
 \cap D^{\leq \dim \gt}(Coh^\LG(\gth))$ by Lemma \ref{XiJT}(d), while
 the condition of Proposition \ref{gdekogo} implies that $\Phi_{IW}^{-1}(\pi^*(\F*J_\la))[r]$
 is concentrated in homological degrees less than $-d$ and greater than $d$.
 Since $Coh^\LG(\gth)$ has homological dimension $\dim(\gt)$, we get the desired vanishing.
 \epf

{\em Proof}\  of Proposition \ref{gdekogo}. In the assumptions of part (a)
Lemma \ref{gdesemii}(a) implies
that for large $\lambda$ the object $\F*J_{\la}$ 
is concentrated  in homological degrees less than $-d+2\dim (\Lg)$ and greater than 
$d-2\dim (\Lg)$. 

We finish the proof by invoking a result of Lusztig \cite{cells} saying that Lusztig's $a$-function 
for the affine Weyl group is bounded by $\dim (G/B)$, thus convolution
of two object in $\P_{II}$ lies in perverse degrees from $-\dim(G/B)$ to $\dim(G/B)$. 
Thus $\F=(\F* J_{\la})*J_{-\la}$ has no cohomology in perverse degree zero provided
that $d> 2\dim (\Lg) + \dim (G/B)$. 
\epf

\subsection{The functor from constructible to coherent category.}
Applying the general construction of section \ref{subseccrit} (see notation introduced prior to Proposition \ref{fingencrit}) in the present situation: $X=St'/\LG$, $\CC=D_{I^0I}$, $L=\O(\la,\mu)$
for strictly dominant weights $\la$, $\mu$,
we get a 
 functor $\Psit$ from $D_{I^0I}$ to $\LG$-equivariant modules over the homogeneous
coordinate ring of $St'$.

\begin{Prop}\label{exaprop} 
For $\F\in Perv_N(G/B)\subset \P$ we  have $\Psit(\F*J_{\rho}[n])=0$ 
for $n\ne 0$.
\end{Prop}

\proof  It suffices to check that for $\F=j_{w!}$, $j_{w*}$, $w\in W_f$ we have
$\Psit(\F*J_\rho[n])=0$ for $n\ne 0$. This reduces to showing that for dominant $\la,\mu,\nu$ with 
 $\mu$ strictly dominant we have $Ext^i(\J_{-\la}*\Xi *Z_\nu*\J_{-\mu},\F)=0$
for $i\ne 0$.
We have $\ell(w\mu)=\ell(\mu)-\ell(w)$, $\ell(\la w)=\ell(\la)+\ell(w)$ for $w\in W_f$.
Thus for such $w$ we have 
$$Ext^i (\J_{-\la}*\Xi *Z_\nu*J_{-\mu},j_{w!})=Ext^i(\J_{-\la}*\Xi *Z_\nu,
j_{w!}j_{\mu*})=Ext^i(\J_{-\la}*\Xi *Z_\nu,
j_{w\mu*}),$$
$$Ext^i (\J_{-\la}*\Xi *Z_\nu* J_{-\mu},j_{w*})=Ext^i(\Xi *Z_\nu*J_{-\mu},
j_{\la*}j_{w*})=Ext^i(\Xi *Z_\nu*J_{-\mu},
j_{\la w*}). $$
Since $\Xi*Z_\nu$ is tilting, $\J_{-\la}*\Xi *Z_\nu$ admits a standard filtration, which shows
that the first $Ext$ group vanishes for $i\ne 0$. Likewise, $\Xi *Z_\nu*J_{-\mu}$ admits
a standard filtration which shows vanishing of the second $Ext$ group for $i\ne 0$. \epf

\begin{Prop}\label{fingenprop}
The module
$\Psit(\F)$ is finitely generated for any $\F\in D$.
\end{Prop}

\proof For $\F$ in the image of $\Phi_{perf}$ this is clear from 
the fact that $\Phi_{perf}$ is a full embedding. 
Every irreducible object in $Perv_N(G/B)$ is a subquotient of $\Xi$.
Then it follows from
the previous Proposition that if $\LL$ is such an irreducible object, $\Psit(\LL*J_\rho)$
is a subquotient of $\Psit(\Xi*J_\rho)$, hence it is finitely generated
(since the homogeneous coordinate ring of Steinberg variety is Noetherian), while
 $\Psit(\LL*J_\rho[n])=0$ for $n\ne 0$. 
 It follows that the 
same is true for any $\LL\in Perv_N(G/B)$.
Now it follows from Proposition \ref{fingencrit} that $\Psit(\J_\la*\F*J_\mu[n])$ is finitely
generated for $\F\in Perv_N(G/B)$ and any $\la$, $\mu\in \La$, $n\in \Zet$.
%
Such objects generate $D$, so the claim follows. \epf

\begin{Prop}\label{boundampl}
There exists $\del$, such that for all $\F\in \P_{I^0I^0}$ we have  
$$Hom_{deeq}^{\LG\times \LT^2}(\Xib, \F[i] )=0$$ for $i\not \in [-\del,
\del]$.
\end{Prop}

\proof We need to check that for some $\del\in \Zet$ we have 
$$Ext^{i}_{\Dhat}(\J_{-\la}*\Zb_\nu*\Xib *\J_{-\mu}, \F)=0,$$
for $i\not \in [-\del,\del]$, $\F\in \P$.
According to a result of Lusztig \cite{cells}, Lusztig's $a$-function
for an affine Weyl group is bounded by $\dim (G/B)$, which implies that the 
convolution of any two objects in $\P_{II}$ is concentrated in perverse degrees
from $-\dim(G/B)$ to $\dim(G/B)$. 
It follows that $\J_\la*\F*\J_\mu\in D^{\geq - 2\dim(G/B)}(\P_{I^0I^0})\cap
D^{\leq 2\dim(G/B)}(\P_{I^0I^0})$.

Thus  it suffices  to show that
for some $\del>0$ we have 
$$
Ext^i_D(\Z_\nu*\Xib, \F)=0\ \ \ \ \ \ \mathrm{for}\ \ \  i\not \in [-\del,\del],\ \ \ \F\in \P_{I^0I^0}. 
$$
Using Proposition \ref{onXi}(a) we can rewrite the right hand side\footnote{In fact,
the main result of \cite{B} (see \cite[Theorem 2]{B}, cf.  a related  statement 
Theorem \ref{suppf}(b) below) 
shows that for an irreducible $\F$ this $Ext$ group is an isotypic component in the cohomology of a coherent IC sheaf
on the nilpotent cone $\N$; in particular, it shows that in the case  required vanishing 
holds with $\del= \frac{1}{2}\dim \N$.}
 as 
 $$Ext^i_{Coh^\LG(\gth)}
(\O\otimes V_\nu, \Phib_{IW}^{-1}(\F)).$$
Now the statement follows
from Lemma \ref{XiJT}(d). 
\epf

\section{The equivalences}\label{theeq}

\subsection{Equivalence \eqref{I0I}}

We use the criterion of Propostion \ref{fullfaithrep}(b) to show that for $\F\in D_{I^0I}$ the functor $M \mapsto Hom(\Phips(\F), M)$
is represented by an object of $D^b(Coh^\LG(St'))$; this object is then  defined uniquely up to a 
unique isomorphism in view of Proposition \ref{fullfaithrep}(a) and we obtain a functor $\Psis: D_{I^0I}\to D^b(Coh^\LG(St'))$ sending $M\in D_{I^0I}$ to the corresponding representing object.

We need to check that conditions of Proposition \ref{fullfaithrep}(b) are satisfied. 
Condition \ref{fullfaithrep}(b)(i) (representability of the restriction to $D_{perf}^{\geq n}$
for all $n$) follows from Propositions  \ref{fingenprop} and \ref{boundampl} (finite generation and bounded amplitude) in view of 
Proposition \ref{fingencrit}.
Condition \ref{fullfaithrep}(b)(ii)
(vanishing on $D_{perf}^{\leq m}$ for $m\ll 0$) follows from Proposition \ref{boundampl}.

Functor $\Psis$ is now defined. 

Proposition 
 \ref{ABsidelinatrube}(b) 
shows it is an equivalence, conditions (A) and (C) are provided respectively by 
Proposition \ref{boundampl} and Proposition \ref{gdekogo}.
  \epf

\bigskip

For future reference we record another favorable property of $\Psis$ in relation  to
the standard $t$-structures on the triangulated categories involved.

\begin{Cor}\label{exactn}
a) For $\F\in Perv_N(G/B)\subset \P$ we  have $\Psis(\F)\in Coh^\LG(St')$.

b) $\Psis(j_{w*})\in Coh^\LG(St')$ for $w\in W^f$ and $\Psis(j_{w!})\in Coh^\LG(St')$ when $w\in W_{f}\nu$, $\nu\in -\Lambda^+$.
\end{Cor}

\proof a) follows from the Proposition \ref{exaprop}. 

b)  follows from a) since $w\in W^f$ can be written as $w=w'\la$, $\la\in \La^+$,
$w'\in W_f$, so that $\ell(w)=\ell(w')+\ell(\la)$.  Then we get $\Psis(j_{w*})=\Psis(j_{w'*}*j_{\la*})=
\Psis(j_{w'*})\otimes \O(0,\la)$. Similarly, if $w=w'\nu$, $\nu\in -\Lambda^+$, then 
$\ell(w)=\ell(w')+\ell(\nu)$. \epf

\subsection{Equivalence \eqref{I0I0}}

We again use the criterion of Proposition \ref{fullfaithrep}(b) to show that the functor $\F \mapsto Hom(\Phipb(\hatt{\F}), M)$
is represented by an object of $D^b(Coh^\LG_\N(St))$, here $\hatt{\F}$ denotes pull back of $\F$ under the morphism $\Sth\to St$.
The representing  object is then  defined uniquely up to a 
unique isomorphism in view of Proposition \ref{fullfaithrep}(a) and we obtain a functor $\Psi: D_{I^0I^0}\to D^b(Coh^\LG_\N(St))$ sending $M\in D_{I^0I^0}$ to the corresponding representing object.

We need to check that conditions of Proposition \ref{fullfaithrep}(b) are satisfied. 
 In view of  Proposition \ref{fingencrit}
condition \ref{fullfaithrep}(b)(i) (representability of the restriction to $D_{perf}^{\geq n}$
for all $n$) follows from Propositions  \ref{fingenprop} and \ref{boundampl} which provide
respectively 
finite generation and bounded amplitude properties (Proposition \ref{fingenprop}
states a similar property for an object of $D$, the case of $D_{I^0I^0}$ follows).

Condition \ref{fullfaithrep}(b)(ii)
(vanishing on $D_{perf}^{\leq m}$ for $m\ll 0$) follows from Proposition \ref{boundampl}.
Since the  log monodromy endomorphisms act  on objects of $D_{I^0I^0}$  nilpotently,
this object is set theoretically supported on the preimage of $\N$ in $St$.
Thus we get the functor $\Psi: D_{I^0I^0}\to D^b(Coh^\LG_\N(St))$.

It also induces a functor between the subcategories in the categories of pro-objects:
$\Psib:\Dhat\to  D^b(Coh^\LG(\Sth))$. Recall that
 $i_{St}$ denotes for the embedding $St'\to \St$.

\begin{Lem}\label{commPsi}
a) The following diagrams commute:

$$ \begin{CD}
\Dhat  @>{\Psib}>> D^b(Coh^\LG(\Sth)) \\     
@V{\pi_*}VV @VV{i_{St}^*}V \\
D_{I^0I} @>{\Psis}>> D^b(Coh^\LG(St'))
\end{CD}
$$

$$ \begin{CD}
D_{I^0I} @>{\Psis}>> D^b(Coh^\LG(St'))\\
@V{\pi^*}VV @VV{i_{St*}}V \\
D_{I^0I^0} @>{\Psi}>> D^b(Coh^\LG_\N(St))   
\end{CD}
$$

b) We have $\Psib\circ \Phid\cong \delta_*$, where $\delta:\gth\to \Sth$ is the diagonal embedding. 

\end{Lem}

\proof Lemma \ref{compaact} implies that both compositions in the first diagram
are compatible with the action of $(D_{perf}^\LG(\Sth), \otimes )$, i.e. if $F_1=i_{St}^*\circ \Psi$,
$F_2=\Psi'\circ \pi_*$,   then $F_i(\F(X))\cong i^*(\F) \otimes F_i(X)$
canonically for 
$\F\in D_{perf}^\LG(\Sth)$, $X\in \Dhat.$
We also have  $F_1(\Xib)\cong \O\cong F_2(\Xib)$, thus we get a functorial isomorphism
$F_1(X)\cong F_2(X)$ for $X$ in the image of $\Phipb$. 


Now given  any $X\in \Dhat$, choose a bounded above complex of equivariant locally free  sheaves
$\F^\bu$ representing the object $\Psib(X)$; then using Proposition \ref{boundampl}
and the fact that $\pi_*$, $i_{St}^*$ have bounded homological dimension we get 
 for $N'\ll N\ll 0$: 
$$
F_1(X)\cong \tau^\bete_{\geq N} \Phipb (\tau^\bete_{\geq N'}(\F^\bu))\cong F_2(X).$$
It is easy to check that the resulting isomorphism does not depend on the choice
of $\F^\bu$ and is functorial in $X$, thus commutativity of the first diagram is established.

The proof for the second diagram is similar, this proves part (a). 

The same observation that all the functors involved commute with the action of $(D_{perf}^\LG(\Sth),
\otimes_\O )$ reduce (b) to checking that 
$$\Psib(\Delta_e)\cong \delta_*(\O_{\gth}).$$

To this end 
  it suffices to construct an isomorphism 
of $\O(\Cb_{\St})$-modules:
\begin{equation}\label{DelO1}
\oplusl_{\la,\mu\in \La^+}Hom_{deeq}^\bu(\Xib, \J_\la*\Del_e *\J_\mu) )\cong 
\oplusl_{\la,\mu\in \La^+}
R\Gamma(\gth, \O(\la+\mu) )
\end{equation}
compatible with the $\O(\Cb_{\St})$ action; here
subscript $_{deeq}$ refers to the $\LG$-deequivariantization
(see \S \ref{deeqsubs}).

Using  \ref{onXi}(a,b) we can rewrite the left hand side of \eqref{DelO1}
as  
$$\oplusl_{\la,\mu\in \La^*}Hom_{deeq}^\bu (Av^{IW}(\Del_e), Av^{IW}( \J_{\la+\mu}) ).
$$
Since $Av^{IW}(\Del_e)=\Phib_{IW}(\O_{\gth})$, $Av^{IW}( \J_{\la+\mu})= \Phib_{IW}(\O_{\gth}(\la+\mu))$, we see that the displayed expression is canonically isomorphic to
$\oplusl_{\la,\mu\in \La^+}Hom_{Coh(\gth)}^\bu (\O_{\gth}, \O_{\gth}(\la+\mu)),$
which yields \eqref{DelO1}; compatibility with the $\O_{\Cb_{\Sth}}$ action is clear from the construction.  
\epf

We are now ready to prove that $\Psi$, and hence $\Psib$ is an equivalence.
Since we know that $\Psis$ is an equivalence and 
the essential image of $i_*:D^b(Coh^\LG(St'))\to D^b(Coh^\LG_\N(St))$
generates the target category, Lemma \ref{commPsi}(a)
shows that the essential image of $\Psi$ generates the target category.
Thus it suffices to check that $\Psi$ is fully faithful. It is enough to see that
$$\begin{CD}
Hom(A,B)@>{\Psi}>> Hom(\Psi(A),\Psi(B))
\end{CD}$$
is an isomorphism when $B$ is obtained from an object $B'\in D_{I^0I}$
by forgetting the equivariance. This follows from the corresponding statement
for $\Psis$ and Lemma \ref{commPsi}(a). \epf

\subsection{Equivalence \eqref{II}}
\label{93}

\subsubsection{Passing from monodromic to equivariant category by killing monodromy}
\label{passing}

Let $X$ be a scheme with an action of an algebraic torus  $A$.
Let $\PP_{mon}$ be the category of unipotently
 monodromic perverse sheaves on $X$.

We have an action of $\fa=Lie(A)$ on $\PP_{mon}$ by log monodromy.
Let $\BK_\fa$ be the  Koszul complex of the vector space $\fa$;
in other words, $\BK_\fa$ is the standard complex for homology
of the abelian algebra $\fa$ with coefficients in the free module $U\fa=Sym(\fa)$. Thus 
$\BK_\fa$ is a graded commutative DG-algebra with $\fa\oplus \fa[1]$
as the space of generators and differential sending $\fa[1]$ to $\fa$ by the identity map.
It is clear that $\BK_\fa$ is
quasi-isomorphic to the base field $k$
and its degree zero part is the enveloping algebra $U\fa$.

We define a DG-category $\PP_{eq}$ as the category of complexes
of objects in $\PP_{mon}$ equipped with an action of $\BK_\fa$, such that
the action of $\fa\subset \BK_\fa^0$ coincides with the log monodromy action.
Let $D(\PP_{eq})=Ho(\PP_{eq})/Ho_{acycl}(\PP_{eq})$ be the quotient of
the homotopy category by the subcategory of acyclic complexes.

We will also write $D(X/A)$ for the $A$-equivariant derived category of constructible
sheaves on $X$ (equivalently, constructible derived category of the stack $X/A$). 

\begin{Lem} \label{toruseqDG}

a) We have a natural equivalence $D(\PP_{eq})\cong D(X/A)$ (the equivalence 
will be denoted by $real_{eq}$).

b) Consider the functors $Forg:\PP_{eq}\to Com(\PP_{mon})$ and $Ind_{U(\fa)}^{\BK_\fa}:Com(\PP_{mon})\to
\PP_{eq}$, where the first one is the functor of forgetting the $\BK_\fa$
action and the second one is the functor of induction from $U(\fa)$ which acts by log monodromy to $\BK_\fa$. 

The induced functors on the derived categories  fit into the following diagrams which commute up to a natural isomorphism:
$$\begin{CD}
D(\PP_{eq}) @>{Forg}>>  D^b(\PP_{mon})\\
@V{real_{eq}}VV @VV{real}V\\
D(X/A) @>{pr^*}>> D(X)
\end{CD}
$$

$$\begin{CD}
D(\PP_{eq})@<{Ind_{U(\fa)}^{\BK_\fa}[-d]}<< D^b(\PP_{mon})\\
@V{real_{eq}}VV @VV{real}V\\
D(X/A) @<{pr_*}<< D(X)
\end{CD}
$$

where $pr$ denotes the projection $X\to X/A$, $real$ denotes  Beilinson's realization functor \cite{Breal} and $d=\dim(\fa)$.

c) Suppose that $\F,\ \G \in \PP_{eq}$ are such that $Ext^{>0}_{D(X)}(\F^i,\G^j)=0$
for all $i,j$. 
Then $$Hom_{Ho(\PP_{eq})}(\F,\G)\iso Hom_{D(X/A)}(real_{eq}(\F),
real_{eq}(\G)).$$

\end{Lem}

\proof a) Assume first that the action of $A$ on $X$ is free and the quotient
$Y=X/A$ is represented by a scheme.
The abelian category $Perv(Y)$ of perverse sheaves on $Y$ admits a full
embedding into the category $Perv_{mon}(X)$ of unipotently monodromic perverse
sheaves on $X$, and the essential image of the embedding consists of sheaves 
with zero action of log monodromy.
Thus we have a natural embedding $Com(Perv(Y))\to Com(\PP_{eq})$ sending
a complex of equivariant sheaves to the same complex equipped with 
zero action of $\fa$ and $\fa[1]$.
We claim that the induced functor $D^b(Perv(Y))\to D(\PP_{eq})$ is an equivalence.

This claim is readily seen to be local on $Y$, i.e. it suffices to check it assuming
that $X=A\times Y$ where $A$ acts on the first factor by translations. In the latter case
the category $Perv_{mon}(X)$ is readily identified with the tensor product of the abelian 
category $Perv(Y)$ and the abelian category of unipotently monodromic local
systems of $A$, the latter is equivalent to the category of modules over the symmetric
algebra $U(\fa)\cong Sym(\fa)$ set-theoretically supported at zero
(see \cite[\S 5]{Del} for the notion of tensor product of abelian categories).
Thus the claim is clear in this case.

Let now $X$ be general. Then an object of $D(X/A)$ is by definition (see \cite{BL})
a collection of objects in $D(\tilde Y)$ given for every $A$-equivariant smooth map 
$\tilde X\to X$ where the action of $A$ on $\tilde X$ is free and $\tilde Y=\tilde X/A$,
subject to certain compatibilities.
We have the pull back functor $\PP_{eq}(X)\to \PP_{eq}(\tilde X)$, composing
it with the functor $\PP_{eq}(\tilde X)\to D^b(Perv({\tilde Y}))\cong D(\tilde Y)$
we get the desired system of objects, the compatibilities are easy to see. 

b) Commutativity of the first diagram is clear from the proof of (a) and commutativity of
the second one follows by passing to adjoint functors (notice that in view
of self-duality of Koszul complex the functor $Ind_{U(\fa)}^{\BK_\fa}[-d]$ is {\em right}
adjoint to the forgetful functor $Forg$). 

c) By a standard argument the condition in (c) implies that 
$$Hom_{Ho(\PP_{mon})}(\F,\G)\cong Hom(Forg(\F),Forg(\G)).$$
We have adjoint pairs of functors compatible with the natural functor
from the homotopy category to the derived category:
$$Ho(\PP_{eq})\overset{Forg}{\To} Ho(\PP_{mon})\overset{Ind}{\To} Ho(\PP_{eq}),$$
$$D(\PP_{eq})\overset{Forg}{\To} D^b(\PP_{mon})\overset{Ind}{\To} D(\PP_{eq}).$$
The composition in each case admits a filtration with associated graded
$Id\otimes \Lambda(\fa[1])$, i.e. for $\F \in Ho(\PP_{eq})$ or $\F\in 
D(\PP_{eq})$ we have
$$Ind\circ Forg(\F)\in  \{ \Lambda^d(\fa)\otimes \F [d] \}* \{ \Lambda^{d-1}
(\fa)\otimes \F[d-1]\} *\cdots * \{ \fa\otimes \F[1]\} *{\F},$$
where we used the notation of \cite{BBD}: $X*Y$ is the set of objects $z$
such that there exists a distinguished triangle $x\to z\to y$, $x\in X$, $y\in Y$. Since $Hom_{Ho(\PP_{eq})}(Ind\circ Forg(\F),\G)\iso
Hom_{D(X/A)}(Ind\circ Forg(\F),\G)$, it follows by induction in $n$ that
$Hom^n_{Ho(\PP_{eq})}(\F,\G) \iso Hom^n_{D(X/A)}(\F,\G)$.
 \epf
 

\begin{Cor}\label{ThII}
Let $\Th_{II}$ denote the DG category whose objects
are finite complexes of objects in $\Th$ equipped with an action of 
$\BK_{\t^2}$ such that the action of $\BK_{\t^2}^0=U(\t^2)$ coincides with the action 
induced by the torus monodromy.
Then the homotopy category $Ho(\Th_{II})$ is naturally equivalent to $D_{II}$.  
\end{Cor}

\proof Lemma \ref{toruseqDG}(c) yields a fully faithful functor $Ho(\Th_{II})\to D_{II}$.
To see that this functor is essentially surjective, notice that Lemma
\ref{toruseqDG}(b) implies that the composition of the natural functors
$Ho(\Th_{II})\to Ho(\Th)\to Ho(\Th_{II})$ contains identity functor as a direct summand
(more precisely, this composition is isomorphic to tensoring with $H^*(T^2)\in D^b(Vect)$).
Thus every object of $D_{II}$ is a direct summand in an object which belongs to the essential
image of the full embedding $Ho(\Th_{II})$. Thus we will be done if we check
that $Ho(\Th_{II})$ is Karoubian (idempotent complete).

Since a direct summand of a free-monodromic tilting object is again free-monodromic
tilting, the category $\Th_{II}$ is idempotent complete.
For $T^\bu\in \Th_{II}$ the space of closed endomorphisms of the complex
commuting with the $\BK_{\t^2}$ action is  a pro-finite dimensional ring
whose quotient by its pro-nilpotent radical is finite dimensional.
The subspace of endomorphisms homotopic to zero is a two-sided  ideal in this ring.
By elementary algebra an idempotent in a quotient of a finite dimensional
algebra by a two-sided ideal can be lifted to an idempotent in the original ring; thus we see
 that every idempotent endomorphism of an object
in $Ho(\Th_{II})$ lifts to an idempotent in the ring of endomorphisms of the corresponding
object in $\Th_{II}$, this shows that $Ho(\Th_{II})$ is idempotent complete. \epf

We are now ready to establish \eqref{II}.

Consider the category of finite complexes of objects in $Coh^\LG(\Sth)$
equipped with an action of $\BK_{\ft^2}$ extending the action of $\t^2=(\Lt^*)^2$ 
coming from the action of linear functions on $\Lt^2$ pulled back under the natural
map $St\to \Lt^2$. (It is easy to see that replacing  $Coh^\LG(\Sth)$ in the previous
sentence by $ Coh^\LG(\St)$ one gets definition of an equivalent category).
Let $Coh_{\BK_{\ft^2}}^\LG(St)$ denote this category
and $Ho(Coh_{\BK_{\ft^2}}^\LG(St))$ be the corresponding homotopy category.

It follows from the definition of the derived coherent category of a DG-scheme 
that there exists a natural functor 
$$real_{coh}:Ho(Coh_{\BK_{\ft^2}}^\LG(St))\to DGCoh^\LG(St\Ltimes_{\Lt^2}
\{0\})=DGCoh^\LG(\Nt\Ltimes _\Lg\Nt).$$

Moreover, given two complexes $\F^\bu,\G^\bu\in Coh_{\BK_{\ft^2}}^\LG(St)$
such that $Ext^{>0}_{Coh^\LG(\Sth)}(\F^i,\G^j)=0$ we have
$$Hom_{Ho(Coh_{\BK_{\ft^2}}^\LG(St))}(\F^\bu, \G^\bu)\iso Hom (real_{coh}(\F^\bu),
real_{coh}(\G^\bu)).$$

Corollary \ref{noExtZT} implies that the functor $\Psib$ sends free-monodromic tilting sheaves
to coherent sheaves;
thus equivalence $\Phi_{I^0I^0}$ and Corollary \ref{ThII} yield a fully faithful functor
$\Psi_{II}:D_{II}\to DGCoh^\LG(\Nt\Ltimes_\Lg \Nt)$. 

It remains to show that $\Psi_{II}$ is essentially surjective.


For $w\in W$ consider the full triangulated subcategory $D_{II}^w$ in $D_{II}$ generated by $j_{w*}$. We first
check that the functor $\Psi_{II}|_{D_{II}^w}$ admits a left adjoint. This is equivalent to existence for any $\F \in DGCoh^\LG(\Nt\Ltimes_\Lg \Nt)$
of an object $\F_w\in D_{II}^w$ and a morphism $s_w:\F\to \Psi_{II}(\F_w)$ such that 
\begin{equation}\label{cones_w} Hom^\bu(Cone(s_w), \Psi_{II}(j_{w*}))=0.
\end{equation}
Notice that for $\F,\G \in DGCoh^\LG(\Nt\Ltimes_\Lg \Nt)$ we have
\begin{equation}\label{Homnz} Hom^\bu (\F,\G)=0 \iff Hom^\bu_{D^bCoh^\LG(St')}(\iota \F,\iota \G)=0, 
\end{equation}
where $\iota$ denotes the functor $\iota:DGCoh^\LG(\Nt\Ltimes_\Lg \Nt)\to D^bCoh^\LG(St')$ of
direct image under the closed embedding of derived schemes.

In fact, a stronger statement holds: $Hom^n(\F,\G)=0$ for $n\ll 0$ and
 if $n$ is the minimal integer such that $Hom^n (\F,\G)\ne 0$
then $n$ is also the minimal integer for which $Hom^n_{D^bCoh^\LG(St')} (\iota \F,\iota \G)\ne 0$;
moreover, the natural map $Hom^n (\F,\G)\to Hom^n_{D^bCoh^\LG(St')}(\iota \F,\iota \G)$ is nonzero.
These statements follow from existence of a spectral sequence with the $E_2$
term being $E_2^{p,q}=Hom^p(\F,\G)\otimes \Lambda^q(\Lt)$ converging
to $Hom^{p+q}_{D^bCoh^\LG(St')}(\iota \F,\iota \G)$.

The functors $\iota$, $\Psi_{II}$ and $\Psi:=\Phi_{I^0I}^{-1}$ fit in the diagram:

\begin{equation}\label{CommD}
\begin{CD}
D_{II} @>{\Psi_{II}}>> DGCoh^\LG(\Nt\Ltimes_\Lg \Nt) \\
@V{Forg}VV @VV{\iota}V \\
D_{I^0I}@>{\Psi}>> D^b(Coh^\LG(\St')) 
\end{CD} 
\end{equation}
whose commutativity follows from Lemma \ref{toruseqDG}(b).

We now show existence of $s_w$ satisfying \eqref{cones_w} by induction in: 

\noindent $d=\dim Hom^\bu_{D^b(Coh^\LG(St'))}(\iota \F,\Psi(j_{w*}))$. Notice that
this space is finite dimensional since $\Psi$ is an equivalence and the corresponding statement in $D$ is clear.
If $d=0$ then $\F_w=0$ satisfies the required property 
in view of \eqref{Homnz} and commutative diagram \eqref{CommD}. If $d\ne 0$ let us choose $f\in Hom^\bu(\F,\Psi_{II}(j_{w*}))$ whose image in
$Hom^\bu_{D^b(Coh^\LG(St'))}(\iota \F,\Psi (j_{w*}))$ is nonzero: this is possible by the observation following \eqref{Homnz} combined with \eqref{CommD}. Set $\F'= Cone(f)[-1]$. 
Then $$\dim Hom^\bu_{D^bCoh^\LG(St')}(\F ',\Psi(j_{w*}))=d-1,$$ so by the induction assumption we can find
a distinguished triangle $C\to \F' \to \F'_w\to C[1]$ where $\F'_w\in \Psi_{II}(D_{II}^w)$ and $C\in ^\perp \Psi_{II}(D_{II}^w)$. 
If $f'$ is the composed arrow $C \to \F' \to \F$ then using the octahedron axiom one checks that
$Cone(f')\in \Psi_{II}(D_{II}^w)$, so the map $s_w:\F \to Cone(f')$
satisfies  \eqref{cones_w}.

For a finite subset $S\subset W$ let $D_{II}^S$ be the full triangulated subcategory generated by 
$j_{w*}$, $w\in S$. 
Using that $Hom^\bu(j_{w*}, j_{v*})=0$ for $v\not \preceq w$ one deduces by induction in $|S|$ that $\Psi_{II}|_{D_{II}^S}$ has a left adjoint or, equivalently,
that for $\F \in DGCoh^\LG(\Nt\Ltimes_\Lg \Nt)$ there exists 
an object $\F_S\in D_{II}^S$ and a morphism $\F \to \Psi_{II}(\F_S)$ whose cone lies in $^\perp \Psi_{II}(D_{II}^S)$.

Now, for $\F\in DGCoh^\LG(\Nt\Ltimes_\Lg \Nt)$ there exist only finitely many $w$
for which $Hom^\bu(\F, \Psi_{II}(j_{w*}))\ne 0$: in view of the commutative diagram 
\eqref{CommD} this follows from \eqref{Homnz} and the corresponding property of $D_{I^0I}$. Given such an $\F$  choose a finite subset  $S\subset W$ containing all such $w$ and closed under descent in the Bruhat order. 
Then  $C=Cone(\F \to  \Psi_{II}(\F_S))$ satisfies $Hom^\bu(C,\Psi_{II}(j_{w*}))=0$ for all $w\in W$: if $w\in S$ this follows from the definition of $\F_S$, otherwise
we have $Hom^\bu(\F,  \Psi_{II}(j_{w*}))=0=Hom^\bu(\F_S, j_{w*})$. By \eqref{Homnz} we see that $Hom^\bu(\iota C, \Psi(j_{w*}))=0$ for all $w$
which implies $\iota C=0$, hence $C=0$, thus $\F$ lies in the essential image of $\Psi_{II}$.  
 \qed

\section{Monoidal structure}\label{monstr}

\subsection{A $DG$-model for convolution of coherent sheaves}

 \begin{Lem}\label{funcomplcorr}
  Let $X$, $Y$ be two algebraic stacks and $F=F_\CK:D^b(Coh(X))\to D^b(Coh(Y))$
 be a functor coming from an object $\CK\in D^b(Coh(X\times Y))$,
 i.e. $F:\F\mapsto pr_{2*}(\CK\Lotimes pr_1^*(\F))$. Let $M\in D^b(Coh(X))$ be represented by a complex  of sheaves $M^\bu$ such that $F(M^i)\in Coh(Y)$. Then $F(M)$ is canonically isomorphic to the object represented by $F(M^\bu)$.

 \end{Lem} 
 
 \proof  A functor as above  lifts
 to a functor between filtered derived categories
 $F^{fil}:DF(Coh(X))\to DF(Coh(Y))$.
 Recall 
 that $DF$ contains the  category
 of bounded complexes  in $Coh(X)$ as a full subcategory, the canonical
 functor from the filtered derived category to the derived category restricted
 to this subcategory coincides with the canonical functor from the category
 of complexes to the derived category. The conditions
 of the Lemma show that $F^{fil}$ sends the object corresponding to the complex $M^\bu$
 to the object corresponding to $F(M^\bu)$, which yields the desired statement.
  \epf
 
 Recall from  \S \ref{dual_side} and \cite{BeRi} that for a proper map  $X\to Y$ of smooth varieties
 convolution yields a monoidal structure on the derived coherent category of the 
 DG-scheme $X\Ltimes _Y X$. If $X\to Y$ is semi-small then $Tor_i^{\O_Y}(\O_X, \O_X)=0$
 for $i>0$ thus we get a convolution monoidal structure on $D^b(Coh(X\times_Y X))$
 and on $D^b(Coh^H(X\times_Y X))$ for an algebraic group $H$ acting compatibly on 
 $X$, $Y$. These monoidal categories act on module
categories $D^b(Coh(X))$, $D^b(Coh^H(X))$ respectively, the action functor
is also denoted by $*$.

 \begin{Cor}\label{convcompl}
  Let $X\to Y$ be a proper semi-small
   morphism of smooth quasiprojective varieties equipped with
 an action of a reductive algebraic group $H$. 

 a)  Let $\F^\bu, \G^\bu$
   be finite complexes  of  $H$-equivariant coherent  sheaves on $X\times_Y X$
  such that the convolution  $\F^i*\G^j$ lies in  $Coh^H(X\times_Y X)$
 for all $i,j$. 
  Let $\F$, $\G$ be the corresponding
  objects in the derived category. Then $\F*\G$ is canonically isomorphic
  to the object represented by the total complex of the bicomplex $\F^i*\G^j$.

b) Assume that three complexes $\F_1^\bu$, $\F_2^\bu$, $\G^\bu$
of $H$-equivariant coherent sheaves on $X\times _Y X$
are such that $\F_1^i*\F_2^j$, $\F_2^j*\G^l$  and $\F_1^i*\F_2^j*\G^l$ 
lie in $Coh^H(X\times _Y X)$ for all $i,j,l$. 
Then the two isomorphisms provided by part (a)
 between $\F_1*\F_2*\G\in D^b(Coh^H(X\times_Y X))$
and the object represented by the complex 
  $\CC^d=\oplusl_{i+j+l=d}\F_1^i*\F_2^j*\G^l$
 coincide.   

c) 
Let    $\F^\bu$ be as in (a) and $\G^\bu$
be a   finite complex  of  $H$-equivariant coherent  sheaves on $X$.
Then  $\F*\G$ is canonically isomorphic
  to the object represented by the total complex of the bicomplex $\F^i*\G^j$ provided
that $\F^i*\G^j \in Coh^H(X)$. 

d) Let $\F_1^\bu,\,  \F_2^\bu$ be as in (b) and $\G$ as in (c).
Assume that $\F_1^i*\F_2^j\in Coh^H(X\times _Y X)$, while
 $\F_2^j*\G^l, \, \F_1^i*\F_2^j*\G^l\in Coh^H(X)$ for all $i,j,l$. 
Then the isomorphism between $\F_1*\F_2*\G\in D^b(Coh^H(X))$
and the object represented by the complex
 $\CC^d=\oplusl_{i+j+l=d}\F_1^i*\F_2^j*\G^l$ 
obtained by applying part (c) twice coincides
with the isomorphism obtained by applying part (a) and part (c).   
  \end{Cor}
 
 \proof The  convolution product comes from a functor 
$$F:D^b(Coh^{H}((X\times_Y X)^2) 
 \to D^b(Coh^H(X\times _Y X))$$
 of  the type considered in Lemma \ref{funcomplcorr}, namely we have $F=F_\CK$,
 where $\CK\in D^b(Coh^{H}(X\times _YX)^3)$ is given by $\CK=\upsilon^*
 \delta_*(\O_{X^3})$; here $\upsilon$ stands for the embedding $(X\times _Y X)^3\to
 (X\times X)^3=X^6$ and $\delta:X^3\to X^6$ is given by 
$(x_1,x_2,x_3)\mapsto (x_3, x_1,x_1,x_2,x_2,x_3)$. Thus statement (a) follows from Lemma \ref{funcomplcorr}.  
 Part b) follows by considering the functor between filtered derived categories
 $DF(Coh^H(X\times _Y X)^3)\to DF (Coh^H(X\times _Y X))$
 corresponding to the triple convolution. 
The proof of  (c,d) is similar to the proof of (a,b) respectively. 
 \epf
 
 \begin{Lem}\label{actfulfaith} 
a)  Let $X\to Y$ be a semi-small proper morphism of smooth quasi-projective varieties equipped with
 an action of a reductive algebraic group $H$.
For $\F\in D^b(Coh^H(X\times_Y X))$ let $a(\F)$ denote the corresponding
functor $D^b(Coh^H(X))\to D^b(Coh^H(X))$. 

For $\F\in Coh^H(X\times_Y X)$, $\F'\in D^b(Coh^H(X\times_Y X))$ any isomorphism
of functors $a(\F)\cong a(\F')$ comes from a unique isomorphism $\F\cong \F'$. 

b) Given an $H$-invariant closed subvariety $Z\subset Y$, the statement in (a) remains
true for $\F,\, \F'\in Coh^H(\widehat{X\times_Y X})$, $a(\F),\, a(\F')\in End(D^b(Coh^H(\widehat{X})))$,
where $\widehat{X}$, $\widehat{X\times_Y X}$ denote formal completions at the preimage of $Z$.
\end{Lem}

\proof In the setting of either part (a) or part (b), an equivariant
 coherent sheaf $\F$ can be reconstructed from the corresponding module $M(\F)$ over the
  homogeneous coordinate ring,

 $M(\F)=\oplusl_{n,m\geq 0} \Gamma(\F\otimes pr_1^*(L^n)
  \otimes pr_2^*(L^m))$,
  where $L$ is an equivariant ample line bundle on $X$. Thus Lemma follows from
the following expression for $M(\F)$ in terms of the functor of convolution by $\F$:
$M(\F)=\oplusl_{m,n}Hom_{deeq}(L^{-n}, \F*L^m)$. \epf





\subsection{Monoidal structure on $\Phi_{I^0I^0}$}\label{monI0I0}

\begin{Lem}\label{compIWact}
The equivalence $\Phib$ is compatible with the action on 
$D^b(Coh^\LG(\gth))\cong \Dhat_{IW}$ via
the equivalence $\Phib_{IW}$, i.e. we have a functorial isomorphism
$$\Phib_{IW}(\F*\G)\cong \Phib(F)*\Phib_{IW}(\G)$$
where  $\F\in D^b(Coh^\LG(\Sth))$, $\G\in D^b(Coh^\LG(\gth))$.
 \end{Lem}
 
 \proof For $\F\in D_{perf}^\LG(\Sth)$
 this is Proposition \ref{compa_act}.

Let now $\F$ be general. For any sufficiently large $N$ we can find $\F' 
\in D_{perf}^\LG(\Sth)$ such that $\F=\tau_{\geq -N}(\F')$. 
The functor $D^b(Coh^\LG(\Sth))\to D^b(Coh^\LG(\gth))$, $\F\mapsto \F*\G$
 has bounded homological
amplitude; the functor $\Dhat \to \Dhat_{IW}$
$X\mapsto X*\Phib_{IW}(\G)$ has bounded homological amplitude and 
by 
Proposition \ref{gdekogo}
the functor
$\Phib$ has homological amplitude bounded above, i.e. it sends
$D^{\leq 0}(Coh^\LG(\Sth))$ to $D^{\leq n} (\Phat)$ for some $n$. 
It follows that for $N\gg m\gg 0$ and $\F'$ as above we have
\begin{multline*}\Phib_{IW}(\F*\G)\cong \Phib _{IW}(\tau_{\geq -m}(\F'*\G))\cong
\tau_{\geq -m}\Phib_{IW}(\F'*\G)\cong \tau_{\geq -m}(\Phib(\F)*\Phib_{IW}(\G))
\cong \\
 \Phib(\F)*\Phib_{IW}(\G),
\end{multline*}
which proves the Lemma.
 \epf

\medskip

We are now ready to equip $\Phib$ with a monoidal structure.
We work with the inverse equivalence $\Psib$.
We need to construct an isomorphism
\begin{equation}\label{moneq}
\Psib (\F*\G)\cong \Psib(\F)*\Psib(\G)
\end{equation} 
compatible with the associativity isomorphism.

Given $\F,\, \G\in \Dhat$ and $\M\in \Dhat_{IW}$, Lemma \ref{compIWact}
 provides isomorphisms
$$\Psib(\F*\G)*\Psib_{IW}(\M)\cong 
\Psib_{IW}(\F*\G*\M)\cong \Psib(\F)*\Psib_{IW}(\G*\M)\cong \Psib(\F)*\Psib(\G)*\Psib_{IW}(\M),$$
where $\Psib_{IW}$ is the equivalence inverse to $\Phib_{IW}$.

Thus we get an isomorphism 
$$a(\Psib(\F*\G))\cong  a(\Psib(\F))\circ a(\Psib(\G)) \cong a(\Psib (\F)*\Psib(\G)),$$
where we used notations of Lemma \ref{actfulfaith}. This
isomorphism  is compatible with the associativity
constraint, since the equivalence $\Phib_{IW}$ sends the corresponding equality to an equality
which holds since monoidal category $\Dhat$ acts on $\Dhat_{IW}$. 

Since $\Psib:\Th\to Coh^\LG(\Sth)$, Lemma \ref{actfulfaith}(b) yields
\eqref{moneq} in the case when $\F,\G\in \Th$, which is 
compatible with the associativity
isomorphism for three objects in $\Th$. Now Corollary \ref{convcompl}(a) compared
to Proposition \ref{tiltprop}(b)
yields \eqref{moneq} in general, while Corollary \ref{convcompl}(b) together with Proposition
\ref{tiltprop}(b) shows that
 the constructed isomorphism  is compatible with associativity constraint. \epf

\subsection{Monoidal structure on $\Phi_{II}$}\label{monII} 
\subsubsection{A monoidal structure on $Ho(\Th_{II})$}
In order to equip $\Phi_{II}$ with a monoidal structure we describe the monoidal structure
on $D_{II}$ in terms of the DG-model $\Th_{II}$ (see Corollary \ref{ThII}).

Let $\Th_{II}^{(2)}$ denote the category of finite complexes of objects in $\Th$
equipped with an action of $\BK_{\t^2}\otimes \Lambda(\t[1])$, and $\Th_{II}^{(3)}$
be the category of  finite complexes of objects in $\Th$
equipped with an action of $\BK_{\t^2}\otimes \Lambda(\t^2[1])$.
In both cases we require that $\t^2\subset \BK_{\t^2}$ acts by logarithm of monodromy.

We have a functor $\Th_{II}\times \Th_{II}\overset{\star}{\To} \Th_{II}^{(2)}$
sending  $(T_1, T_2)$ to the convolution $T_1*T_2$; the latter complex is equipped
with two actions of $\BK_\t$ coming respectively from the left action on $T_1$
and the right action on $T_2$. To define the action 
of $\Lambda(\t[1])$ observe that  the right monodromy action on $T_1$ and the left
monodromy action on $T_2$ induce the same action on $T_1*T_2$, the diagonal
action of $\BK_{\t}$ kills the augmentation ideal of $\BK_\t^0=Sym(\t)$, thus it factors
through an action of $\Lambda(\t[1])$.

Similarly, we have a functor $\Th_{II}\times \Th_{II}\times \Th_{II}\overset{\star_2}{\To} \Th_{II}^{(3)}$
sending $(T_1,T_2,T_3)$ to $T_1*T_2*T_3$ where the two actions of $\BK_\t$ come respectively from the left action on $T_1$
and the right action on $T_3$, and the two actions of $\Lambda(\t[1])$ come
from the diagonal action of $\BK_{\t}$ on the first and the second factor, and 
 the diagonal action of $\BK_{\t}$ on the second and the third factor respectively.
 We use the same notation $\star$, $\star_2$ for the corresponding functors on the homotopy
 categories. 
 
Furthermore, we have functors $\mu: Ho( \Th_{II}^{(2)})\to Ho(\Th_{II})$, $\mu:M \mapsto \M\Lotimes_{\La(\t[1])}
k$ and $\mu^{(2)}:Ho(\Th_{II}^{(3)})\to Ho(\Th_{II})$, $\mu^{(2)}:M\mapsto M\Lotimes_{\La(\t^2[1])}
k$.

The following Proposition obviously yields a monoidal structure on the equivalence
\eqref{II}.

\begin{Prop}
a) The product $(M_1,M_2)\mapsto \mu (M_1\star M_2)$ makes
$Ho(\Th_{II})$ into a monoidal category, where the associativity constraint
comes from the natural isomorphisms
\begin{equation}\label{mu2}
(M_1\ot M_2)\ot M_3\cong \mu^{(2)}\star_2(M_1,M_2,M_3)\cong M_1\ot(M_2\ot M_3).
\end{equation}

b) The equivalence $real_{eq}:Ho(\Th_{II}) \cong D_{II}$ is naturally enhanced to a monoidal functor.

c) The equivalence $Ho(\Th_{II})\cong DGCoh^\LG(\Nt\Ltimes_\Lg \Nt)$ 
is naturally enhanced to a monoidal functor.
\end{Prop}

\proof To check (a) and (b) it suffices to provide a bi-functorial isomorphism
$$real_{eq}(M_1\star M_2)\cong real_{eq}(M_1)*real_{eq}(M_2)$$
sending the isomorphism \eqref{mu2} to the associativity constraint in $D_{II}$.
This follows from the next Lemma \ref{52}.

c) follows from the definition of convolution in $DGCoh^\LG(\Nt\Ltimes_\Lg \Nt)$.
\epf

In order to state the next Lemma we return to the setting of \ref{passing}. Let $X$ be an algebraic variety equipped with an action
of an algebraic torus $A$ and let $f:X/A\to Y$ be a map where $Y$ is an algebraic variety
and $X/A$ is the stack quotient. Let $pr:X\to X/A$ be the projection and set
 $\tilde f=f\circ pr:X\to Y$.

\begin{Lem} \label{52}
a) Let $M^\bu\in \PP_{eq}$ be a complex of monodromic perverse sheaves on $X$
equipped with a $\BK_\fa$ action and let $\bar M$ be the corresponding object
in $D(X/A)$ (see Lemma \ref{toruseqDG}).

Assume that $\tilde f_*(M^i)$ is a perverse sheaf for all $i$.

We then have a canonical isomorphism 
$$f_*(\bar{M})\cong real_{eq}(\tilde f_*(M^\bu)\Lotimes _{\Lambda
(\fa[1])} k[\dim (A)]).$$

b) Assume also that
a torus $A'$ acts on $X,\, Y$ so that $f$ is $A'$-equivariant and the action on $X$ commutes
with $A$. Let $\bar{f}$ be the  morphism $X/(A\times A')\to Y/A'$.

 Let $M^\bu\in \PP_{eq}$ be a complex of monodromic perverse sheaves on $X$
equipped with a $\BK_{ \fa\oplus \fa'}$ action and let $\bar M$ be the corresponding object
in $D(X/(A\times A'))$.

Assume that $\tilde f_*(M^i)$ is a perverse sheaf for all $i$.

We then have a canonical isomorphism of objects in $ D(Y/A')$:
$$\bar{f}_* (\bar{M})\cong real_{eq}(\tilde f_*(M^\bu)\Lotimes _{\Lambda (\fa[1])} k[\dim (A)]).$$

\end{Lem}

\proof a) is a particular case of b), which we will presently deduce  from the following two statements:

I) The equivalence of Lemma \ref{toruseqDG}(a) satisfies the following functoriality.
Consider an $A$-equivariant map of schemes $f:X\to Y$ and use Lemma 
\ref{toruseqDG}(a)
to identify $D_A(X)\cong D(\P_{eq}(X))$,  $D_A(Y)\cong D(\P_{eq}(Y))$. 
Then for $\F^\bu\in \P_{eq}(X)$ such that $f_*(\F^i)\in Perv(Y)$ the object
of $\P_{eq}(Y)$ obtained from $\F^\bu$ by term-wise application of $f_*$ corresponds to the object
$f_*(\F)\in D_A(Y)$. 

The special case of this functoriality where the group $A$ is trivial
 is checked in \cite{Breal}, the general case is similar.


II) For a subtorus $ A'$ of $A$  the functor $Res_{\BK_{\fa'}}^{\BK_{\fa}}:
\P_{eq}^A(X)\to \P_{eq}^{ A'}(X)$ corresponds 
under the equivalence of Lemma \ref{toruseqDG}(a) to the restriction of equivariance
 functor $Res^A_{A'}:D_A(X) \to D_{A'}(X)$, while the functor 
 $Ind_{\BK_{\fa'}}^{\BK_{\fa}}[\dim A'-\dim A]: \P_{eq}^{ A'}(X)\to
\P_{eq}^A(X)$ corresponds to the functor of direct image under the morphism of stacks
$X/A'\to X/A$. 
 
 This is a straightforward generalization of Lemma \ref{toruseqDG}(b).
 
 Now  (I) applied to the torus $A\times A'$ acting compatibly on $X$,
  $Y$ yields a description of an object in $\P_{eq}^{A\times A'}(Y)$ representing  
  $\bar{f}_*(\bar{M})$: 
 we have 
  $\bar{f}_*(\bar{M})\cong real_{eq}(\tilde f_*(M^\bu))$, where 
$\tilde f_*(M^\bu)$ is equipped with the action of $\BK_{\fa\oplus \fa'}$ inherited from the action on $M^\bu$. However, $\bar{f}_*(\bar{M})$ can also be rewritten as the 
 direct image of $f_*(\bar{M})$
under the morphism $Y/A'\to Y/(A\times A')$. Using (II) we 
get an isomorphism
$$ \Lambda(\fa[1])\otimes _k  (real_{eq}^{A'})^{-1} ( f_*(\bar{M}))[-\dim (A)] \cong
(real_{eq}^{A\times A'})^{-1}(\bar{f}_*(\bar{M})),$$
where we have adorned $real_{eq}$ with an additional superscript  making clear in which
equivariant category it lands. Applying the functor $-\otimes_{\Lambda(\fa[1])} k[\dim (A)]$
to both sides of the last isomorphism we get the Lemma.
 \epf

\subsection{Compatibility of \eqref{I0I} with the action of categories from \eqref{I0I0}, \eqref{II}}

To finish the proof of Theorem \ref{mainthm} it remains to establish compatibility
of equivalence \eqref{I0I} with the structure of a module category over the monoidal 
categories appearing in \eqref{I0I0} and \eqref{II}.  

To check compatibility with the action of $D_{I^0I^0}\cong D^b(Coh^\LG_\N(St))$
we pass to the pro-completions and check compatibility of \eqref{I0I} with the action
of $\Dhat \cong D^b(Coh^\LG(\Sth))$.  We have an action of the monoidal category of 
free monodromic tilting complexes $\Th$ on the category of tilting objects $\TT\subset \P$
which induces a structure of a module category for $Ho(\Th)$ on $Ho(\TT)$.
In view of Proposition \ref{tiltprop}(c)
this module structure is compatible with one arising from the equivalences
$Ho(\TT)\cong D$, $Ho(\Th)\cong \Dhat$. On the other hand,
using Corollary \ref{convcompl}(c,d) we see that 
  the equivalence
$Ho(\TT)\cong D^b(Coh^\LG(St'))$
is compatible with the action of $Ho(\Th)\cong D^b(Coh^\LG(\Sth))$.
This yields compatibility with the action of categories in \eqref{I0I0}.
  
To check compatibility with the action of $D_{II}\cong DGCoh^\LG(\Nt\Ltimes_\N \Nt)$
we use Lemma \ref{toruseqDG} to identify $D_{I^0I}$ with the homotopy
category of complexes in $\Th$ equipped with an action of $\BK_{\t}$ compatible
with the right log monodromy action.  This category of complexes carries a natural
action of the monoidal DG-category of complexes in $\Th$ with a compatible action 
of $\BK_{\t^2}$. The resulting triangulated module category is module equivalent
to both $D^b(Coh^\LG(St'))$ and $D_{I^0I}$ by an argument parallel to that of section 
\ref{monII}. 

This establishes the  compatibilities thereby completing the proof of Theorem \ref{mainthm}.

\subsection{Compatibility with projections}

We finish the section by recording another useful  compatibility between equivalences of
Theorems \ref{mainthm} and \ref{mainthm2}.

\begin{Prop}\label{commupr}
The following diagrams commute 

$$\begin{CD}
D^b(Coh^\LG(\Sth)) @>{\Phib}>> \Dhat \\
@V{pr_{Spr,1*}}VV @VV{Av_{IW}}V \\
D^b(Coh^\LG(\gt)  @>{\Phib_{IW}}>>\Dhat_{IW}
\end{CD}$$

$$\begin{CD}
D^b(Coh^\LG(St')) @>{\Phib}>> D_{I^0I} \\
@V{pr'_{Spr,1*}}VV @VV{Av_{IW}}V \\
D^b(Coh^\LG(\Nt)  @>{\Phi_{IW}}>>D_{IW}^I
\end{CD}$$
\end{Prop}

\proof  These diagrams are obtained from the last two diagrams in Lemma \ref{YYY}
by passing to adjoint functors. Alternatively, the first diagram above follows from compatibility 
of the equivalences \eqref{hatI0I0} and $\Phi^{I^0}_{IW}$ with the action of the monoidal category on the module category, since 
$ pr_{Spr,1*}(\F)=\F*\O_{\gt}$ and $Av_{IW}(\G)=\G*\Delta_i^{IW}$ for $\F\in 
D^b(Coh^\LG(St))$, $\G\in \Dhat$. \qed 


\section{Further properties}\label{sec_fur}
In this section we mention  further properties and possible generalizations of the constructed
equivalences.

\subsection{Frobenius compatibility}
As pointed out in the Introduction, our main result is  inspired by different geometric realizations
of the affine Hecke algebra. However, the Grothendieck group of the  categories in Theorem
\ref{mainthm} is isomorphic to a less interesting ring $\Zet[W]$.
 A possible "upgrade" of the Theorem involving
categories whose  Grothendieck group  is related to the affine Hecke algebra is an equivalence
between $D^b(Coh^{\LG\times \Gm}_\N(St))$ and an appropriately defined  mixed version of
$D_{I^0I^0}$. However, for many applications (cf. \cite{BM}, \cite{BHum}) the following simpler
version is sufficient.

Fix a finite field $\Fq$ and assume that the base field $k=\overline{\Fq}$. Then the categories
in the left hand side of \eqref{I0I0}--\eqref{PhiIWmon}
carry an automorphism coming from the Frobenius automorphism
of $k$. 

Let $\bq: St\to St$ be the map given by $\bq: (x,\b_1,\b_2)\mapsto (qx,\b_1,\b_2)$.
We use the same letter to denote the induced automorphisms of $St'$, $\Sth$ etc.

\begin{Prop}
The equivalences in Theorems \ref{mainthm}, \ref{mainthm2} intertwine Frobenius automorphism with the functor
$\bq^*$ acting on the derived categories of coherent sheaves.
\end{Prop}

The proof is parallel to the proof of \cite[Proposition 1]{AB}.

\subsection{Category $\P$ and the noncommutative Springer resolution}\label{CatPS}
Recall that the main result of \cite{BM}\footnote{Note the difference of notation: the
group denoted here by $\LG$ is denoted by $G$ in \cite{BM}.}
 is a construction of a certain
noncommutative $\O(\Lg)$ algebra $A$ and its quotient $A^0$ with  derived equivalences
$ D^b(Coh(\gt))\cong D^b(A-mod_{fg})$,
$ D^b(Coh(\Nt))\cong D^b(A^0-mod_{fg})$, see \cite[\S 1.5.3; Theorem 1.5.1(b)]{BM}. 
In fact, the algebras are defined and $A=End(\EE)^{op}$, $A^0=End(\EE_0)^{op}$,
$\EE_0=\EE|_\Nt$
for a tilting vector bundle $\EE$ on $\gt$, where the equivalence is given by $\F \mapsto 
RHom(\EE,\F)$ (respectively, $RHom(\EE_0,\F)$).

The algebras come equipped with a natural $\LG$ action and equivalences admit an
equivariant version. Furthermore, applying a version of \cite[Theorem 1.5.1(b)]{BM}
to the group $\LG\times \LG$ one gets an equivalence
$$ D^b(Coh^{\LG} (St'))\cong D^b(A\otimes _{\O(\Lg)} A^0-mod_{fg}^{\LG}) $$
given by $\F\mapsto RHom(pr_1^*(\EE)\otimes pr_2^*(\EE_0),\F)$.

We will also need a modification of this construction: the  bundle 
$\EE_0^*$  dual to $\EE_0$ is also tilting and we have an equivalence
$$ D^b(Coh^{\LG} (St'))\cong D^b(A\otimes _{\O(\Lg)} (A^0)^{op}-mod_{fg}^{\LG}) $$
given by $\F\mapsto RHom(pr_1^*(\EE)\otimes pr_2^*(\EE_0^*),\F)$.

Composing it with equivalence \eqref{I0I} and recalling that $D_{I^0I}\cong
D^b(\P)$ we get an equivalence
\begin{equation}\label{AA}
D^b(\P) \cong D^b(A\otimes _{\O(\Lg)} (A^0)^{op}-mod_{fg}^{\LG}) .
\end{equation}

We now describe the relation between the natural $t$-structures on the two sides
of \eqref{AA}. 

For a nilpotent orbit $O\subset \N$ consider
the full subcategory of 
complexes such that each cohomology module considered as a module over the 
center $\O(\N)\subset A\otimes_{\O(\Lg)} (A^0)^{op}$ 
is set theoretically supported on the closure of $O$. These subcategories define a filtration
by thick subcategories on the triangulated category $D^b(A\otimes _{\O(\Lg)} (A^0)^{op}-mod_{fg}^{\LG})$
indexed by the partially ordered set of nilpotent orbits. We will refer to this filtration as the support
filtration. 

For an orbit $O\subset \N$ we let $j^O:O\to \N$ be the embedding, $d_O=\frac{\dim O}{2}$,
$d=\frac{\dim \N}{2}$. Recall the {\em perverse $t$-structure of middle perversity} on $D^b(Coh^\LG(\N))$,
\cite[Example 4.15]{AriB} and the minimal extension functor $j^O_{!*}$ from equivariant
perverse coherent sheaves on $O$ to those on $\N$, \cite[\S 4]{AriB}. A straightforward generalization of {\em loc. cit.} produces a functor on (a subcategory of) the derived
category of $\LG$-equivariant modules over a finite $\O(\N)$ algebra equipped with
a compatible $\LG$-action. 

\begin{Thm}\label{suppf}
a) The support filtration is compatible with the image of the tautological $t$-structure
on $D^b(\P) $ under the equivalence \eqref{AA}.  The induced $t$-structure on the
associated graded category corresponding to the nilpotent orbit $O$ coincides
with $t$-structure coming from the tautological $t$-structure on 
$D^b(A\otimes _{\O(\Lg)} (A^0)^{op}-mod_{fg}^{\LG}) $ shifted by $d-d_O$.

b) Let $\F\in D_{I^0I}\cong D^b(\P)$ be an object and 
$\M\in D^b(A\otimes _{\O(\Lg)} (A^0)^{op}-mod_{fg}^{\LG})$ be its image under \eqref{AA}. Then the following are equivalent:

\begin{enumerate} \item $\F\in \P$.

\item $Forg(\M)[-d]$ is a perverse coherent sheaf for the middle perversity. Here $Forg: A\otimes _{\O(\Lg)} (A^0)^{op}-mod_{fg}^{\LG}
\to Coh^\LG(\N)$ is the forgetful functor. 

\end{enumerate}

c) For $\F$, $\M$ as in (b) the following are equivalent:
\begin{enumerate} \item $\F$ is an irreducible object in $\P$.

\item There exists an orbit $O\subset \N$ and an irreducible object $\L$ in the category
of $\LG$-equivariant $A\otimes (A^0)^{op}|_O$-modules, such that $\M[d]=j^O_{!*}(\L[d_O])$.
\end{enumerate}
\end{Thm}

\proof Part (c) follows from (b) by a straightforward generalization of \cite[Proposition 4.11]{AriB}, 
while (b) follows from part (a) by comparing it with the definition of the perverse coherent 
$t$-structure in \cite{AriB}.

 The proof of part (a) is mostly parallel to the proof of 
\cite[Theorem 6.2.1]{BM} which asserts the similar property of equivalence
\eqref{PhiIW}. More precisely, the arguments of {\em loc. cit.}  show that the image
of the tautological $t$-structure is uniquely characterized by two  properties:
perverse normalization and braid positivity. Here by perverse normalization we mean that the functor
of direct image to $\N$ is $t$-exact where the target category is equipped with the perverse coherent $t$-structure
of middle perversity. The content of the braid positivity condition is that the 
Coxeter generators of the group $B_{aff}^2$ acting, respectively, by the left convolution with $\F_\alpha=\Phi_{I^0I^0}^{-1}(\nabla_{s_\alpha})$
and right convolution with $\F_\alpha'=\Phi_{II}^{-1}(j_{s_\alpha *})$ are acting by right exact functors.

The arguments of \cite{AriB} show existence of a $t$-structure on $D^b(Coh^\LG(St'))$ whose heart is characterized by property (2) in part (b).
Perverse normalization holds for that $t$-structure since $\EE$ contains
the structure sheaf $\O$ as a direct summand. The braid positivity condition follows once we check the same property for the $t$-structure 
given by the tilting bundle $\EE \boxtimes \EE|_{\Nt}^*$ on $\gt \times \Nt$.

In \cite[\S 1.8]{BM} one finds a description of a collection of $t$-structures on $D^b(Coh(\gt))$, $D^b(Coh(\Nt))$
indexed by alcoves. The vector bundle $\EE$ (respectively, $\EE|_{\Nt}$) is a tilting generator for the heart of the $t$-structure
$\tau_{\aleph_0}$ corresponding to the fundamental alcove $\aleph_0$. It satisfies the braid positivity condition
for the $B_{aff}$ action as in \cite[Theorem 1.1.1]{BM},  we will call it the canonical action. One can check that
 $\F_\alpha'$ is identified with the kernels of the generators of the action on  $D^b(Coh(\Nt))$
 while $\F_\alpha$ is identified with the completion of the kernel of the action on $D^b(Coh(\gt))$ 
 at the preimage of $\N$,  
 see \cite[\S 5.5]{BeLo}, thus our equivalences are compatible with the $B_{aff}$ actions.

 It follows by \cite[Theorem 1.8.2(e)]{BM}
that  the dual bundle $\EE^*$ (respectively, $\EE^*|_{\Nt}$)  is a tilting generator for the heart of the $t$-structure
$\tau_{-\aleph_0}$.
 That $t$-structure satisfies the positivity condition for another  $B_{aff}$ action
conjugate to the canonical one by the element $\tilde w_0^{-1}$ \cite[Theorem 1.8.2(a)(2)]{BM}; 
we will refer to this as the {\em transposed} action of $B_{aff}$. 

Let $\sigma$ be the involution of $St$ or $\Nt\Ltimes_\Lg \Nt$ 
given by $\sigma(\b_1,\b_2,x)=(\b_2,\b_1,x)$. Then
we claim that the 
generators of the transposed action are identified with $\sigma^*(\F_\alpha)$, 
 $\sigma^*(\F_\alpha')$;
 since
right convolution with a complex $\F$  coincides with  left convolution with $\sigma^*(\F)$ this implies the claim.

For a finite simple root $\alpha$ we have $\sigma^*(\F_\alpha)\cong \F_\alpha$ and same for $\F_\alpha '$, this
follows from the explicit  description of $\F_\alpha$ in \cite[Theorem 1.3.2(a,i)]{BM}. On the
other hand, conjugation with $\tilde w_0$ 
permutes elements $\tilde s_\alpha$ where $\alpha$ runs over the set of finite simple roots. 

Consider now an affine simple root $\alpha_0$. We have $\tilde s_{\alpha_0} = \theta_{\beta} \tilde s$
where $s\in W_f$ is the reflection at the highest root $\beta$. Notice that $\tilde w_0$ commutes with
$\tilde s$, while $\tilde w_0 \tilde s^{-1}$ commutes with $\tilde s_{\alpha_0}$, this follows
from the corresponding obvious identities in $W$ since $\ell( w_0  s)=\ell(sw_0)=
\ell( w_0)-\ell(  s)= \ell( w_0  s s_{\alpha_0})-1.$
Thus $$\tilde w_0^{-1} \tilde s_{\alpha_0} \tilde w_0=\tilde s  \tilde s_{\alpha_0} \tilde s^{-1}=
\tilde s \theta_{\beta}.$$
Clearly $\sigma$ induces an anti-involution on the image of $B_{aff}$ which fixes each of the generators
 $\tilde s_\alpha$ for a finite $\alpha$ and  $\theta_\lambda$, $\lambda\in \Lambda$, thus the last  equality
 implies the claim.
 \qed

Similar properties hold for the rest of the equivalences \eqref{I0I0}--\eqref{PhiIWmon}.

\subsection{Lusztig's cells} In order to simplify the statement in this subsection we assume
that $G$ is simply-connected, thus $W$ is a Coxeter group. 
Recall the notion of a {\em two sided cell} in $W$. These are certain subsets
in $W$. In \cite{cells4} Lusztig has established a bijection between 2-sided
cells in $W$ and the set $\N/\LG$ of nilpotent conjugacy classes in $\Lg$. The set of two sided cells
is equipped with a partial order.
It has been conjectured by Lusztig and proved in \cite{B} that this order
matches the adjacency order on the set of nilpotent orbits under the bijection between two-sided cells and $\N/\LG$.  We now present a stronger statement relating the 2-sided cells to 
the support filtration introduced in \S \ref{CatPS}.


\begin{Thm}\label{cellsThm} Let $\unc$ be a two sided cell in $W$ and 
$O_{\unc}\subset \N$ be the corresponding nilpotent orbit. 

Let  $D_{\leq \unc}\subset  D_{I^0I}$ 
be the thick subcategory generated
by irreducible objects $IC_w\in \P$,  $w\in \unc'\leq \unc$.

Let $D^b(Coh^\LG_{O_\unc}(St'))$ be the full subcategory in $D^b(Coh^\LG(St'))$
consisting of complexes whose cohomology is set-theoretically supported
on the preimage of the closure of $O_{\unc}$. 

Then $D^b(Coh^\LG_{O_\unc}(St'))$ is the image of   $D_{\leq \unc}$  
under
the equivalence $\Psi_{I^0I}$.

\end{Thm} 

\proof Fix an orbit $O$ and let $W_O$ be the set of all $w\in W$ such that
$IC_w\in \Phi 
 ( D^b(Coh^\LG_{O}(St')) )$. 
Theorem \ref{suppf}(c) implies that $\Phi
( D^b(Coh^\LG_{O}(St')) )$
is generated as a triangulated category by $IC_w$, $w\in W_O$. 

Thus the Theorem follows once we check that $W_{O_\unc}=\cupl_{\unc'\leq \unc} \unc'$,
equivalently that for $w\in \unc$ the closure of the orbit 
$\overline{O_{\unc}}$ coincides with  the image in $\N$  of the support of $\Psi_{I^0I}
(IC_w)$.  Let $S_w$ denote that image.

We deduce this  from \cite[Theorem 4(a)]{B}, which provides the similar statement for the equivalence
\eqref{PhiIW} of Theorem \ref{mainthm2}. (In fact, {\em loc. cit.} deals with the category $^f\P$, a quotient category
of the category of Iwahori equivariant perverse sheaves on $\Fl$; however, $^f\P$ is equivalent
to the category in the left hand side of \eqref{PhiIW} by \cite[Theorem 2]{AB}).
Thus we see that 
\begin{equation}\label{cellIW}
\Phi^I_{IW}: D^b(Coh^\LG_{O_\unc}(\Nt))\iso (D^I_{IW})_{\leq \unc},
\end{equation}
where  $(D^I_{IW})_{\leq \unc}$ is the image of $D_{\leq \unc}$ under the functor the Whittaker
averaging functor $Av^{IW}:D_{I^0I}\to D^{IW}_I$ (the same notation was used above
for the Whittaker averaging functor $D_{I^0I^0}\to D^{IW}_{I^0}$).

The second commutative diagram in Proposition \ref{commupr} yields
\begin{equation}\label{compapr}
\begin{array}{ll}
(\Psi^I_{IW})^{-1}(IC_w^{IW})\cong p_{Spr,2*}'(\Psi_{I^0I}(IC_w))\ \ \ \ \ \ \ \ \  \ \ \ 
{\mathrm{for}}\ \ \ \ w\in W^f,\\
  p_{Spr,2*}'(\Psi_{I^0I}(IC_w))=0 \ \ \ \ \ \ \ \ \  \ \ \ 
{\mathrm{for}}\ \ \ \ w\not \in W^f;
\end{array}
\end{equation}
here $IC_w$, $w\in W$ and $IC_w^{IW}$, $w\in W^f$ are irreducible objects in 
$\P$ and $\P_I^{IW}$ respectively. 

It is clear that 
\begin{equation}\label{cell1}
S_w=\cupl_{\la,\mu} \supp (p_{Spr*} \Phi_{I^0I}^{-1}(IC_w)(\la,\mu)),
\end{equation}
\begin{equation}\label{cell2}
S_w\supset  \supp ( \Phi_{I^0I}^{-1}(IC_w)( M*IC_w*N)). 
\end{equation}
Since $D_{\leq \unc}$ is invariant under both left and right convolution, we  see that
\eqref{cell1} combined with  \eqref{cellIW}, \eqref{compapr} shows that $S_w\subseteq \overline
{O_\unc}$. Also, for any $w_1,\, w_2\in \unc$ the object $IC_{w_1}$ is a direct summand
in the convolution $X*IC_{w_2}*Y$ for some $X$, $Y$; thus \eqref{cell2} shows that
$S_w\supseteq \overline{O_\unc}$. \epf

\subsection{Exactness and Hodge $D$-modules}
Recall that in view of Corollary \ref{exactn}(a),
the restriction of the functors $\Psi_{I^0I}$, $\Psi_{I^0I^0}$ to the subcategory
of sheaves supported on the finite dimensional flag variety $G/B\subset \Fl$
is $t$-exact, i.e. it sends a perverse sheaf to a coherent sheaf.

On the other hand, a well known result in representation theory asserts that  the category 
$O$ for Langlands dual Lie algebras are equivalent, i.e. we have an equivalence
of abelian categories $$\Upsilon:Perv_{\LU}(\LG/\LB)\iso Perv_U(G/B)=
Perv_{I^0}(G/B).$$
This allows to state a relation between the restriction of our
equivalence $\Phi_{I^0I}$ to $Perv_{I^0}(G/B)\subset \P_{I^0 I}$ and
Hodge $D$-module theory.

Notice that the stack $St'/\LG$ can be interpreted as the {\em cotangent} to 
the stack $\LU\bs \LG/\LB$. Thus for a $\LU$-equivariant $D$-module $M$ on 
$\LG/\LB$ equipped with a $\LU$-equivariant good filtration we get $gr(M)\in
Coh^\LG(St')$. 

Let $\MH_\LU(\LG/\LB)$ be the category of mixed Hodge modules on $\LG/\LB$ 
equivariant with respect to $\LU$. 
We have the 
forgetful functor $Forg:\MH_\LU(\LG/\LB)\to D-mod_\LU(\LG/\LB)\cong Perv_\LU(\LG/\LB)$
where the second equivalence is the Riemann-Hilbert functor. 
Recall that a part of the data of a mixed Hodge structure on a $D$-module is a good
filtration, i.e. for $\tilde M\in \MH_\LU(\LG/\LB)$ the $D$-module $M=Forg(\tilde M)$ 
is equipped with a canonical good filtration. Thus we 
get a functor $gr:\MH_\LU(\LG/\LB)\to Coh^\LG(St')$.

\begin{Conj}
For $\tilde M\in \MH _\LU(\LG/\LB)$ we have a canonical isomorphism 
$$gr(\tilde M)\otimes \O(-\rho) \cong\Psi_{I^0 I}( \Upsilon(M)).$$
\end{Conj}

This Conjecture can be compared to the results of Ben-Zvi and Nadler \cite{BZN}.

\begin{Ex}
Recall that the finite Weyl group $W_f$ acts on the open subvariety $\gt^{reg}\subset \gt$.

For $w\in W_f$ let $\Ga_w\subset St$ be the closure of the graph of $w$.
Let $\Ga_w'$ be the scheme theoretic intersection $\Ga_w\cap St'$.
One can show that:
$$\Psi_{I^0I}:\Xis\mapsto \O_{St'},$$
$$\Psi_{I^0I}:j_{w*}\mapsto \O_{\Gamma_w'},$$ 
$$\Psi_{I^0I}:j_{w!} \mapsto \Omega_{\Ga_w'},$$
where $\Omega_{\Ga_w'}$ is the dualizing sheaf for the Cohen-Macaulay variety $\Ga_w'$
(the Cohen-Macaulay property is proven in \cite{BeRi}).
 Parallel results for associated graded of Hodge $D$-modules 
will be shown in \cite{BeRiHo}.
\end{Ex}

\medskip

We finish by sketching some generalizations of the equivalences described
in the paper. We expect they can be obtained by similar methods.

\subsection{Nonunipotent monodromy}
Consider the category of $\bI^2$ monodromic sheaves on $\Fltil$ with a fixed generalized
eigenvalues of monodromy. The latter corresponds to a tame rank one local system
on $\LT^2$, such local systems are in bijection with elements of $\LT^2$
(or a subset of that in the $l$-adic setting). For $\theta_1, \, \theta_2 
\in \LT$ let $D_{\theta_1,\theta_2}$ be the category of monodromic sheaves on 
$\Fltil$ with corresponding generalized eigenvalues of monodromy.

Let $\tii \LG\subset \LG\times \LG/\LB$ be the closed subvariety given
by $\tii \LG=\{(g,x)\ |\ g(x)=x\}$. We have a projection $\tii \LG \to \LT$.
Set $St_{grp}=\tii \LG\times_\LG\tii \LG$, and for
  $t_1,\, t_2\in \LT$ let $St_{grp}^{t_1,t_2}$ be the preimage of $(t_1,t_2)$ under the
  projection $St_{grp}\to \LT\times \LT$.

\begin{Conj}
We have a canonical equivalence of triangulated categories:
$$D_{\theta_1,\theta_2}\cong D^b\left( Coh^\LG_{St_{grp}^{t_1,t_2}}
(St_{grp}) \right).$$
\end{Conj}

An equivariant isomorphism between the variety of unipotent elements in $\LG$ and $\N$ 
and its extension to the formal neighborhoods in $\LG$ (respectively, $\Lg$) can be used
to identify the category $ Coh^\LG_{St_{grp}^{1,1}}(St_{grp})$ with $Coh^\LG_\N(St)$, thus in the 
special case $t_1=t_2=1$ the Conjecture amounts to equivalence
\eqref{I0I0}. One can also
state similar generalizations of \eqref{I0I}, \eqref{II}. Appearance of a group rather than a Lie algebra element here agrees with Langlands duality  where the element
is interpreted as the image of a topological generator of the tame ramification subquotient 
of the Galois group of the local field. On the other hand, working with the Lie algebra as we did
in the present article, makes it easier to describe a graded version of the category, Koszul duality etc. 

\subsection{Parabolic-Whittaker categories}
Let $P$ be a parabolic subgroup in $G$, and let $\bI_P\subset \GO$ be the parahoric subgroup
which is the preimage of $P$ under the projection $\GO\to G$. Let $\Fl_P=\GK/\bI_P$ 
be the corresponding partial affine flag variety. 

Let $Q$ be another parabolic subgroup and let $\psi_Q$ be an additive character of $\bI^0$
vanishing on the finite simple roots which are not in the Levi subgroup of $Q$ as well as on the
affine roots and not vanishing on the simple roots in the Levi of $Q$.
Let $D_{IW_Q}(\Fl_P)$ be the corresponding category of partial Whittaker sheaves.

Let $\LQ$, $\LP$ be the corresponding parabolic subgroups in $\LG$. Define $\Nt_\LQ
\subset \Nt_\LQ' \subset \gt_\LQ
\subset \LG/\LQ\times \Lg$ and $\Nt_\LP\subset \LG/\LP\times \Lg$ by:
$\gt_\LQ=\{(\fq,x)\ |\  x\in \fq\}$, $\Nt_\LP=\{(\fp,x)\ |\ x\in rad(\fp)\}$,
$\Nt_\LQ'=\{(\fq,x)\ |\ x\in rad(\fq) + {\mathfrak z}( \fq/rad(\fq))\}$, where we used the identification
between $\LG/\LQ$, respectively $\LG/\LP$ and the corresponding conjugacy class
of parabolic subalgebras,  $rad$ stands for the nilpotent radical and $\mathfrak z$ denotes the center.

\begin{Conj}\label{conjPQ}
We have canonical equivalences
\begin{equation}\label{conj1}
D_{IW_Q}(\Fl_P)\cong D^b(Coh^\LG(\gt_\LQ\times _\Lg \Nt_\LP)),
\end{equation}
\begin{equation}\label{conj2}
D_{\bI_Q}(\Fl_P)\cong DGCoh^\LG(\Nt_\LQ\Ltimes _\Lg \Nt_\LP)),
\end{equation}
\begin{equation}\label{conj3}
D_{\bI_Q'}(\Fl_P)\cong DGCoh^\LG(\Nt'_\LQ\Ltimes _\Lg \Nt_\LP)),
\end{equation}
where $\bI_Q'$ is the derived subgroup of $\bI_Q$.
\end{Conj}

There are natural pull-back, push-forward and Iwahori-Whittaker averaging functors between
the categories of constructible sheaves which should correspond to the functors
between the derived categories of coherent sheaves given by the natural correspondences,
 Proposition \ref{commupr} is an example of such a compatibility. 

\begin{Ex} Some special cases of Conjecture \ref{conjPQ} follow from results found in the literature. 

Let $P=Q=G$. Then the right hand side of \eqref{conj1} is $D^b(Coh^\LG(\{0\}))=
D^b(Rep(\LG))$. In this case \eqref{conj1} is essentially equivalent to the so-called
geometric Casselman-Shalika formula established in \cite{FGV}.

The right hand side of  \eqref{conj2} for $P=Q=G$ is $DGCoh^\LG(\{0\}\Ltimes_{\Lg} \{0\})$.
In view of Koszul duality, this special case of \eqref{conj2} follows from the second equivalence in  \cite[Theorem 5]{BeF}, it is discussed in detail (along with equivalences for various 
Ind-completions of the two categories) in \cite[\S 12]{AG1}.

For $Q=B$ and $P=G$ the left hand side in \eqref{conj3} is the derived category of $\bI^0$
equivariant sheaves on the affine Grassmannian, while the right hand is $DGCoh^\LG(\gt\Ltimes_{\Lg}\{0\})$. This special case of \eqref{conj3} 
amounts to one of the main results of  \cite{ABG}; again one needs to apply linear Koszul
duality \cite{MiR} to pass from the coherent side of the equivalence in  {\em loc. cit.} to the right
hand side of \eqref{conj3}; see also  \cite[\S 2.4]{BLa}.
\end{Ex}

Finally, let us mention the {\em Koszul duality} functors which give equivalences between
the graded version of $D_{IW_Q}(\Fl_P)$ and $D_{IW_P}(\Fl_Q)$, see \cite{BY}.
Under the first equivalence of Conjecture \ref{conjPQ} these should correspond
to {\em linear Koszul duality}  \cite{MiR}. In the special case when $P=Q=B$ is a Borel subgroup,
 this would provide a categorification
of the main result of \cite{MiR1}.


\end{document}